\documentclass{article}
\usepackage{CJK,CJKnumb,CJKulem,times,dsfont,ifthen,mathrsfs,latexsym,amsfonts,color}
\usepackage{amsmath,amsthm,makeidx,fontenc,amssymb,bm,graphicx,psfrag,listings,curves,extarrows,enumitem}
\usepackage{hyperref}
\usepackage[title]{appendix}
\usepackage{leftidx}

\allowdisplaybreaks[4] 

\usepackage{geometry}
\usepackage{indentfirst}
\usepackage{url}
\makeatletter
\def\url@leostyle{%
  \@ifundefined{selectfont}{\def\UrlFont{\sf}}{\def\UrlFont{\small\ttfamily}}}
\makeatother
\usepackage{amssymb}
\makeatletter

\newcommand{\Rmnum}[1]{\expandafter\@slowromancap\romannumeral #1@}
\makeatother

\newtheorem{theorem}{Theorem}[section]   
\newtheorem{definition}[theorem]{Definition} 	 
\newtheorem{lemma}[theorem]{Lemma}	 
\newtheorem{corollary}[theorem]{Corollary}	 
\newtheorem{proposition}[theorem]{Proposition}	 
\newtheorem{remark}{Remark}[section]

\newcommand{\al}{\alpha}
\newcommand{\ga}{\gamma}
\newcommand{\Ga}{\Gamma}

\newcommand{\e}{\varepsilon}

\newcommand{\iy}{\infty}
\newcommand{\q}{\theta}

\newcommand{\la}{\lambda}
\newcommand{\vp}{\varphi}
\newcommand{\pa}{\partial}

\newcommand{\lab}{\label}
\newcommand{\f}{\frac}
\newcommand{\bt}{\begin{theorem}}
\newcommand{\et}{\end{theorem}}
\newcommand{\bl}{\begin{lemma}}
\newcommand{\el}{\end{lemma}}
\newcommand{\bpr}{\begin{proposition}}
\newcommand{\epr}{\end{proposition}}
\newcommand{\bd}{\begin{definition}}
\newcommand{\ed}{\end{definition}}
\newcommand{\bc}{\begin{corollary}}
\newcommand{\ec}{\end{corollary}}
\newcommand{\bp}{\begin{proof}}
\newcommand{\ep}{\end{proof}}
\newcommand{\bx}{\begin{example}}
\newcommand{\ex}{\end{example}}
\newcommand{\bi}{\begin{exercise}}
\newcommand{\ei}{\end{exercise}}
\newcommand{\br}{\begin{remark}}
\newcommand{\er}{\end{remark}}
\newcommand{\be}{\begin{equation}}
\newcommand{\ee}{\end{equation}}
\newcommand{\bal}{\begin{align}}
\newcommand{\bn}{\begin{enumerate}}
\newcommand{\en}{\end{enumerate}}
\newcommand{\ba}{\begin{align}}
\newcommand{\ea}{\begin{align}}
\newcommand{\bg}{\begin{align*}}
\newcommand{\eg}{\end{align*}}
\newcommand{\bcs}{\begin{cases}}
\newcommand{\ecs}{\end{cases}}

\newcommand{\CR}{{\cal C}}

\newcommand{\LR}{{\cal L}}

\newcommand{\SR}{{\cal S}}

\newcommand{\R}{{\mathbb R}}
\newcommand{\bean}{\begin{eqnarray*}}
\newcommand{\eean}{\end{eqnarray*}}
\newcommand{\loc}{\operatorname{\rm loc}}

\newcommand{\sbr}[1]{\left(#1\right)}
\newcommand{\mbr}[1]{\left[#1\right]}
\newcommand{\lbr}[1]{\left\{#1\right\}}
\newcommand{\abr}[1]{\left\langle#1\right\rangle}
\newcommand{\rd}{\mathrm d}
\newcommand{\abs}[1]{\left\lvert#1\right\rvert}

\newcommand{\wt}{\widetilde}
\newcommand{\nm}[1]{\left\|#1\right\|}

\setitemize{itemindent=38pt,leftmargin=0pt,itemsep=-0.4ex,listparindent=26pt,partopsep=0pt,parsep=0.5ex,topsep=-0.25ex}

\numberwithin{equation}{section}

\begin{document}
\theoremstyle{plain}

\title{\bf Sharp estimates, uniqueness and nondegeneracy of positive solutions of the Lane-Emden system in planar domains }

\date{}
\author{ {\bf Zhijie Chen,\quad Houwang Li\quad \&\quad Wenming Zou} }

\maketitle
\begin{center}
\begin{minipage}{120mm}
\begin{center}{\bf Abstract }\end{center}

	We study the Lane-Emden system
	$$\bcs
		-\Delta u=v^p,\quad u>0,\quad\text{in}~\Omega,\\
		-\Delta v=u^q,\quad v>0,\quad\text{in}~\Omega,\\
		u=v=0,\quad\text{on}~\pa\Omega,
	\ecs$$
	where $\Omega\subset\mathbb{R}^2$ is a smooth bounded domain. In a recent work, we studied the concentration phenomena of positive solutions as $p,q\to+\infty$ and $|q-p|\leq \Lambda$.
	In this paper, we obtain sharp estimates of such multi-bubble solutions, including sharp convergence rates of local maxima and scaling parameters, and accurate approximations of solutions. As an application of these sharp estimates, we show that when $\Omega$ is convex, then the solution of this system is unique and nondegenerate for large $p, q$.



\vskip0.23in

\end{minipage}
\end{center}

\vskip0.4in
\section{Introduction}
	In this paper, we study the Lane-Emden system
\be\lab{equ1}
	\bcs
		-\Delta u=v^p,\quad u>0,\quad &\text{in}~\Omega,\\
		-\Delta v=u^q,\quad v>0,\quad &\text{in}~\Omega,\\
		u=v=0,&\text{on}~\pa\Omega,
	\ecs
\ee
where $\Omega$ is a smooth bounded domain in $\mathbb{R}^2$ and $p,q>0$. Without loss of generality, we may assume $q\geq p$.

One of our motivations to study the Lane-Emden system \eqref{equ1} is that \eqref{equ1} is one of the simplest Hamiltonian type elliptic systems. Due to the wide applications in physics and astrophysics, Hamiltonian-type elliptic systems have received great interests for decades, and abundant results on the existence and multiplicity of solutions are available in the literature, see e.g. the survey paper \cite{survey-1} for Hamiltonian elliptic systems and the survey paper \cite{survey-2} for especially the Lane-Emden system. On the other hand, the uniqueness of solutions is a central but challenging problem in the field of elliptic PDEs.
For the Lane-Emden system, it is only known (see \cite{unique-1}) that \eqref{equ1} has a unique positive solution in sublinear cases, i.e., $pq<1$. While for superlinear cases $pq>1$, it was proved in \cite{unique-2} that \eqref{equ1} has a unique least energy solution for $\Omega=\lbr{x\in\R^2:|x-x_0|< r}$ being a ball. To the best of our knowledge, the uniqueness of solutions for the Lane-Emden system \eqref{equ1} is still largely open.  In this paper, we will use blow up techniques to prove the uniqueness of positive solutions on convex domains for large $p,q>1$.

Another motivation of studying the Lane-Emden system \eqref{equ1} comes from the Lane-Emden equation
\be\lab{equ2}
	\bcs
	-\Delta u=u^p, \quad\text{in}~\Omega,\\
	u>0,\quad\text{in}~\Omega\\
	u=0,\quad\text{on}~\pa\Omega,
	\ecs
\ee
which arises from the astrophysics (see \cite{bg-1}) and has been widely studied in the literature. One attractive problem is the concentration phenomena of solutions when $p\to\iy$. In this direction, plenty of outstanding results have been given in the last decade, see a series of papers \cite{asy2-1,asy2-2,asy2-3,asy2-4,asy2-5,asy2-6,LE-1} and the references therein.
In particular, as an application of this concentration phenomena, it was proved in \cite{asy2-5,LE-1} that when $\Omega\subset\R^2$ is convex, the solution of the Lane-Emden equation \eqref{equ2} is unique and nondegenerate for large $p$, which partially answers a long standing conjecture.

Recently in \cite{Chen-Li-Zou}, we succeeded to construct the concentration phenomena of the Lane-Emden system \eqref{equ1}, which says that as $p,q\to\iy$ and $|q-p|\leq \Lambda$, the positive solutions behave as isolated bubbles around multiple points in $\Omega$. However, to obtain further the uniqueness, nondegeneracy and Morse index for the Lane-Emden system, we have to obtain sharp estimates of the multi-bubble solutions, which is one of the main purposes of this paper. The main results of this paper can be seen as generalizations of those in \cite{asy2-5,LE-1} for the Lane-Emden equation \eqref{equ2} to the Lane-Emden system \eqref{equ1}.

\subsection{Sharp estimates}

For any sequence $p_n\to\iy$, we take
\be\lab{assume-q}
	q_n=p_n+\q_n,\quad \q_n\ge0\quad \text{and}\quad \q_n\to\q~~\text{as}~n\to\iy.
\ee
Let $(u_n,v_n)$ be a solution sequence of \eqref{equ1} with $(p,q)=(p_n,q_n)$, which also satisfies the energy condition
\be\lab{assume-energy}
	\limsup_{n\to\iy} p_n\int_\Omega \nabla u_n\cdot \nabla v_n\rd x< \iy.
\ee
Remark that very recently Kamburov and Sirakov \cite{priori-1} proved that this energy condition \eqref{assume-energy} holds automatically for star-shaped domains; we will recall this important result in Lemma \ref{tem-100}.
Then the results in \cite{Chen-Li-Zou} tell us that there exist $k\in\mathbb N\setminus\{0\}$ and a finite set $\SR=\lbr{x_{\iy,1},\cdots,x_{\iy,k}}\subset\Omega$ such that up to a subsequence, the solution sequence $(u_n,v_n)$ will behave as bubble solutions concentrating around the set $\SR$. To be precise, for small $r>0$ and any $i=1,\cdots,k$, let $x_{n,i}$ be the local maximum point of $v_{n,i}$ defined by
\be
	v_n(x_{n,i})=\max_{B_{2r}(x_{\iy,i})}v_n,
\ee
then as $n\to\iy$ there holds
	$$x_{n,i}\to x_{\iy,i},\quad v_n(x_{n,i})\to\sqrt e, \quad u_n(x_{n,i})\to\sqrt e,$$
and
	$$ u_n=O(\f{1}{p_n}),~v_n=O(\f{1}{p_n}),\quad\text{in}~\CR_{loc}^2(\overline{\Omega}\setminus\SR).$$
The solution sequence $(u_n,v_n)$ are called multi-bubble solutions or $k$-bubble solutions ($1$-bubble solutions for $k=1$),  and the points $x_{\iy,1},\cdots,x_{\iy,k}$ are called blow up points (of $p_nu_n$ and  $p_nv_n$).

We introduce the Green function $G(x,y)$ of $-\Delta$ in $\Omega$ with the Dirichlet boundary condition:
$$\bcs
	-\Delta_x G(x,y)=\delta_y 	\quad \text{in}~\Omega,\\
	G(x,y)=0					\quad \text{on}~\pa\Omega,
\ecs$$
where $\delta_y$ is the Dirac function. It has the following form
\be\label{c1-g}  G(x,y)=-\f{1}{2\pi}\ln|x-y|-H(x,y),\quad(x,y)\in\Omega\times\Omega,  \ee
where the function $H(x,y)$ is the regular part of $G(x,y)$. It is well known that $H$ is a smooth function in $\Omega\times\Omega$, both $G$ and $H$ are symmetric in $x$ and $y$, and there is some constant $C>0$ such that
	$$|G(x,y)|\leq C(1+|\ln|x-y||),\qquad\forall \, x,y\in \Omega,~x\neq y,	$$
	$$|\nabla_xG(x,y)|\leq \f{C}{|x-y|},\qquad\forall \, x,y\in \Omega,~x\neq y.	$$
We recall the Robin function defined as
\be R(x)=H(x,x).\ee
For $\boldsymbol{x}=(x_1,\cdots,x_k)$ with $x_i\in\Omega$, we define the Kirchoff-Routh function $\Phi_k(\boldsymbol{x}):\Omega^k\to\R$ as
\be\lab{function-Phik}
	\Phi_k(\boldsymbol{x}):=\sum_{i=1}^k \Phi_{k,i}(\boldsymbol{x})  \quad\text{with}\quad
		\Phi_{k,i}(\boldsymbol{x})=R(x_i)-\sum_{j=1,j\neq i}^k G(x_i,x_j).
\ee

Recall the $k$-bubble solution $(u_n,v_n)$ and the blow up points $x_{\iy,1},\cdots,x_{\iy,k}$ given above. Then the blow up points fulfill the system
\be\lab{tem-35}
	\nabla_{\boldsymbol{x}}\Phi_k(\boldsymbol{x}_\iy)=0,\quad\text{where}~\boldsymbol{x}_\iy=(x_{\iy,1},\cdots,x_{\iy,k})\in\Omega^k.
\ee
Define the scaling parameters
\be 	
	\mu_{n,i}=\sbr{p_n v_n^{p_n-1}(x_{n,i})}^{-\f{1}{2}}. 	
\ee
Since $v_n(x_{n,i})\to\sqrt e$, we know $\mu_{n,i}\to0$ as $n\to\iy$. Define
\be\lab{function-wz}
	\bcs
		w_{n,i}(x)=\f{p_n}{v_n(x_{n,i})}(u_n(x_{n,i}+\mu_{n,i}x)-v_n(x_{n,i})),\\
		z_{n,i}(x)=\f{p_n}{v_n(x_{n,i})}(v_n(x_{n,i}+\mu_{n,i}x)-v_n(x_{n,i})),
	\ecs
\ee
then $(w_{n,i},z_{n,i})\to (U-\f{\q}{2},U)$ in $\CR_{loc}^2(\R^2)$, which means $w_{n,i},z_{n,i}$ behave as bubbles, where $\q$ is the limit of $\q_n$ in \eqref{assume-q} and 
\be\lab{function-U} 	U(x)=-2\ln (1+\f{1}{8}|x|^2)	 	\ee
is a solution of the Liouville equation
\be\lab{Liouville} \bcs
		-\Delta u=e^u\quad\text{in}~\R^2,\\
		\int_{\R^2}e^u\rd x=8\pi.
\ecs \ee
We refer to Section 2.2 for all the details about the asymptotic results, in particular see Theorem \ref{thm0}.

\vskip 0.1in
Based on the above results, we have $v_n(x_{n,i})\to\sqrt e$, $u_n(x_{n,i})\to\sqrt e$, $\mu_{n,i}\to0$, $w_{n,i}\to U-\f{\q}{2}$ and $z_{n,i}\to U$, which describes a rough picture of the asymptotic behaviors of $u_n,v_n$. However, in order to obtain further properties of $k$-bubble solutions, the rough picture is far from enough. One of our main results is the following sharper estimates of the convergence rates:

\bt\lab{thm1}
	Let $k\ge1$ and $(u_n,v_n)$ be a sequence of $k$-bubble solutions to \eqref{equ1}, then it holds that
\begin{itemize}[fullwidth,itemindent=0em]
\item[(a)]Convergence rates of local maxima:
	\be\lab{rate-1} \bcs
		u_n(x_{n,i})=\sqrt e\mbr{ 1 -\f{\ln p_n}{p_n-1} +\f{1}{p_n}\sbr{ 4\pi\Phi_{k,i}(\boldsymbol{x}_n)+3\ln2+2-\f{1}{4}\q_n }+O(\f{1}{p_n^{2-\delta}}) },\\
		v_n(x_{n,i})=\sqrt e\mbr{ 1 -\f{\ln p_n}{p_n-1} +\f{1}{p_n}\sbr{ 4\pi\Phi_{k,i}(\boldsymbol{x}_n)+3\ln2+2+\f{1}{4}\q_n }+O(\f{1}{p_n^{2-\delta}}) },
	\ecs \ee
	where $\Phi_{h,i}$ is defined by \eqref{function-Phik}, $\boldsymbol{x_n}=(x_{n,1},\cdots,x_{n,k})$ and $\delta>0$ is any fixed small constant.
\item[(b)]Convergence rates of scaling parameters:
	\be\lab{rate-2}
		\mu_{n,i}=e^{-\f{p_n}{4}}\mbr{e^{-(2\pi\Phi_{k,i}(\boldsymbol{x}_n)+\f{3}{2}\ln2+\f{3}{4}+\f{1}{8}\q_n) }+O(\f{1}{p_n^{1-\delta}})}.
	\ee
\item[(c)]Convergence rates of solutions:
	\be\lab{rate-3} \bcs
		w_{n,i}=U-\q_n\ln v_n(x_{n,i})+\f{s_{n,i}^*}{p_n}+O(\f{1}{p_n^{2-\delta}}),\\
		z_{n,i}=U+\f{t_{n,i}^*}{p_n}+O(\f{1}{p_n^{2-\delta}}),
	\ecs \ee
	in $\CR_{loc}^2(\R^2)$, where $s_{n,i}^*,t_{n,i}^*$ are radially symmetric functions given in \eqref{function-st*}.
\end{itemize}
\et

Note that when $\q_n=0$, i.e. $p_n=q_n$ in \eqref{equ1}, then
$$\int_\Omega|\nabla (u_{p_n}-v_{p_n})|^2=\int_\Omega (v_{p_n}^{p_n}-u_{p_n}^{p_n})(u_{p_n}-v_{p_n})\le 0,$$
which means $u_{p_n}=v_{p_n}$, so the Lane-Emden system \eqref{equ1} reduces to the (scalar) Lane-Emden equation \eqref{equ2}.
An analogue of Theorem \ref{thm1} for the Lane-Emden equation \eqref{equ2} was proved by Grossi, Ianni, Luo and Yan \cite{LE-1}, and Theorem \ref{thm1} is a generalization of their result to the system case. Remark that this type of theorems is much harder for the Lane-Emden system than for the scalar equation. Comparing with the scalar equation, the Lane-Emden system brings more difficulties mainly in two aspects. One is that we cannot estimate directly the errors between $(w_{n,i},z_{n,i})$ and its limit $(U-\f{\q}{2},U)$, and cannot estimate directly the errors between $(s_{n,i},t_{n,i})$ and its limit $(s_{\iy,i},t_{\iy,i})$ either (see Remark \ref{rmk-42}), where $(s_{n,i},t_{n,i})$ is defined by \eqref{function-stn}. We fix this problem by approximating the bubbles using the parameterized limit functions. Another difficulty is due to the difference of linearized equations: for the Lane-Emden equation, the solution space of the corresonding linearized equation is $3$ dimensional, see \cite[Lemma 2.3]{Lin=CPAM=2002}; but for the Lane-Emden system, we will prove in Lemma \ref{linear} that the solution space of the corresponding linearized system is $4$ dimensional. To overcome this difficulty, we have to estimate the sharp convergence rates of $v_n(x_{n,i})-u_n(x_{n,i})$ first, which will be done in Section 3.

We believe that this sharp estimates can be used to study properties of $k$-bubble solutions of the Lane-Emden system, such as local uniqueness, nondegeneracy, Morse index, etc. We also think that it can also bring some inspirations for constructing $k$-bubble solutions using reduction methods.

\subsection{Uniqueness and nondegeneracy}
By using the sharp estimates in Theorem \ref{thm1}, we prove the uniqueness and nondegeneracy for $1$-bubble solutions, i.e. $k=1$. Notice that $\Phi_1=R$, namely the Kirchoff-Routh function reduces to the Robin function in this case.

Let $(u_n,v_n)$ be a sequence of $1$-bubble solutions of system \eqref{equ1}, whose blow up point is $x_\iy\in\Omega$. Then \eqref{tem-35} implies that $x_\iy$ is a critical point of the Robin function $R$. We have the following uniqueness result.

\bt\lab{thm3}
	Let $(u_n^{(1)},v_n^{(1)})$ and $(u_n^{(2)},v_n^{(2)})$ be two sequences of $1$-bubble solutions of the Lane-Emden system \eqref{equ1} with the same $(p_n, q_n)$ and concentrate at the same point $x_\iy$. If $x_\iy$ is a nondegenerate critical point of $R(x)$, then there exists a $n_0>0$ such that
		$$u_n^{(1)}=u_n^{(2)},\quad v_n^{(1)}=v_n^{(2)},\quad\text{for any}~n\ge n_0.$$
\et

Recall the linearized system of \eqref{equ1} at $(u_n,v_n)$
\be\lab{linear-0} \bcs
	-\Delta \xi=p_nv_n^{p_n-1}\eta,&\text{in}~\Omega,\\
	-\Delta \eta=q_nu_n^{q_n-1}\xi,&\text{in}~\Omega,\\
	\xi=\eta=0,&\text{on}~\partial\Omega.
\ecs \ee
The solution $(u_n,v_n)$ is called {\it nondegenerate}, if the linearized system \eqref{linear-0} has only the trivial solution $(0,0)$. We have the following nondegenerate result.

\bt\lab{thm2}
	Let $(u_n,v_n)$ be a sequence of $1$-bubble solutions concentrating at $x_\iy$. If $x_\iy$ is a nondegenerate critical point of $R(x)$, then there exists a $n_0>0$ such that $(u_n,v_n)$ is nondegenerate for any $n\ge n_0$.
\et

We emphasize that Theorems \ref{thm3}-\ref{thm2} hold for any smooth bounded domains under the assumption of the nondegeneracy of the critical points of the Robin function. The exact number and nondegeneracy of the critical points of the Robin function is a widely studied but challenging problem, see \cite{cp-1,cp-2,cp-3,cp-4,cp-5} and the references therein.

In the particular case when $\Omega$ is convex, there can only exist $1$-bubble solutions for system \eqref{equ1}. Indeed, suppose $(u_n,v_n)$ is a sequence of $k$-bubble solutions of system \eqref{equ1} with $k\ge 1$ different blow up points $x_{\iy,1},\cdots,x_{\iy,k}$. Then it holds $$\nabla_{\boldsymbol{x}}\Phi_k(\boldsymbol{x}_\iy)=0,\quad\text{where}~\boldsymbol{x}_\iy=(x_{\iy,1},\cdots,x_{\iy,k}).$$
But the result in \cite{cp-1} says that in a convex domain, $\nabla_{\boldsymbol{x}}\Phi_k(\boldsymbol{x})\neq0$ unless
	$$\boldsymbol{x}\in\lbr{(x_{1},\cdots,x_{k})\in\Omega^k:~x_i=x_j~\text{for some}~i\neq j},$$
which implies $k=1$. Moreover, it is shown that in a convex domain, the Robin function has only one critical point, which is also nondegenerate (\cite{cp-3}). Besides, it follows from \cite{priori-1} that the energy condition \eqref{assume-energy} holds automatically for convex domains.  Hence, combining these conclusions with Theorem \ref{thm1} and \ref{thm2}, we obtain the following uniqueness result.
\bt\lab{thm4}
	Let $\Omega\subset\R^2$ be a smooth bounded convex domain. Then for any $\Lambda>0$, there exists a $p_0=p_0(\Omega,\Lambda)>0$ such that for any $p,q\ge p_0$ satisfying $|p-q|\le \Lambda$, the Lane-Emden system
	\be\lab{tem-95} \bcs
		-\Delta u=v^p,\quad u>0,\quad\text{in}~\Omega,\\
		-\Delta v=u^{q},\quad v>0,\quad\text{in}~\Omega,\\
		u=v=0,\quad\text{on}~\pa\Omega,
	\ecs \ee
	has a unique solution, which is also nondegenerate.
\et


\vskip 0.2in
This paper is organized as follows. In Section 2, we characterize the solutions of the linearized system which is applied frequently in the sharp estimates, and we also give some results used in the following sections. In Section 3, we give a sharp estimate of $v_n(x_{n,i})-u_n(x_{n,i})$, which plays an essential role in the proof of Theorem \ref{thm1}. In Section 4, we construct the sharp estimates of $w_{n,i},z_{n,i}$ and give the sharp estimates of $\mu_{n,i}$ and $v_n(x_{n,i})$, and then finish the proof of Theorem \ref{thm1}. In Section 5 and Section 6, we study the nondegeneracy and uniqueness of $1$-bubble solutions, where we use blow up techniques and local Pohozaev identities. Finally, in Section 7, we consider the special case when $\Omega$ is convex and prove Theorem \ref{thm4}.

Throughout this paper, we use the following notations.
\begin{itemize}

\item We use $C, C_1, \cdots$ to denote various positive constants that is independent of $(p_n,q_n)$.
Denote the ball $B_r(x):=\{y\in\mathbb{R}^2 | |y-x|<r\}$.
\item For a function $f(x,y)$, $x,y\in\R^2$, we use $\pa_i$ and $\nabla$ to denote the derivatives of $f(x,y)$ with respect to the first variable $x$, and use $D_i$ and $D$ to denote the derivatives of $f(x,y)$ with respect to the second variable $y$.

\item For a function $f(x)$, $x\in\R^2$, we also write $f(x)=f(|x|, x/|x|)$ for convenience. In particular, we write $f(x)=f(|x|)$ if $f(x)$ is radially symmetric.
\end{itemize}

\section{Preliminary}

\subsection{The linearized system}

In later sections we need to use many facts about the linearized system in our proofs. First we recall Chen and Lin's remarkable result from
their seminal work \cite{Lin=CPAM=2002}.

\bl\cite[Lemma 2.3]{Lin=CPAM=2002}\lab{linear-scalar} Recall $U(x)$ in \eqref{function-U} and let
	 $u(x)$ be a solution of the linearized equation
	\be\lab{lineareq-0}\bcs
		-\Delta u=e^U u &\text{in}~\R^2,\\
		|u(x)|\le C(1+|x|)^\tau &\text{in}~\R^2,
	\ecs\ee
	for some $\tau\in[0,1)$. Then
		$$u(x)=\sum_{i=0}^2c_i\phi_i(x),$$
	for some constants $c_i\in\R$, where
		\begin{equation}\label{eq1}\phi_0(x)=\f{8-|x|^2}{8+|x|^2}\quad\text{and}\quad \phi_i(x)=\f{x_i}{8+|x|^2} \quad\text{for}~i=1,2.\end{equation}
\el

Since we deal with the Lane-Emden system rather than the scalar equation, Lemma \ref{linear-scalar} is not enough for us. Here we need the following new result for the linearized system.
\bl\lab{linear}
	Recall $U(x)$ in \eqref{function-U} and let $(u,v)$ be a solution of the linearized system
	\be\lab{linearsys-0}\bcs
		-\Delta u=e^U v &\text{in}~\R^2,\\
		-\Delta v=e^U u &\text{in}~\R^2,\\
		|u(x)|,|v(x)|\le C(1+|x|)^\tau&\text{in}~\R^2,
	\ecs\ee
	for some $\tau\in[0,1)$. Then
		$$(u(x),v(x))=\sum_{i=0}^2c_i(\phi_i(x),\phi_i(x))+c_3(\phi_3(x),-\phi_3(x)),$$
	for some constants $c_i\in\R$, where $\phi_0, \phi_1, \phi_2$ are given by \eqref{eq1},
	and $\phi_3$ is a positive radial function satisfying $\phi_3(0)=1$ and \begin{equation}\label{eq2}\phi_3(x)=\f{e^{\f{\sqrt7}{2}\pi}+e^{-\f{\sqrt7}{2}\pi}}{2\pi}\ln(1+\f{1}{8}|x|^2)+O(1),
\quad\text{as $|x|\to\iy$},\end{equation}
\begin{equation}\label{eq3}\nabla\phi_3(x)=\f{e^{\f{\sqrt7}{2}\pi}+e^{-\f{\sqrt7}{2}\pi}}{\pi}\f{x}{|x|^2}+o(\f{1}{|x|}),
\quad\text{as $|x|\to\iy$}.\end{equation}
\el

We will see that Lemma \ref{linear} plays a crucial role in our proofs of Theorems \ref{thm1}-\ref{thm4}.
To prove Lemma \ref{linear}, we write
	$$u=\f{u+v}{2}+\f{u-v}{2} \quad\text{and}\quad v=\f{u+v}{2}-\f{u-v}{2}.$$
Then $\frac{u+v}{2}$ solves \eqref{lineareq-0}, so Lemma \ref{linear-scalar} implies
	$$\f{u(x)+v(x)}{2}=\sum_{i=0}^2c_i\phi_i(x).$$
for some constants $c_i$.
On the other hand, $\f{u-v}{2}$ satisfies
\be\lab{tem-1}\bcs
	-\Delta \vp+e^U\vp=0\quad\text{in}~\R^2,\\
	|\vp(x)|\le C(1+|x|)^\tau\quad\text{in}~\R^2.
\ecs\ee
Therefore, Lemma \ref{linear} follows immediately from the following result.

\bl\lab{linear-1}
	Let $\vp(x)$ be a solution of \eqref{tem-1}. Then $\vp(x)=\vp(0)\phi_3(x)$, where $\phi_3(x)$ is a positive radial function satisfying $\phi_3(0)=1$ and \eqref{eq2}-\eqref{eq3}.
\el

\bp Let $\vp(x)$ be a solution of \eqref{tem-1}. Take the angular momentum decomposition
\be\lab{tem-4} 	\vp(x)=\sum_{j=0}^\iy \sbr{ \vp_{j,1}(|x|)\sin(j\f{x}{|x|})+\vp_{j,2}(|x|)\cos(j\f{x}{|x|}) }, 	\ee
where 
\be\lab{tem-5}	\vp_{j,1}(r)=\f{1}{2\pi}\int_0^{2\pi}\vp(r,\theta)\sin(j\theta)\rd \theta 		\ee
and
\be 	\vp_{j,2}(r)=\f{1}{2\pi}\int_0^{2\pi}\vp(r,\theta)\cos(j\theta)\rd \theta. 	\ee
Since $\sin{j\theta},\cos{j\theta}$ are the corresponding eigenfunctions of $-\Delta_{S^1}$ with respect to the eigenvalue $\la_j=j^2$ for nonnegative integers $j\ge 0$ (here $S^1=\{x\in\mathbb{R}^2\,|\, |x|=1\}$, see e.g. \cite{Berezin-Shubin}), we have that $\vp_{j,i}$, $i=1,2$, are solutions of the following ODE
\be\lab{tem-2} 	\bcs
	\phi''+\f{1}{r}\phi'=\f{1}{(1+\f{1}{8}r^2)^2}\phi+\f{j^2}{r^2}\phi,\quad r>0,\\
|\phi(r)|\leq C(1+r)^{\tau},\quad r\geq 0,
\ecs \ee
for any $j\ge 0$. Here we used $e^{U(x)}=\f{1}{(1+\f{1}{8}|x|^2)^2}$. Then
\begin{equation}\label{2-11}
(r\phi')'=\sbr{\f{j^2}{r}+\f{r}{(1+\f{1}{8}r^2)^2}}\phi.
\end{equation}

{\bf Step 1:} We prove that
$\vp_{j,i}\equiv0$ for any $j\ge1$ and $i=1,2$.

Suppose for some $j\ge1$ that $\phi\not\equiv0$ is a solution of \eqref{tem-2}. It is easy to see from \eqref{tem-2} that any local maxima of $\phi$ must be nonpositive and any local minima of $\phi$ must be nonnegative.

First, we claim that $\phi$ is a monotone function. If not, then by replacing $\phi$ with $-\phi$ if necessary, we may assume that $\phi$ has a local minimum point $r_0>0$. Then $\phi(r_0)\geq 0$ and $\phi(r)$ is increasing in $[r_0, +\infty)$. Take $r_1>r_0$ such that $\phi(r_1)>0$ and $\phi'(r_1)>0$. Then \eqref{2-11} gives $r\phi'(r)\geq r_1\phi'(r_1)>0$ for any $r\geq r_1$, which implies
\begin{equation}\label{2-12}\phi(r)\geq r_1\phi'(r_1) (\ln r-\ln r_1)+\phi(r_1)\to +\infty,\quad \text{as }r\to+\infty.\end{equation}
Furthermore, by multiplying \eqref{2-11} with $r\phi'$ and integrating from $r_1$ to $r$, we get
\be\lab{2-13} 	\sbr{r\phi'(r)}^2\ge j^2\sbr{\phi(r)^2-\phi(r_1)^2},\quad \forall r\geq r_1.	\ee
By $\tau \in [0,1)$ and \eqref{2-12} we can take $r_2>r_1$ such that $\phi(r)\ge \phi(r_1)\sbr{1-\f{(1+\tau)^2}{4j^2}}^{-\f{1}{2}}$ for any $r\ge r_2$. Then we see from \eqref{2-13} that $\sbr{r\phi'(r)}^2\ge \frac{(1+\tau)^2}{4}\phi(r)^2$ for $r\geq r_2$, i.e.
$$\sbr{\ln\phi(r)}'=\f{\phi'(r)}{\phi(r)}\ge \f{1+\tau}{2r},\quad\text{for}~r\ge r_2, $$
from which we conclude that
\be\label{2-14}\phi(r)\ge \phi(r_2)\sbr{\f{r}{r_2}}^{\f{1+\tau}{2}}, \quad\text{for}~r\ge r_2,\ee
clearly a contradiction with \eqref{tem-2} and $\tau\in [0,1)$.

Thus $\phi$ is a monotone function. Again by replacing $\phi$ with $-\phi$ if necessary, we may assume $\phi$ is increasing in $[0,+\infty)$. If $\phi(0)<0$, then \eqref{2-11} implies $\phi'(r)<0$ for $r>0$ small, a contradiction. So $\phi(0)\ge 0$, then the same argument as \eqref{2-12}-\eqref{2-14} yields a contradiction with \eqref{tem-2}.

This proves $\vp_{j,i}\equiv0$ for any $j\ge1$, $i=1,2$.
Then we see from \eqref{tem-4} that $\vp(x)=\vp(|x|)$ is a smooth radial function which satisfies
\be\lab{tem-6} 	\bcs
	\vp''+\f{1}{r}\vp'=\f{1}{(1+\f{1}{8}r^2)^2}\vp,\quad r>0,\\
|\vp(r)|\leq C(1+r)^{\tau},\quad r\geq 0,\\
	\vp'(0)=0.
\ecs \ee

{\bf Step 2: } We solve the equation \eqref{tem-6} by using the well-known Hypergeometric equation.

By letting $t=1+\f{1}{8}r^2$ and $\psi(t)=\vp(r)$, we get
\be\lab{tem-7} 	t^2(t-1)\frac{\rd^2\psi}{\rd t^2}+t^2\frac{\rd\psi}{\rd t}-2\psi=0,\quad t>1.  \ee
Recall the Hypergeometric equation (see e.g. \cite[Chapter 2]{GP} or \cite[Chapter 15]{2F1})
	$$z(1-z)\frac{\rd^2w}{\rd z^2}+\sbr{d_3-(d_1+d_2+1)z}\frac{\rd w}{\rd z}-d_1d_2w=0,\quad z\in\mathbb{C},$$
where $d_1,d_2,d_3$ are complex constants. It has a solution of the following form $$w(d_1,d_2;d_3;z)=(1-z)^{-d_1-d_2+d_3}z^{d_1-d_3}\cdot\leftidx{_2}F_1(-d_1+d_3,1-d_1;-d_1-d_2+d_3+1;1-\f{1}{z}),$$
where $\leftidx{_2}F_1(d_1,d_2;d_3;z)$ is the Hypergeometric function defined by
\begin{equation}\label{2-17}\leftidx{_2}F_1(d_1,d_2;d_3;z)
:=\sum_{n=0}^{+\infty}\frac{(d_1)_{n}(d_2)_{n}}{(d_3)_{n}(1)_n}z^n,\end{equation}
with
\[(d)_n:=\begin{cases}1&\text{if }n=0,\\
d(d+1)\cdots(d+n-1) &\text{if }n\in\mathbb{N}_{\geq 1}.\end{cases}\]
Note that $(d)_n$ can be expressed by the Gamma function $\Gamma(z)$ as follows
\[(d)_n=\frac{\Gamma(d+n)}{\Gamma(d)},\quad \text{if }d\neq 0, -1,-2,\cdots,\]
and in particular, $(1)_n=n!$.

Set $f_d(z):=z^{d}\cdot w(d,d;2d;z)=\leftidx{_2}F_1(d,1-d;1;1-\f{1}{z})$, then by direct computations we obtain
\be\lab{tem-8} 	
	z^2(z-1)\frac{\rd^2f_d}{\rd z^2}+z^2\frac{\rd f_d}{\rd z}+d(d-1)f_d=0,\quad z\in\mathbb{C}.
\ee
Comparing \eqref{tem-7} and \eqref{tem-8}, by taking $d=\f{1+i\sqrt 7}{2}$, we obtain a solution $f_d(z)=\leftidx{_2}F_1(d,1-d;1;1-\f{1}{z})$ of \eqref{tem-7} in the complex plane. For any $t\in\R$, it follows from \eqref{2-17} and $\bar d= 1-d$ that
	$$\overline {f_d(t)}=\leftidx{_2}F_1(\bar d,\overline{1-d};1;1-\f{1}{t})=\leftidx{_2}F_1(1-d,d;1;1-\f{1}{t})=f_d(t),$$
i.e. $f_d(t)\in\R$ as long as $t\in\R$, so $f_d(t)$ is a solution of \eqref{tem-7}. Thus we obtain a positive solution of \eqref{tem-6} as follows
	\begin{align*}\vp(r)
	&=f_d(1+\f{1}{8}r^2)=\leftidx{_2}F_1(d,1-d;1;\f{r^2}{8+r^2})\\
	&=\sum_{n=0}^{+\infty}\frac{(d)_{n}(1-d)_{n}}{(1)_{n}(1)_n}\left(\f{r^2}{8+r^2}\right)^n
		=\sum_{n=0}^{+\infty}\frac{(d)_{n}(\bar{d})_{n}}{(1)_{n}(1)_n}\left(\f{r^2}{8+r^2}\right)^n\\
	&=\sum_{n=0}^{+\infty}\frac{|(d)_{n}|^2}{(1)_{n}(1)_n}\left(\f{r^2}{8+r^2}\right)^n>0.
\end{align*}
Clearly this implies
	$$\vp(r)=1+\f{1}{4}r^2+O(r^4),\quad\text{as}~r\to0.$$
Next we give the asymptotic estimates of $\vp(r)$ and $\vp'(r)$ as $r\to\iy$. By using the formula (15.3.10) in \cite{2F1}, we get 
\be\lab{tem723-1}
	\begin{aligned}
	\leftidx{_2}F_1(d,1-d;1;z)=\f{1}{\Ga(d)\Ga(\bar d)}\sum_{n=0}^\iy \f{(d)_n(\bar d)_n}{(n!)^2}&(1-z)^n [2\psi^{(0)}(n+1)-\psi^{(0)}(n+d) \\
	&-\psi^{(0)}(n+\bar d)-\ln(1-z)  ]
	\end{aligned}
\ee
for $|1-z|<1$ and $|\arg(1-z)|<\pi$, where $\psi^{(0)}$ is the Digamma function (see e.g. \cite[Chapter 6]{2F1}). From the properties of $\psi^{(0)}$, we have
	$$\psi^{(0)}(1)=-\boldsymbol{\ga},\quad \psi^{(0)}(n+z)=\psi^{(0)}(z)+\sum_{j=0}^{n-1}\f{1}{j+z},$$
where $\boldsymbol{\ga}$ is the Euler-Mascheroni constant. Let 
	$$C_n:=\f{(d)_n(\bar d)_n}{(n!)^2}[2\psi^{(0)}(n+1)-\psi^{(0)}(n+d)-\psi^{(0)}(n+\bar d)].$$
Then for $n\ge2$, it holds
	$$\begin{aligned}
	\f{C_n}{C_{n-1}}
	&=\f{(d+n-1)(\bar d+n-1)}{n^2} \times \f{-2\boldsymbol{\ga}-\psi^{(0)}(d)-\psi^{(0)}(\bar d)+\sum_{j=0}^{n-1}\sbr{\f{2}{j+1}-\f{1}{j+d}-\f{1}{j+\bar d}}  }{-2\boldsymbol{\ga}-\psi^{(0)}(d)-\psi^{(0)}(\bar d)+\sum_{j=0}^{n-2}\sbr{\f{2}{j+1}-\f{1}{j+d}-\f{1}{j+\bar d}} } \\
	&=\sbr{1-\f{1}{n}+\f{2}{n^2}} \times \f{ 2\boldsymbol{\ga}+\psi^{(0)}(d)+\psi^{(0)}(\bar d)+\sum_{j=0}^{n-1}\f{j-3}{(j+1)(j^2+j+2)} }{ 2\boldsymbol{\ga}+\psi^{(0)}(d)+\psi^{(0)}(\bar d)+\sum_{j=0}^{n-2}\f{j-3}{(j+1)(j^2+j+2)}  }.
	\end{aligned}$$
which means $\lim_{n\to\iy}\f{C_n}{C_{n-1}}=1$. Thus the power series $\sum_{n=1}^\iy C_n z^{n-1}$ is convergent for $|z|<1$. So it follows from \eqref{tem723-1} that for $0<z<1$ and $z\to1$ it holds
	$$\begin{aligned}
		\leftidx{_2}F_1(d,1-d;1;z)
		&=\f{1}{\Ga(d)\Ga(\bar d)}\mbr{2\psi^{(0)}(1)-\psi^{(0)}(d)-\psi^{(0)}(\bar d)-\ln(1-z)}\\
			&\quad -\f{(1-z)\ln(1-z)}{\Ga(d)\Ga(\bar d)}\sum_{n=1}^\iy\f{(d)_n(\bar d)_n}{(n!)^2}(1-z)^{n-1}+\f{(1-z)}{\Ga(d)\Ga(\bar d)}\sum_{n=1}^\iy C_n(1-z)^{n-1}\\
		&=\f{-\ln(1-z)}{\Ga(d)\Ga(\bar d)}+O\sbr{1+(1-z)\ln(1-z)}.
	\end{aligned}$$
Noting that 
	$$\Ga(d)\Ga(\bar d)=\abs{\Ga\sbr{\f{1}{2}+i\f{\sqrt 7}{2}}}^2=\f{2\pi}{e^{\f{\sqrt7}{2}\pi}+e^{-\f{\sqrt7}{2}\pi}},$$
we get
	$$\begin{aligned}
		\vp(r)=\leftidx{_2}F_1(d,1-d;1;\f{r^2}{8+r^2})=\f{e^{\f{\sqrt7}{2}\pi}+e^{-\f{\sqrt7}{2}\pi}}{2\pi}\ln(1+\f{1}{8}r^2)+O(1),\quad\text{as}~r\to\iy.
	\end{aligned}$$
By using the formulas (15.2.1) and (15.3.12) in \cite{2F1}, we also get 
\be\lab{tem723-2}
	\begin{aligned}
	&\f{\rd}{\rd z}\leftidx{_2}F_1(d,1-d;1;z)
	=d\bar d \cdot\leftidx{_2}F_1(d+1,\bar d+1;2;z)\\
	&=\f{2}{\Ga(d+1)\Ga(\bar d+1)}(1-z)^{-1} - \f{2}{\Ga(d)\Ga(\bar d)}\sum_{n=0}^\iy \f{(d+1)_n(\bar d+1)_n}{n!(n+1)!}(1-z)^n \\
		&~~~~~~\times [\psi^{(0)}(n+1)+\psi^{(0)}(n+2)-\psi^{(0)}(n+d+1)-\psi^{(0)}(n+\bar d+1)-\ln(1-z) ]
	\end{aligned}
\ee
for $|1-z|<1$ and $|\arg(1-z)|<\pi$. Repeating the above analysis, we obtain that for $0<z<1$ and $z\to1$ it holds
	$$\f{\rd}{\rd z}\leftidx{_2}F_1(d,1-d;1;z)=\f{2}{\Ga(d+1)\Ga(\bar d+1)}(1-z)^{-1}+O(\ln(1-z)).$$
Since 
	$$\Ga(d+1)\Ga(\bar d+1)=\abs{\Ga(d+1)}^2=\abs{d\cdot\Ga(d)}^2=\f{4\pi}{e^{\f{\sqrt7}{2}\pi}+e^{-\f{\sqrt7}{2}\pi}},$$
we get
	$$\begin{aligned}
		\phi'(r)
		&=\f{\rd}{\rd z}\leftidx{_2}F_1(d,1-d;1;\f{r^2}{8+r^2})\cdot\f{16r}{(8+r^2)^2}\\
		&=\f{e^{\f{\sqrt7}{2}\pi}+e^{-\f{\sqrt7}{2}\pi}}{4\pi}\sbr{1+\f{r^2}{8}+O(\ln r)}\sbr{\f{16}{r^3}+O(\f{1}{r^5})}\\
		&=\f{e^{\f{\sqrt7}{2}\pi}+e^{-\f{\sqrt7}{2}\pi}}{\pi}\f{1}{r}+o(\f{1}{r}),\quad\text{as}~r\to\iy.
	\end{aligned}$$
So setting $\phi_3(x)=\vp(|x|)=f_d(1+\f{1}{8}|x|^2)>0$, we obtain that $\phi_3(x)$ is a smooth radial solution of \eqref{tem-1} satisfying $\phi_3(0)=1$ and as $|x|\to\iy$,
		$$\phi_3(x)=\f{e^{\f{\sqrt7}{2}\pi}+e^{-\f{\sqrt7}{2}\pi}}{2\pi}\ln(1+\f{1}{8}|x|^2)+O(1),$$
		$$\nabla\phi_3(x)=\f{e^{\f{\sqrt7}{2}\pi}+e^{-\f{\sqrt7}{2}\pi}}{\pi}\f{x}{|x|^2}+o(\f{1}{|x|}).$$
By the uniqueness of ODE equations, $\vp(x)=\vp(0)\phi_3(x)$. Then we finish the proof.
\ep

\subsection{Asymptotic behavior}
For any sequence $p_n\to\iy$, we take $q_n=p_n+\q_n$ with $\q_n\ge0$ and $\q_n\to\q$ as $n\to\iy$. Let $(u_n,v_n)$ be a solution sequence of \eqref{equ1} with $(p,q)=(p_n,q_n)$, which also satisfies the energy condition
\be
	\limsup_{n\to\iy} p_n\int_\Omega \nabla u_n\cdot \nabla v_n\rd x< \iy.
\ee
Then
\be
	p_n\int_\Omega u_n^{q_n+1}=p_n\int_\Omega v_n^{p_n+1}=p_n\int_\Omega\nabla u_n\cdot\nabla v_n\le C,
\ee
and so
\be\lab{est-v-1}
	p_n\int_\Omega v_n^{p_n}
	\le \sbr{p_n\int_\Omega v_n^{p_n+1}}^{\f{p_n}{p_n+1}}(p_n|\Omega|)^{\f{1}{p_n+1}}
	\le C,
\ee
\be\lab{est-u-1}
	p_n\int_\Omega u_n^{q_n}
	\le \sbr{p_n\int_\Omega u_n^{q_n+1}}^{\f{q_n}{q_n+1}}(p_n|\Omega|)^{\f{1}{q_n+1}}
	\le C
\ee
for $p_n$ large.

For the finite set $\SR=\lbr{x_{\iy,1},\cdots,x_{\iy,k}}\subset\Omega$ in Theorem \ref{thm0} below,
we take $d_0>0$ such that for all $i,j=1,\cdots,k$ and $i\neq j$ it holds
\be\lab{constant-d0}
	B_{2d_0}(x_{\iy,i})\subset\Omega,\quad B_{2d_0}(x_{\iy,i})\cap B_{2d_0}(x_{\iy,j})=\emptyset.
\ee
In \cite{Chen-Li-Zou}, we have proved that
\bt(\cite{Chen-Li-Zou})\lab{thm0}
There exist $k\in\mathbb N\setminus\{0\}$ and a finite set of concentration points $\SR=\lbr{x_{\iy,1},\cdots,x_{\iy,k}}\subset\Omega$ such that up to a subsequence, $(u_n,v_n)$ satisfies the following properties:
	for fixed $r\in(0,d_0]$ and $i=1,\cdots,k$, set
		$$v_n(x_{n,i})=\max_{B_{2r}(x_{\iy,i})}v_n\quad \text{and}\quad u_n(y_{n,i})=\max_{B_{2r}(x_{\iy,i})}u_n,$$
	then
\begin{itemize}[fullwidth,itemindent=2em]
\item[(1)] 	$v_n(x_{n,i})\to\sqrt e$, $u_n(y_{n,i})\to\sqrt e$ and $x_{n,i}\to x_{\iy,i}$, $y_{n,i}\to x_{\iy,i}$ as $n\to\iy$.
\item[(2)]	as $n\to\iy$,
 				$$p_nu_n(x),p_nv_n(x)\to 8\pi\sqrt e \sum_{i=1}^kG(x,x_{\iy,i})\quad\text{in}~\CR_{loc}^2(\overline{\Omega}\setminus\SR).$$
\item[(3)]	the concentration points $\{x_i\}_{i=1}^k$ satisfy
				$$\nabla R(x_{\iy,i})-2\sum_{j=1,j\neq i}^k\nabla G(x_{\iy,i},x_{\iy,j})=0,\quad\text{for every $i=1,\cdots,k$}.$$
\item[(4)]	let $\mu_{n,i}:=\sbr{p_nv_n^{p_n-1}(x_{n,i})}^{-\f{1}{2}}$ and
\be\label{2-21} \bcs
		w_{n,i}(x):=\f{p_n}{v_n(x_{n,i})}(u_n(x_{n,i}+\mu_{n,i}x)-v_n(x_{n,i})),\\
		z_{n,i}(x):=\f{p_n}{v_n(x_{n,i})}(v_n(x_{n,i}+\mu_{n,i}x)-v_n(x_{n,i})),
\ecs \ee
then $(w_{n,i},z_{n,i})\to (U-\f{\q}{2},U)$ in $\CR_{loc}^2(\R^2)$. While let $\tilde \mu_{n,i}:=\sbr{p_nu_n^{p_n-1}(y_{n,i})}^{-\f{1}{2}}$ and
\be \bcs
		\tilde w_{n,i}(x):=\f{p_n}{u_n(y_{n,i})}(u_n(y_{n,i}+\tilde \mu_{n,i}x)-u_n(y_{n,i})),\\
		\tilde z_{n,i}(x):=\f{p_n}{u_n(y_{n,i})}(v_n(y_{n,i}+\tilde \mu_{n,i}x)-u_n(y_{n,i})),
\ecs\ee
then $(\tilde{w}_{n,i},\tilde{z}_{n,i})\to (U_\q-\f{\q}{2},U_\q)$ in $\CR_{loc}^2(\R^2)$, where
	$$U_\q(x)=\ln\f{e^{\f{\q}{2}}}{(1+\f{1}{8}e^{\f{\q}{2}}|x|^2)^2}$$
is a solution of the Liouville equation \eqref{Liouville}.
\item[(5)]  $p_n\int_\Omega\nabla u_n\cdot\nabla v_n\to k8\pi e$, $p_n\int_\Omega|\nabla u_n|^2\to k8\pi e$ and $p_n\int_\Omega|\nabla v_n|^2\to k8\pi e$ as $n\to\iy$.
\end{itemize}
\et

In this paper, we focus on  the solution sequence $(w_{n,i},z_{n,i})$, and the properties of the solution sequence $(\wt w_{n,i},\wt z_{n,i})$ can be obtained by similar arguments and hence the details are omitted.

Note that the $(w_{n,i},z_{n,i})$ defined in \eqref{2-21} satisfies
\be\label{c2-29}
	\bcs
		-\Delta w_{n,i}=(1+\f{z_{n,i}}{p_n})^{p_n} \quad\text{in}~\Omega_{n,i}:=\frac{\Omega-x_{n,i}}{\mu_{n,i}},\\
		-\Delta z_{n,i}=v_{n}(x_{n,i})^{\q_{n}}(1+\f{w_{n,i}}{p_n})^{q_n} \quad\text{in}~\Omega_{n,i},\\
		w_{n,i}(x)=z_{n,i}(x)=-p_n\quad\text{on}~\pa\Omega_{n,i},\\
w_{n,i}(x)>-p_n,\; z_{n,i}(x)>-p_n\quad\text{in}~\Omega_{n,i},\\
z_{n,i}(0)=0,\quad \nabla z_{n,i}(0)=0.
	\ecs
\ee
Furthermore, it follows from \eqref{est-v-1}-\eqref{est-u-1} that for $p_n$ large,
\be\label{c2-25}\int_{\Omega_{n,i}}\sbr{1+\f{z_{n,i}(x)}{p_n}}^{p_n}dx=\frac{p_n}{v_{n}(x_{n,i})}
\int_{\Omega}v_n(y)^{p_n}dy\leq C,\ee
\be\label{c2-26}\int_{\Omega_{n,i}}\sbr{1+\f{w_{n,i}(x)}{p_n}}^{q_n}dx=\frac{p_n}{v_{n}(x_{n,i})^{\theta_n+1}}
\int_{\Omega}u_n(y)^{q_n}dy\leq C.\ee

We need the following estimates, which will be used frequently in the later arguments, in particular when we apply the Dominated Convergence Theorem.

\begin{lemma}\label{decay-1} Fix any $\ga\in (0,4)$. Then there exist $C_\gamma>0$ and $n_\gamma$ such that
\begin{equation}\label{c2-27}
0\le \sbr{1+\f{w_{n,i}(x)}{p_n}}^{q_n},\; e^{w_{n,i}(x)}\le \f{C}{1+|x|^\ga}, \quad\text{for any }|x|\leq \frac{d_0}{\mu_{n,i}},
\end{equation}
\[0\le \sbr{1+\f{z_{n,i}(x)}{p_n}}^{p_n}\leq e^{z_{n,i}(x)}\le \f{C}{1+|x|^\ga}, \quad\text{for any }|x|\leq \frac{d_0}{\mu_{n,i}},\]
for any $n\geq n_\gamma$ and $1\leq i\leq k$. Here $d_0$ is defined in \eqref{constant-d0}.
\end{lemma}

\begin{proof}
First, we already proved in \cite[Lemma 4.3]{Chen-Li-Zou} that there exist small $r_\ga>0$, large $R_\ga>1$, $\tilde{n}_\ga>1$ and $\tilde{C}_\ga>0$ such that
\be  	\max\lbr{ w_{n,i}(x),~z_{n,i}(x)} \le \ga\ln\f{1}{|x|}+\tilde{C}_\ga, \quad\text{if}~2R_\ga\le |x|\le \f{r_\ga}{\mu_{n,i}}	\ee
for any $n\ge \tilde{n}_\ga$ and $i=1,\cdots,k$. Consequently, we have for $2R_\ga\le |x|\le \f{r_\ga}{\mu_{n,i}}$,
	$$\sbr{1+\f{w_{n,i}(x)}{p_n}}^{q_n}=e^{q_n\ln(1+\f{w_{n,i}(x)}{p_n})}\le e^{\f{q_n}{p_n}w_{n,i}(x)}
			\le \f{C}{|x|^{\f{q_n}{p_n}\ga}}\le \f{C}{|x|^{\ga}},$$
\[e^{w_{n,i}(x)}\le \f{C}{|x|^{\ga}}.\]
On the other hand, Theorem \ref{thm0}-(2) says that $u_n=O(\f{1}{p_n})$ in $B_{d_0}(x_{n,i})\setminus B_{r_\gamma}(x_{n,i})$. From here and $\mu_{n,i}=\sbr{p_nv_n^{p_n-1}(x_{n,i})}^{-\f{1}{2}}$, $v_n(x_{n,i})\to\sqrt e$, $q_n\geq p_n\to+\infty$, it is easy to see that $$0\leq \sbr{1+\f{w_{n,i}(x)}{p_n}}^{q_n}=\sbr{\f{u_n(x_{n,i}+\mu_{n,i}x)}{v_n(x_{n,i})}}^{q_n}\le \left(\f{C}{p_n}\right)^{q_n}\le C\mu_{n,i}^{\ga}\le \f{C}{|x|^\ga},$$
\[e^{w_{n,i}(x)}=e^{p_n(\f{u_n(x_{n,i}+\mu_{n,i}x)}{v_n(x_{n,i})}-1)}
\leq e^{p_n(\frac{C}{p_n}-1)}=e^{C-p_n}\le C\mu_{n,i}^{\ga}\le \f{C}{|x|^\ga}\]
hold for $\f{r_\ga}{\mu_{n,i}}\le |x|\le \f{d_0}{\mu_{n,i}}$ and $n$ large, where $\gamma<4$ is used to obtain $e^{C-p_n}\le C\mu_{n,i}^{\ga}$ for $n$ large.

Meanwhile, by $w_{n,i}\to U-\theta/2$ in $\CR_{loc}^2(\R^2)$, we have that for $|x|\le 2R_\ga$ and $n$ large,
 	$$\max\left\{\sbr{1+\f{w_{n,i}(x)}{p_n}}^{q_n}, e^{w_{n,i}(x)}\right\}\le Ce^{U(x)}=\f{C}{(1+\f{1}{8}|x|^2)^2}\le \f{C}{1+|x|^\ga}.$$
Therefore, we conclude that \eqref{c2-27} holds
for $n$ large.
The other estimates can be proved similarly.
\end{proof}

\begin{remark}\label{decay-10} In Sections 5-6, we will study the $1$-bubble solutions $(u_n, v_n)$, i.e. $k=i=1$, so we will omit the subscript $i$ in this case. Note that since $u_n, v_n=O(\f{1}{p_n})$ in $\Omega\setminus B_{r_\gamma}(x_{n})$, the proof of Lemma \ref{decay-1} actually implies that
\[0\le \sbr{1+\f{w_{n,i}(x)}{p_n}}^{q_n},\; e^{w_{n,i}(x)}\le \f{C}{1+|x|^\ga}, \quad\text{for any }x\in \frac{\Omega-x_n}{\mu_n},
\]
\[0\le \sbr{1+\f{z_{n,i}(x)}{p_n}}^{p_n}\leq e^{z_{n,i}(x)}\le \f{C}{1+|x|^\ga}, \quad\text{for any }x\in \frac{\Omega-x_n}{\mu_n},\]
for $n$ large.
\end{remark}

\begin{lemma}\label{decay-12} For any $1\leq i\leq k$, there exist $C>0$ and $n_0$ such that for any $n\geq n_0$,
\begin{equation}\label{c02-27}
|w_{n,i}(x)|,\; |z_{n,i}(x)|\le C\ln(2+|x|), \quad\text{if }|x|\leq \frac{d_0}{\mu_{n,i}}-1.
\end{equation}
\end{lemma}

\begin{proof} We only prove the estimate for $w_{n,i}(x)$ and the proof for $z_{n,i}(x)$ is similar and is omitted.

Note from $w_{n,i}\to U-\theta/2$ in $\CR_{loc}^2(\R^2)$ that $|w_{n,i}(x)|\leq C$ for any $|x|\leq 2$ and $n$ large. Thus we only need to consider $2\leq |x|\leq \frac{d_0}{\mu_{n,i}}-1$.

By \cite[Lemma 4.3]{Chen-Li-Zou}, we have
$$
	\begin{aligned}
		w_{n,i}(x)
		&=\f{1}{2\pi}\int_{\Omega_{n,i}} \ln\f{|y|}{|y-x|}\sbr{1+\f{z_{n,i}(y)}{p_n}}^{p_n}\rd y\\
			&\quad -\int_{\Omega_{n,i}} \mbr{H(x_{n,i}+\mu_{n,i}x,x_{n,i}+\mu_{n,i}y)-H(x_{n,i},x_{n,i}+\mu_{n,i}y)}\sbr{1+\f{z_{n,i}(y)}{p_n}}^{p_n}\rd y\\
		&=:\Rmnum{1}+\Rmnum{2},
	\end{aligned}
$$
where $H$ is the regular part of the Green function $G$.
Since $H$ is smooth and $|\mu_{n,i}x|\le d_0$, it follows from the differential mean value theorem and \eqref{c2-25} that $\Rmnum{2}=O(1)$.

Recall Lemma \ref{decay-1} that for $n$ large, we have \be\label{c2-23} 0\leq \sbr{1+\f{z_{n,i}(y)}{p_n}}^{p_n}\leq C\quad\text{ for any }\;|y|\leq \frac{d_0}{\mu_{n,i}}.\ee Then
\begin{align*}
&\left|\int_{\Omega_{n,i}} \ln\f{|y|}{|y|+1}\sbr{1+\f{z_{n,i}(y)}{p_n}}^{p_n}\rd y\right|\\
\leq &C\int_{|y|\leq 2} \left|\ln\f{|y|}{|y|+1}\right|\rd y+
C\int_{\Omega_{n,i}\cap\{|y|> 2\}} \sbr{1+\f{z_{n,i}(y)}{p_n}}^{p_n}\rd y\leq C,
\end{align*}
\begin{align*}
&\left|\int_{\Omega_{n,i}\cap\{|x-y|\leq 1\}} \ln\f{|y|+1}{|y-x|}\sbr{1+\f{z_{n,i}(y)}{p_n}}^{p_n}\rd y\right|\\
\leq &C\int_{|y-x|\leq 1} \ln\f{1}{|y-x|}\rd y+
C\ln |x|\int_{\Omega_{n,i}\cap\{|y-x|\leq 1\}} \sbr{1+\f{z_{n,i}(y)}{p_n}}^{p_n}\rd y\\
\leq & C+C\ln |x|.
\end{align*}
Note that for $|x-y|>1$ and $|x|\geq 2$, we have (see e.g. \cite[Lemma 8.23]{CL})
\[\left|\frac{\ln(|y|+1)-\ln |x-y|+\ln |x|}{\ln |x|}\right|\leq C,\]
i.e. $|\ln(|y|+1)-\ln |x-y||\leq C\ln |x|$, so
\begin{align*}
&\left|\int_{\Omega_{n,i}\cap\{|x-y|>1\}} \ln\f{|y|+1}{|y-x|}\sbr{1+\f{z_{n,i}(y)}{p_n}}^{p_n}\rd y\right|\\
\leq &
C\ln |x|\int_{\Omega_{n,i}\cap\{|y-x|>1\}} \sbr{1+\f{z_{n,i}(y)}{p_n}}^{p_n}\rd y\\
\leq & C\ln |x|.
\end{align*}
Inserting these estimates into
\begin{align*}
\Rmnum{1}=\f{1}{2\pi}\int_{\Omega_{n,i}} \ln\f{|y|}{|y|+1}\sbr{1+\f{z_{n,i}(y)}{p_n}}^{p_n}\rd y+\f{1}{2\pi}\int_{\Omega_{n,i}} \ln\f{|y|+1}{|y-x|}\sbr{1+\f{z_{n,i}(y)}{p_n}}^{p_n}\rd y,
\end{align*}
we conclude $|\Rmnum{1}|\leq C\ln|x|+C$ and so $|w_{n,i}(x)|\leq C\ln|x|+C$ for any $2\leq |x|\leq \frac{d_0}{\mu_{n,i}}-1$.
\end{proof}

Next, we expand $u_n,v_n$ at the local maximum point $x_{n,i}$. For any fixed $0<r<d_0$ and $i=1,\cdots,k$, we define
\be  	C_{n,i}=C_{n,i}(r):=\int_{B_r(x_{n,i})}v_n^{p_n}(x)\rd x\quad\text{and}\quad \wt C_{n,i}=\wt C_{n,i}(r):=\int_{B_r(x_{n,i})}u_n^{q_n}(x)\rd x.  \ee
From Theorem \ref{thm0}, Lemma \ref{decay-1} and the Dominated Convergence Theorem, we get \begin{align}\label{2-27}C_{n,i}&=\f{v_n(x_{n,i})}{p_n}\int_{B_{\f{r}{\mu_{n,i}}}(0)}\sbr{1+\f{z_{n,i}(y)}{p_n}}^{p_n}\rd y\nonumber\\
&=\f{\sqrt e}{p_n}\int_{\R^2}e^{U}dx+o(\f{1}{p_n})=\f{8\pi\sqrt e}{p_n}+o(\f{1}{p_n}),\end{align}
and similarly,
	\begin{equation}\label{2-28}\wt C_{n,i}=\f{v_n(x_{n,i})^{\q_n+1}}{p_n}\int_{B_{\f{r}{\mu_{n,i}}}(0)}\sbr{1+\f{w_{n,i}(y)}{p_n}}^{q_n}\rd y=\f{8\pi\sqrt e}{p_n}+o(\f{1}{p_n}).\end{equation}
We have
\bl\lab{expansion-1}
	For fixed $r\in(0,d_0)$, we have
	\be\bcs
		u_n(x)=\sum_{i=1}^k C_{n,i}G(x_{n,i},x)+o(\f{\mu_n}{p_n}),\\
		v_n(x)=\sum_{i=1}^k \wt C_{n,i} G(x_{n,i},x)+o(\f{\mu_n}{p_n}),
	\ecs\qquad \text{in}~\CR^1(\Omega\setminus\cup_{i=1}^kB_{2r}(x_{n,i})),
	\ee
	where $\mu_n=\max\lbr{\mu_{n,1},\cdots,\mu_{n,k}}$.
\el
\bp
For any $x\in\Omega\setminus\cup_{i=1}^kB_{2r}(x_{n,i})$, by the Green's representation formula we have
\be\label{c2-30}\begin{aligned}
	u_n(x)
	&=\int_\Omega G(y,x)v_n^{p_n}(y)\rd y\\
	&=\sum_{i=1}^k\int_{B_{r}(x_{n,i})} G(y,x)v_n^{p_n}(y)\rd y +\int_{\Omega\setminus\cup_{i=1}^kB_{r}(x_{n,i})} G(y,x)v_n^{p_n}(y)\rd y.
\end{aligned}\ee
Form Theorem \ref{thm0}-(2), we see that \begin{equation}\label{c2-31}u_n=O(\f{1}{p_n}),\quad v_n=O(\f{1}{p_n})\quad\text{ in }\;\overline{\Omega}\setminus\cup_{i=1}^kB_{r}(x_{n,i}),\end{equation} so
	$$\int_{\Omega\setminus\cup_{i=1}^kB_{r}(x_{n,i})} G(y,x)v_n^{p_n}(y)\rd y=O(\f{C^{p_n}}{p_n^{p_n}})\int_{\Omega} G(y,x)\rd y=O(\f{C^{p_n}}{p_n^{p_n}})=o(\f{\mu_n}{p_n}).$$
Since $\nabla^2 G(x,y)=O(1)$ for any $x\in\Omega\setminus B_{2r}(x_{n,i})$ and $y\in B_r(x_{n,i})$,  we obtain by the Taylor's expansion that
\be\begin{aligned}
	\int_{B_{r}(x_{n,i})} G(y,x)v_n^{p_n}(y)\rd y
&=C_{n,i}G(x_{n,i},x)+\int_{B_{r}(x_{n,i})} (G(y,x)-G(x_{n,i},x))v_n^{p_n}(y)\rd y \\
	&=C_{n,i}G(x_{n,i},x)+\sum_{h=1}^2\pa_hG(x_{n,i},x)\int_{B_{r}(x_{n,i})} (y-x_{n,i})_hv_n^{p_n}(y)\rd y \\
	&\quad +O\sbr{ \int_{B_{r}(x_{n,i})} |y-x_{n,i}|^2v_n^{p_n}(y)\rd y }.
\end{aligned}\ee
Using the Dominated Convergence Theorem and the property of odd functions, we get $$\lim_{n\to+\infty}\int_{B_{\f{r}{\mu_{n,i}}}(0)}y_h\sbr{1+\f{z_{n,i}(y)}{p_n}}^{p_n}\rd y=\int_{\R^2}y_he^{U(y)}\rd y=0,$$
so
	$$\int_{B_{r}(x_{n,i})} (y-x_{n,i})_hv_n^{p_n}(y)\rd y=\f{\mu_{n,i}v_n(x_{n,i})}{p_n}\int_{B_{\f{r}{\mu_{n,i}}}(0)}y_h\sbr{1+\f{z_{n,i}(y)}{p_n}}^{p_n}\rd y=o(\f{\mu_{n,i}}{p_n}).$$
Using Lemma \ref{decay-1} with $\gamma=4-\delta$ where $\delta\in (0,1)$ is small, we also get
$$\begin{aligned}
\int_{B_{r}(x_{n,i})} |y-x_{n,i}|^2v_n^{p_n}(y)\rd y
&=\f{\mu_{n,i}^2v_n(x_{n,i})}{p_n}\int_{B_{\f{r}{\mu_{n,i}}}(0)}|y|^2\sbr{1+\f{z_{n,i}(y)}{p_n}}^{p_n}\rd x\\
&\le C\f{\mu_{n,i}^2}{p_n}\int_{B_{\f{r}{\mu_{n,i}}}(0)}\f{|y|^2}{(1+|y|)^{4-\delta}}\rd x\\
&=O(\f{\mu_{n,i}^{2-\delta}}{p_n})=o(\f{\mu_n}{p_n}).
\end{aligned}$$
Inserting these estimates into \eqref{c2-30} leads to
\be\lab{tem-10}	u_n(x)=\sum_{i=1}^k C_{n,i}G(x_{n,i},x)+o(\f{\mu_n}{p_n}).  \ee
Similarly, for $h=1,2$, by similar arguments we can prove
	$$\pa_h u_n(x)=\int_\Omega D_hG(y,x)v_n^{p_n}(y)\rd y=\sum_{i=1}^k C_{n,i}D_hG(x_{n,i},x)+o(\f{\mu_n}{p_n}), $$
so \eqref{tem-10} holds in $\CR^1(\Omega\setminus\cup_{i=1}^kB_{2r}(x_{n,i}))$. The assertion for $v_n(x)$ can be proved similarly.
\ep

\bl
For fixed $r\in(0,d_0)$ and $i=1,\cdots,k$, it holds
\be\lab{tem-17}
	\int_{B_{\f{r}{\mu_{n,i}}}(0)} \sbr{1+\f{z_{n,i}}{p_n}}^{p_n}\rd x-v_n^{\q_n}(x_{n,i})\int_{B_{\f{r}{\mu_{n,i}}}(0)} \sbr{1+\f{w_{n,i}}{p_n}}^{q_n}\rd x=O(\f{1}{p_n}).
\ee
\el
\bp
By the Green's representation formula and \eqref{c2-31}, we have
\be\label{c2-32}\begin{aligned}
	u_n(x_{n,i})
	&=\int_\Omega G(x_{n,i},y)v_n^{p_n}(y)\rd y\\
	&=\sum_{j=1}^k \int_{B_r(x_{n,j})} G(x_{n,i},y)v_n^{p_n}(y)\rd y +O(\f{C^{p_n}}{p_n^{p_n}}).
\end{aligned}\ee
Note that we already proved in \cite[Lemma 4.4]{Chen-Li-Zou} that
$$\begin{aligned}
	\int_{B_r(x_{n,i})} G(x_{n,i},y)v_n^{p_n}(y)\rd y
	&=\f{v_n(x_{n,i})}{p_n}\int_{B_{\f{r}{\mu_{n,i}}}(0)} G(x_{n,i},x_{n,i}+\mu_{n,i}y)\sbr{ 1+\f{z_{n,i}(y)}{p_n} }^{p_n}\rd y\\
	&=-\f{v_n(x_{n,i})}{2\pi}\f{\ln \mu_{n,i}}{p_n}\int_{B_{\f{r}{\mu_{n,i}}}(0)} \sbr{1+\f{z_{n,i}(y)}{p_n}}^{p_n}\rd y+O(\f{1}{p_n}).
\end{aligned}$$
While for $j\neq i$, it holds $G(x_{n,i},y)=O(1)$ for $y\in B_r(x_{n,j})$ and hence \eqref{est-v-1} implies
	$$\int_{B_r(x_{n,j})} G(x_{n,i},y)v_n^{p_n}(y)\rd y=O\sbr{ \int_{B_r(x_{n,j})} v_n^{p_n}(y)\rd y }=O(\f{1}{p_n}).$$
Inserting these estimates into \eqref{c2-32}, we get
\begin{align}\lab{tem-16}
	-\f{1}{2\pi}\f{\ln \mu_{n,i}}{p_n}\int_{B_{\f{r}{\mu_{n,i}}}(0)} \sbr{1+\f{z_{n,i}(y)}{p_n}}^{p_n}\rd y&=\frac{u_n(x_{n,i})}{v_n(x_{n,i})}+O(\frac1{p_n})\nonumber\\
&=1+\f{w_{n,i}(0)}{p_n}+O(\f{1}{p_n})=1+O(\f{1}{p_n}).
\end{align}

Similarly, we have
\be \label{c2-33}\begin{aligned}
v_n(x_{n,i})
&=\int_\Omega G(x_{n,i},y)u_n^{q_n}(y)\rd y\\
&=\sum_{j=1}^k \int_{B_r(x_{n,j})} G(x_{n,i},y)u_n^{q_n}(y)\rd y +O(\f{C^{q_n}}{p_n^{q_n}}),
\end{aligned}\ee
$$\begin{aligned}
	\int_{B_r(x_{n,i})} G(x_{n,i},y)u_n^{q_n}(y)\rd y
	&=\f{v_n^{\q_n+1}(x_{n,i})}{p_n}\int_{B_{\f{r}{\mu_{n,i}}}(0)} G(x_{n,i},x_{n,i}+\mu_{n,i}y)\sbr{ 1+\f{w_{n,i}(y)}{p_n} }^{q_n}\rd y, \\
	&=-\f{v_n^{\q_n+1}(x_{n,i})}{2\pi}\f{\ln \mu_{n,i}}{p_n} \int_{B_{\f{r}{\mu_{n,i}}}(0)} \sbr{ 1+\f{w_{n,i}(y)}{p_n} }^{q_n}\rd y +O(\f{1}{p_n}),
\end{aligned}$$
and
	$$\int_{B_r(x_{n,j})} G(x_{n,i},y)u_n^{q_n}(y)\rd y=O\sbr{ \int_{B_r(x_{n,j})} u_n^{q_n}(y)\rd y }=O(\f{1}{p_n}),\quad \forall j\neq i.$$
Thus we get
\be\lab{tem-18}
	-\f{v_n^{\q_n}(x_{n,i})}{2\pi}\f{\ln \mu_{n,i}}{p_n} \int_{B_{\f{r}{\mu_{n,i}}}(0)} \sbr{ 1+\f{w_{n,i}(y)}{p_n} }^{q_n}\rd y =1+O(\f{1}{p_n}).
\ee
Comparing \eqref{tem-16} and \eqref{tem-18} and using
\[\f{\ln \mu_{n,i}}{p_n}=-\frac12\left(\frac{\ln p_n}{p_n}+\frac{p_n-1}{p_n}\ln v_n(x_{n,i})\right)\to -\frac14,\] we obtain \eqref{tem-17}.
\ep

\section{Sharp Estimates of \texorpdfstring{$v_n(x_{n,i})-u_n(x_{n,i})$}{} }
By Theorem \ref{thm0}-(4), we have \begin{equation}\label{c3-1}w_{n,i}(0)=p_n\left(\frac{u_n(x_{n,i})}{v_n(x_{n,i})}-1\right)\to U(0)-\f{\q}{2}=-\f{\q}{2},\end{equation} i.e.
	$$v_n(x_{n,i})-u_n(x_{n,i})=O(\f{1}{p_n}).$$
In this section, we give a sharper estimate of $v_n(x_{n,i})-u_n(x_{n,i})$, which plays an essential role in the proof of Theorem \ref{thm1}.
\begin{proposition}\lab{sharpest-0}
For any $i=1,\cdots,k$, it holds
\be\lab{keyest-1}  v_n(x_{n,i})-u_n(x_{n,i})=\f{1}{p_n}\q_n v_n(x_{n,i})\ln v_n(x_{n,i})+O(\f{1}{p_n^2}). \ee
\end{proposition}

Define notations
\be\lab{gamma} 		\sigma_{n,i}:=\q_n\ln v_n(x_{n,i}), \ee
and
\be\lab{epsilon} \e_{n,i}^{-1}:=p_n\sbr{\f{u_n(x_{n,i})}{v_n(x_{n,i})}-1}+\sigma_{n,i}=w_{n,i}(0)+\sigma_{n,i}. \ee
Then $\sigma_{n,i}\to\f{\q}{2}$ and so \eqref{c3-1} implies $\e_{n,i}\to\iy$ as $n\to\iy$.

We will prove Proposition \ref{sharpest-0} by contradiction. Suppose \eqref{keyest-1} does not hold for some $i\in\{1,\cdots,k\}$, i.e. up to a subsequence,
\be\lab{tem-101}  -\frac{p_n}{\e_{n,i}}v_n(x_{n,i})=p_n^2\sbr{v_n(x_{n,i})-u_n(x_{n,i})-\f{1}{p_n}\q_n v_n(x_{n,i})\ln v_n(x_{n,i})}\to\iy,\;\text{as}~n\to\iy, \ee
and we assume \eqref{tem-101} holds in this whole section. Clearly this implies $\e_{n,i}/p_n\to 0$ and so
\be\label{c3-07}\lim_{n\to\infty}\e_{n,i}\mu_{n,i}^{\tilde{\tau}}=
\lim_{n\to\infty}\frac{\e_{n,i}}{p_n}\frac{p_n}{p_n^{\frac12\tilde{\tau}}v_n(x_{n,i})^{\frac{p_n-1}{2}\tilde{\tau}}}=0,\quad\forall \tilde{\tau}>0.\ee

Let
\be\label{c3-6} 	s_{n,i}:=\e_{n,i}(w_{n,i}+\sigma_{n,i}-U)\quad\text{and}\quad t_{n,i}:=\e_{n,i}(z_{n,i}-U). 	\ee
where $U(x)=-2\ln(1+\f{1}{8}|x|^2)$. Then
\be\label{c3-7} s_{n,i}(0)=1,\quad t_{n,i}(0)=0,\quad \nabla t_{n,i}(0)=0.\ee

\bl\lab{limitsys-0}
	For fixed $r\in(0,d_0)$, we have
	\be\bcs
		-\Delta s_{n,i}=t_{n,i}(e^U+g_{n,i})+h_{n,i},\quad\text{in}~B_{\f{r}{\mu_{n,i}}}(0),\\
		-\Delta t_{n,i}=s_{n,i}(e^U+\wt g_{n,i})+\wt h_{n,i},\quad\text{in}~B_{\f{r}{\mu_{n,i}}}(0),
	\ecs\ee
	where
	\be 	g_{n,i}=O\sbr{\f{|z_{n,i}-U|}{1+|x|^{4-\delta}}}, \quad \wt g_{n,i}=O\sbr{\f{|w_{n,i}+\sigma_{n,i}-U|}{1+|x|^{4-\delta}}}\quad\text{in }B_{\f{r}{\mu_{n,i}}}(0), \ee
	\be 	\max\lbr{|h_{n,i}|,|\wt h_{n,i}|}=O\sbr{\f{\e_{n,i}}{p_n}\f{1}{1+|x|^{4-\delta}}}=o\sbr{\f{1}{1+|x|^{4-\delta}}} \quad\text{in }B_{\f{r}{\mu_{n,i}}}(0).	\ee
for any given small $\delta>0$.
\el
\bp
By \eqref{c2-29} and \eqref{c3-6}, we have
\be\label{c3-11}\begin{aligned}
-\Delta s_{n,i}
&=\e_{n,i}\sbr{ \sbr{1+\f{z_{n,i}}{p_n}}^{p_n}-e^U }\\
&=\e_{n,i}\sbr{ e^{z_{n,i}}-e^U }+\e_{n,i}\sbr{ \sbr{1+\f{z_{n,i}}{p_n}}^{p_n}-e^{z_{n,i}} }\\
&=t_{n,i}(e^U+g_{n,i})+h_{n,i},
\end{aligned}\ee
where we denote
\be \label{c3-2}	
g_{n,i}:=\f{\e_{n,i}\sbr{e^{z_{n,i}}-e^U}-t_{n,i}e^U}{t_{n,i}} \quad\text{and}\quad h_{n,i}:=\e_{n,i}\sbr{ \sbr{1+\f{z_{n,i}}{p_n}}^{p_n}-e^{z_{n,i}} }.
\ee
Now we estimate $g_{n,i}$ and $h_{n,i}$. By the Taylor's expansion and Lemma \ref{decay-1} with $\tau=4-\delta$, we get
\be\lab{tem-11-1} 	\begin{aligned}
\e_{n,i}\sbr{ e^{z_{n,i}}-e^U }
&=\e_{n,i}e^U\sbr{ z_{n,i}-U+O\sbr{|z_{n,i}-U|^2e^{\la(z_{n,i}-U)}} }\quad\text{(where $\lambda\in (0,1)$)}\\
&=t_{n,i}e^U+O\sbr{\f{t_{n,i}(z_{n,i}-U)}{1+|x|^{4-\delta}}},
\end{aligned}\ee
which implies $g_{n,i}=O\sbr{\f{|z_{n,i}-U|}{1+|x|^{4-\delta}}}$ in $B_{\f{r}{\mu_{n,i}}}(0)$.

To estimate $h_{n,i}$, we denote \begin{equation}\label{c3-15}y_{n,i}(x):=p_n\left(e^{\frac{z_{n,i}(x)}{p_n}}-1\right)\quad\text{ such that
}\quad e^{z_{n,i}(x)}=\left(1+\frac{y_{n,i}(x)}{p_n}\right)^{p_n}.\end{equation} Remark from Theorem \ref{thm0} (1) and (4) that for $n$ large,
\be\label{c3-031}-1<\frac{z_{n,i}(x)}{p_n}, \; \frac{w_{n,i}(x)}{p_n}\leq 1, \quad \forall |x|\leq r/\mu_{n,i}.\ee
(Indeed $z_{n,i}(x)\leq 0$.) Then
\begin{align}\label{c3-17}
y_{n,i}(x)-z_{n,i}(x)
&=p_n\left(\frac{z_{n,i}(x)}{p_n}+\frac12\frac{z_{n,i}(x)^2}{p_n^2}+O(\frac{z_{n,i}(x)^3}{p_n^3})\right)-z_{n,i}(x)\nonumber\\
&=\frac12\frac{z_{n,i}(x)^2}{p_n}+O(\frac{z_{n,i}(x)^3}{p_n^2})\\
&=O(\frac{z_{n,i}(x)^2}{p_n}),\quad \forall |x|\leq r/\mu_{n,i}\nonumber.
\end{align}
Since the function $(1+\frac{s}{p_n})^{p_n}$ is convex for $s>-p_n$, we see from Lemmas \ref{decay-1}-\ref{decay-12} that for $n$ large,
\begin{align}\label{tem-12-1}
\left|\frac{h_{n,i}(x)}{\e_{n,i}}\right|&=\left|\sbr{1+\f{z_{n,i}}{p_n}}^{p_n}
-\left(1+\frac{y_{n,i}}{p_n}\right)^{p_n}\right|\nonumber\\
&\leq \max\left\{\sbr{1+\f{z_{n,i}}{p_n}}^{p_n-1}, \left(1+\frac{y_{n,i}}{p_n}\right)^{p_n-1}\right\}|y_{n,i}(x)-z_{n,i}(x)|\nonumber\\
&\leq C\max\left\{\sbr{1+\f{z_{n,i}}{p_n}}^{p_n-1}, e^{\frac{p_n-1}{p_n}z_{n,i}}\right\}\frac{z_{n,i}(x)^2}{p_n}\nonumber\\
&\leq C\left(\frac{C}{1+|x|^{4-\frac{\delta}{2}}}\right)^{\frac{p_n-1}{p_n}}\frac{(\ln(|x|+2))^2}{p_n}
\leq \f{1}{p_n}\f{C}{1+|x|^{4-\delta}}
\end{align}
holds for any $|x|\leq r/\mu_{n,i}$. This proves the estimate for $h_{n,i}(x)$.

Also by \eqref{c2-29}, \eqref{c3-6} and \eqref{gamma}, we have
\be\begin{aligned}
-\Delta t_{n,i}
&=\e_{n,i}\sbr{ v_n^{\q_n}(x_{n,i})\sbr{1+\f{w_{n,i}}{p_n}}^{q_n}-e^U }\\
&=\e_{n,i}\sbr{ e^{w_{n,i}+\sigma_{n,i}}-e^U }+\e_{n,i}e^{\sigma_{n,i}}\sbr{ \sbr{1+\f{w_{n,i}}{p_n}}^{q_n}-e^{w_{n,i}} }\\
&=s_{n,i}(e^U+\wt g_{n,i})+\wt h_{n,i},
\end{aligned}\ee
where
\be 	\wt g_{n,i}:=\f{\e_{n,i}\sbr{ e^{w_{n,i}+\sigma_{n,i}}-e^U }-s_{n,i}e^U}{s_{n,i}} \ee
and
\be 	\wt h_{n,i}:=\e_{n,i}e^{\sigma_{n,i}}\sbr{ \sbr{1+\f{w_{n,i}}{p_n}}^{q_n}-e^{w_{n,i}} }. \ee
The desired estimates for $\wt g_{n,i}$ and $\wt h_{n,i}$ can be obtained similarly as \eqref{tem-11-1}-\eqref{tem-12-1} and we omit the details.
\ep

To get the convergence information, we have to prove the a priori estimate of $(s_{n,i},t_{n,i})$. Actually we need to establish the following lemma:
\bl\lab{bound-0}
For fixed $r\in(0,d_0)$ and any $\tau\in(0,1)$, there exists a $C_\tau>0$ such that
\be
	\max\lbr{|s_{n,i}(x)|,~|t_{n,i}(x)|}\le C_\tau(1+|x|)^\tau\quad\text{for any}~x\in B_{\f{r}{\mu_{n,i}}}(0).
\ee
\el
The proof of Lemma \ref{bound-0} is complex and technical. The idea of the proof origins to Estimate C in \cite{Lin=CPAM=2002}, and we also refer to \cite[Proposition 3.5]{LE-1} for the case of the scalar equation. Since all the arguments are in local sense, to simplify the notations, we use $z_n,w_n,s_n,t_n,\mu_n,x_n,\e_n,\ga_n$ to denote $z_{n,i},w_{n,i}$, etc, i.e. we omit the subscript $i$. Let
\be 	N_n:=\max_{|x|\le\f{r}{\mu_n}}\lbr{ \f{|s_n(x)|}{(1+|x|)^\tau},~\f{|t_n(x)|}{(1+|x|)^\tau} }. \ee
To prove Lemma \ref{bound-0}, we have to show that $N_n$ is bounded. We prove it by contradiction. Assume $N_n\to\iy$ as $n\to\iy$, then we have the following claim:
\bl\lab{claim-1}
If $N_n\to\iy$, then $N_n^*=o(1)N_n$, where
\be  N_n^*:=\max_{|x|\le\f{r}{\mu_n}}\max_{|x'|=|x|}\lbr{ \f{|s_n(x)-s_n(x')|}{(1+|x|)^\tau},~\f{|t_n(x)-t_n(x')|}{(1+|x|)^\tau} }. \ee
\el
\bp
Assume by contradiction that up to a subsequence, $N_n^*\ge CN_n$ for some $C>0$. Without loss of generality, we take $x_n'$ and $x_n''$ such that $|x_n'|=|x_n''|\le\f{r}{\mu_n}$ and
	$$N_n^*=\f{|s_n(x_n')-s_n(x_n'')|}{(1+|x_n'|)^\tau}.$$
By rotation, we assume that $x_n'=(x_{n,1}', x_{n,2}')$ with $x'_{n,2}>0$ and $x_n''=x_n'^-$, where $x^-:=(x_1,-x_2)$ for $x=(x_1,x_2)$. Set
\begin{align}\label{c3-3} 	&\delta_n^*(x):=s_n(x)-s_n(x^-)=\varepsilon_n(w_n(x)-w_n(x^-)),\nonumber\\ &\ga_n^*(x):=t_n(x)-t_n(x^-)=\varepsilon_n(z_n(x)-z_n(x^-)), \end{align}
and\[\delta_n(x):=\f{\delta_n^*(x)}{(1+x_2)^\tau},\quad\quad \ga_n(x):=\f{\ga_n^*(x)}{(1+x_2)^\tau}.\]
Define
\be \label{c3-23}	N_n^{**}:=\max_{|x|\le\f{r}{\mu_n},x_2\ge0}\lbr{ |\delta_n(x)|,~|\ga_n(x)| }, \ee
then
\[N_n^{**}\geq \f{|\delta_n^*(x_n')|}{(1+x_{n,2}')^\tau}\geq \f{|\delta_n^*(x_n')|}{(1+|x_n'|)^\tau}=N_n^*\ge C N_n\to +\infty.\]
Without loss of generality we take $x_n^*\in \lbr{x\in\R^2:|x|\le\f{r}{\mu_n},x_2>0}$ such that $N_n^{**}=|\delta_n(x_n^*)|\geq |\ga_n(x_n^*)|$, where $x_{n,2}^*>0$ follows from $\delta_n(x_1,0)=\ga_n(x_1,0)\equiv 0$.

{\bf Step 1.} We prove $|x_n^*|\le \f{r}{2\mu_n}$.

By contradiction, we assume $\f{r}{2\mu_n}\le |x_n^*|\le \f{r}{\mu_n}$. Using \eqref{c3-6}, \eqref{2-21}, $U(x_n^*)=U(x_n^{*-})$ and Lemma \ref{expansion-1}, we have
\be\lab{tem-28} \begin{aligned}
\delta_n^*(x_n^*)
&=s_n(x_n^*)-s_n(x_n^{*-})\\
&=\f{\e_np_n}{v_n(x_{n,i})}(u_n(x_{n,i}+\mu_{n,i}x_n^*)-u_n(x_{n,i}+\mu_{n,i}x_n^{*-}))\\
&=\f{\e_np_n}{v_n(x_{n,i})}\sum_{j=1}^k C_{n,j}(G(x_{n,j},x_{n,i}+\mu_{n,i}x_n^*)-G(x_{n,j},x_{n,i}+\mu_{n,i}x_n^{*-}))+o(\e_n\mu_n).
\end{aligned}\ee
Note from \eqref{2-27} that $p_nC_{n,j}=O(1)$ for all $j$. Since $\f{r}{2\mu_n}\le |x_n^*|\le \f{r}{\mu_n}$, we obtain that for $j\neq i$,
	$$G(x_{n,j},x_{n,i}+\mu_{n,i}x_n^*)-G(x_{n,j},x_{n,i}+\mu_{n,i}x_n^{*-})
=O(\mu_{n,i}|x_n^*-x_n^{*-}|)=O(\mu_{n,i}x_{n,2}^*).$$
While for $j=i$, we have
	$$\begin{aligned}
	&\quad G(x_{n,i},x_{n,i}+\mu_{n,i}x_n^*)-G(x_{n,i},x_{n,i}+\mu_{n,i}x_n^{*-})\\
	&=-H(x_{n,i},x_{n,i}+\mu_{n,i}x_n^*)+H(x_{n,i},x_{n,i}+\mu_{n,i}x_n^{*-})\\
	&=O(\mu_{n,i}|x_n^*-x_n^{*-}|)=O(\mu_{n,i}x_{n,2}^*).
	\end{aligned}$$
Thus $\delta_n^*(x_n^*)=O(\e_n\mu_nx_{n,2}^*)+o(\e_n\mu_n)$ and then we see from $x_{n,2}^*\leq \frac{r}{\mu_{n}}$ and \eqref{c3-07} that
\be\label{c3-28} 	N_n^{**}=\f{|\delta_n^*(x_n^*)|}{(1+x_{n,2}^*)^\tau}=O(\e_n\mu_n(x_{n,2}^*)^{1-\tau})+o(\e_n\mu_n)=o(1), 	\ee
which is a contradiction with $N_n^{**}\to\iy$. So $|x_n^*|\le \f{r}{2\mu_n}$.

{\bf Step 2.} We prove $|x_n^*|\le C$.

By contradiction, we assume $|x_n^*|\to\iy$. Since $N_n^{**}=|\delta_n(x_n^*)|\to\iy$, we may assume $\delta_n(x_n^*)>0$ for large $n$ (The case $\delta_n(x_n^*)<0$ can be discussed similarly and we omit the details). By \eqref{c3-23}, $|x_n^*|\le \f{r}{2\mu_n}$ and $x_{n,2}^*>0$, we get
	$$\nabla \delta_n(x_n^*)=0 \quad\text{and}\quad \Delta \delta_n(x_n^*)\le 0.$$
Let
$$\LR :=-\Delta-\f{2\tau}{1+x_2}\f{\pa}{\pa x_2}+\f{\tau(1-\tau)}{(1+x_2)^2},$$
then
\be\label{c3-26} -\Delta [(1+x_2)^{\tau}\varphi(x)]=(1+x_2)^{\tau}\LR \varphi(x).\ee
From here, \eqref{c3-11} and $U(x)=U(x^-)$ we have
\begin{align} \label{c3-027} 	\LR \delta_n &=\frac{1}{(1+x_2)^{\tau}}(-\Delta s_n(x)+\Delta s_n(x^-))\nonumber\\
&=\frac{\varepsilon_n}{(1+x_2)^{\tau}}\left[\left(1+\frac{z_{n}(x)}{p_n}\right)^{p_n}
-\left(1+\frac{z_{n}(x^-)}{p_n}\right)^{p_n}\right]
=:\f{g^*_n(x)}{(1+x_2)^\tau}.	\end{align}
Then we obtain
\be\lab{tem-27} 	0<\f{\tau(1-\tau)}{(1+x_{n,2}^*)^2}\delta_n(x_n^*)\le \LR \delta_n(x_n^*)= \f{g^*_n(x_n^*)}{(1+x_{n,2}^*)^\tau}, \ee
where $\tau\in (0,1)$ is used.
On the other hand, by \eqref{c3-3}, $\ln(1+s)\leq s$ and Lemma \ref{decay-1},
\begin{align}\label{c3-029}
0\leq\frac{g^*_n(x)}{\ga_n^*(x)}&=\frac{(1+\frac{z_{n}(x)}{p_n})^{p_n}
-(1+\frac{z_{n}(x^-)}{p_n})^{p_n}}{z_{n}(x)-z_{n}(x^-)}\nonumber\\
&=\left(1+\frac{\lambda z_{n}(x)+(1-\lambda)z_{n}(x^-)}{p_n}\right)^{p_n-1}\;\;(\text{where $\lambda\in (0,1)$})\\
&\leq e^{\frac{p_n-1}{p_n}(\lambda z_{n}(x)+(1-\lambda)z_{n}(x^-))}\nonumber\\
&\leq \left(\frac{C}{1+|x|^{4-\delta}}\right)^{\frac{p_n-1}{p_n}},\quad\forall |x|\leq \frac{r}{\mu_n}.\nonumber
\end{align}
From here, $\delta_n(x_n^*)=|\delta_n(x_n^*)|\geq |\ga_n(x_n^*)|$ and $|x_n^*|\to\iy$,  we have for $n$ large,
	$$\f{g^*_n(x_n^*)}{(1+x_{n,2}^*)^\tau}\le \left(\frac{C}{1+| x_n^*|^{4-\delta}}\right)^{\frac{p_n-1}{p_n}}\frac{\ga_n(x_n^*)}{(1+x_{n,2}^*)^\tau} \le \f{1}{2}\f{\tau(1-\tau)}{(1+x_{n,2}^*)^2}\delta_n(x_n^*),$$
a contradiction with \eqref{tem-27}. So $|x_n^*|\le C$.

{\bf Step 3.} We complete the proof by obtaining a contradiction.

Define $\bar\delta_n:=\f{\delta_n}{N_n^{**}}$ and $\bar\ga_n:=\f{\ga_n}{N_n^{**}}$. Then
\be\label{c3-29} \bcs
	\LR\bar\delta_n=\f{g_n^* }{N_n^{**}(1+x_2)^\tau},\quad\text{in}~B_{\f{r}{\mu_n}}(0)\cap\lbr{x\in\R^2:x_2\ge0},\\
	\LR\bar\ga_n=\f{\wt g_n^* }{N_n^{**}(1+x_2)^\tau},\quad\text{in}~B_{\f{r}{\mu_n}}(0)\cap\lbr{x\in\R^2:x_2\ge0},\\
	|\bar\delta_n(x_n^*)|=1=\max\lbr{|\bar\delta_n|,|\bar\ga_n|},
\ecs \ee
where similarly as \eqref{c3-027} and \eqref{c3-029},
\begin{align}\label{c3-030} 	\wt g_n^*(x):=&\e_nv_n(x_n)^{\theta_n}\left[\left(1+\frac{w_{n}(x)}{p_n}\right)^{q_n}
-\left(1+\frac{w_{n}(x^-)}{p_n}\right)^{q_n}\right]\nonumber\\
=&v_n(x_n)^{\theta_n}\delta_n^*(x)\frac{\left(1+\frac{w_{n}(x)}{p_n}\right)^{q_n}
-\left(1+\frac{w_{n}(x^-)}{p_n}\right)^{q_n}}{w_{n}(x)-w_{n}(x^-)}\nonumber\\
=&\frac{q_n}{p_n}v_n(x_n)^{\theta_n}\delta_n^*(x)\left(1+\frac{\tilde{\lambda} w_{n}(x)+(1-\tilde{\lambda})w_{n}(x^-)}{p_n}\right)^{q_n-1}\;\;(\text{where $\tilde{\lambda}\in (0,1)$}).\end{align}
Since Theorem \ref{thm0} says that $z_{n}(x), z_n(x^-)\to U(x)$ and $w_n(x), w_n(x^-)\to U-\frac{\theta}{2}$ in $\mathcal{C}_{\loc}^2(\mathbb{R}^2)$, we see from \eqref{c3-029} and \eqref{c3-030} that
\[\frac{\frac{g_n^* }{N_n^{**}(1+x_2)^\tau}}{\bar\ga_n}=\frac{g^*_n(x)}{\ga_n^*(x)}\to e^{U(x)}\quad\text{in }\mathcal{C}_{\loc}(\mathbb{R}^2),\]
\[\frac{\frac{\tilde{g}_n^* }{N_n^{**}(1+x_2)^\tau}}{\bar\delta_n}=\frac{\tilde{g}^*_n(x)}{\delta_n^*(x)}\to e^{\frac{\theta}{2}}e^{U(x)-\frac{\theta}{2}}=e^{U(x)},\quad\text{in }\mathcal{C}_{\loc}(\mathbb{R}^2).\]
From here and \eqref{c3-29}, we conclude from the standard elliptic estimates that
up to a subsequence, $(\bar\delta_n,\bar\ga_n)\to(\bar\delta_\iy,\bar\ga_\iy)$ in $\CR_{loc}^2(\R^2_+)$ and
\be \bcs
	\LR\bar\delta_\iy= e^U\bar\ga_\iy,\quad\text{in}~\R^2_+,\\
	\LR\bar\ga_\iy= e^U\bar\delta_\iy,\quad\text{in}~\R^2_+,\\
\|\bar\delta_\iy\|_{L^\infty},\|\bar\ga_\iy\|_{L^\infty}\leq 1,\\
\bar\delta_\iy(x_1, 0)=\bar\ga_\iy(x_1, 0)=0,
\ecs \ee
where $\R^2_+=\lbr{x\in\R^2:x_2\ge0}$.
Since $|x_n^*|\le C$ we obtain $\max_{\R_+^2} \bar\delta_\iy=1$, which means $\bar\delta_\iy \not\equiv0$. Define $(\delta_\iy,\ga_\iy):=((1+x_2)^\tau\bar\delta_\iy,(1+x_2)^\tau\bar\ga_\iy)$ and consider its odd extension in $x_2$, then by \eqref{c3-26} we deduce that $(\delta_\iy,\ga_\iy)$ is a solution of
\be\lab{linearsys-01}\bcs
		-\Delta u=e^U v,\;-\Delta v=e^U u &\text{in}~\R^2,\\
		|u(x)|,|v(x)|\le C(1+|x|)^\tau&\text{in}~\R^2.
	\ecs\ee
Applying Lemma \ref{linear} leads to
	$$(\delta_\iy,\ga_\iy)=\sum_{j=0}^2 c_j(\phi_j,\phi_j)+c_3(\phi_3,-\phi_3).$$
Since $\nabla v_n(x_n)=0$, we obtain from \eqref{c3-3} and \eqref{2-21} that
	$$\nabla \ga_n^*(0)=\e_n(\nabla z_n(x)-\nabla z_n(x^-))|_{x=0}=\sbr{0,~\f{2\e_np_n\mu_n}{v_n(x_n)}\f{\pa}{\pa x_2}v_n(x_n)}= 0,$$
i.e. $\nabla \sbr{ (1+x_2)^\tau\bar\ga_n(x)}|_{x=0}=0$. Since $\ga_n(0)=0$, we obtain $\nabla \ga_n(0)=0$ and so $\nabla \ga_\iy(0)=0$. It follows from \eqref{eq1} that $c_1=c_2=0$. Moreover $\ga_\iy(x_1,0)\equiv0$ for any $x_1\in\R$ implies $c_0=c_3=0$. Thus $\delta_\iy\equiv0$, a contradiction with $\bar\delta_\iy\not\equiv0$.
The proof is complete.
\ep

Now we are in the position to complete the proof of Lemma \ref{bound-0}.
\bp[Proof of Lemma \ref{bound-0}] Fix any $r_0\in (0,d_0)$.
Suppose by contradiction that \[N_n=\max_{|x|\le\f{r_0}{\mu_n}}\lbr{ \f{|s_n(x)|}{(1+|x|)^\tau},~\f{|t_n(x)|}{(1+|x|)^\tau} }\to +\infty\quad\text{as }n\to\iy.\]
Note that
\be\label{c3-33}|s_n(x)|, |t_n(x)|\leq N_n (1+|x|)^\tau,\quad \forall |x|\leq r_0/\mu_n.\ee
Set
\be \label{c3-36}	\vp_n(r):=\f{1}{2\pi}\int_0^{2\pi}s_n(r,\rho)\rd \rho \quad\text{and}\quad \psi_n(r):=\f{1}{2\pi}\int_0^{2\pi}t_n(r,\rho)\rd \rho. \ee
Since
\[0\leq\max_{\rho\in[0,2\pi]}|s_n(r,\rho)|-|\vp_n(r)|\leq \max_{\rho_1,\rho_2\in[0,2\pi]}|s_n(r,\rho_1)-s_n(r,\rho_2)|,\] 
we obtain from Lemma \ref{claim-1} that
	$$\max_{0\le r\le \f{r_0}{\mu_n}}\lbr{ \f{|\vp_n(r)|}{(1+r)^\tau},~\f{|\psi_n(r)|}{(1+r)^\tau} }=N_n(1+o(1))\to +\infty.$$
Without loss of generality, we take $r_n\in [0,\f{r_0}{\mu_n}]$ such that
	\be\label{c3-35}\f{|\vp_n(r_n)|}{(1+r_n)^\tau}=\max_{0\le r\le \f{r_0}{\mu_n}}\lbr{ \f{|\vp_n(r)|}{(1+r)^\tau},~\f{|\psi_n(r)|}{(1+r)^\tau} }.\ee

{\bf Step 1. }
We claim that $r_n$ is uniformly bounded.

To see it, let $\phi_0(x)=\frac{8-|x|^2}{8+|x|^2}$ be defined in \eqref{eq1}, then $|\phi_0(x)|\leq 1$ and $-\Delta \phi_0=e^U\phi_0$. Recalling Lemma \ref{limitsys-0} and Lemma \ref{decay-12} that
\[g_{n}=O\sbr{\f{|z_{n}-U|}{1+|x|^{4-\delta}}}=O\sbr{\f{\ln(2+|x|)}{1+|x|^{4-\delta}}}
=O\sbr{\f{1}{1+|x|^{4-2\delta}}},\]
\[h_n=O\sbr{\frac{\e_n}{p_n}\f{1}{1+|x|^{4-\delta}}}=o\sbr{\f{1}{1+|x|^{4-\delta}}},\]
for $|x|\leq r_0/\mu_n$,  it follows from \eqref{c3-33} that
\be \begin{aligned}
\int_{|x|\le r}\phi_0\Delta s_n-s_n\Delta\phi_0\rd x
&=\int_{|x|\le r} e^U(s_n-t_n)\phi_0-t_ng_n\phi_0-h_n\phi_0\rd x\\
&\le CN_n\int_0^{r}\f{s(1+s)^\tau}{1+s^{4-2\delta}}\rd s=O(N_n).
\end{aligned}\ee
On the other hand, we have
\be \begin{aligned}
\int_{|x|\le r}\phi_0\Delta s_n-s_n\Delta\phi_0\rd x
&=\int_{|x|=r}\phi_0\pa_\nu s_n-s_n\pa_\nu\phi_0\rd \sigma\\
&=\f{8-r^2}{8+r^2}2\pi r\vp_n'(r)+\f{32r}{(8+r^2)^2}2\pi r\vp_n(r).
\end{aligned}\ee
From here and $\vp_n(r)\leq (1+r)^\tau N_n(1+o(1))$ we get
\be 	|\vp_n'(r)|\le C \sbr{ \f{N_n}{r}+\f{\vp_n(r)}{r^3} }=\f{O(N_n)}{r},\quad \forall 3\leq r\leq r_0/\mu_n, \ee
It follows that
	$$(1+r_n)^\tau N_n(1+o(1))=|\vp_n(r_n)|\le \int_3^{r_n}|\vp'(r)|\rd r+|\vp_n(3)|\le C(1+\ln r_n)N_n,$$
which gives that $r_n$ is uniformly bounded.

{\bf Step 2. }We complete the proof by obtaining a contradiction.

Let $\bar\vp_n(x):=\f{\vp_n(|x|)}{\vp_n(r_n)}$ and $\bar\psi_n(x):=\f{\psi_n(|x|)}{\vp_n(r_n)}$. Then \eqref{c3-35} implies $|\bar\vp_n(x)|\le (1+|x|)^\tau$ and $|\bar\psi_n(x)|\le (1+|x|)^\tau$ for $|x|\leq r_0/\mu_n$. By Lemma \ref{limitsys-0} and \eqref{c3-36}, we have
\be \bcs
	-\Delta\bar\vp_n=\bar\psi_n\sbr{ e^U+O\sbr{\f{|z_n-U|}{1+|x|^{4-\delta}}} }+O\sbr{\f{1}{N_n(1+|x|^{4-\delta})}},\quad\text{in}~B_{\f{r}{\mu_n}}(0),\\
	-\Delta\bar\psi_n=\bar\vp_n\sbr{ e^U+O\sbr{\f{|w_n+\sigma_n-U|}{1+|x|^{4-\delta}}} }+O\sbr{\f{1}{N_n(1+|x|^{4-\delta})}},\quad\text{in}~B_{\f{r}{\mu_n}}(0),\\
\ecs \ee
Since $z_n-U\to 0$ and $w_n+\sigma_n-U\to 0$ in $\mathcal{C}_{loc}^2(\mathbb{R}^2)$,
by standard elliptic estimates we obtain that up to a subsequence, $(\bar\vp_n,\bar\psi_n)\to(\bar\vp_\iy,\bar\psi_\iy)$ in $\CR_{loc}^2(\R^2)$, where $(\bar\vp_\iy,\bar\psi_\iy)$ is a solution of the linearzed system \eqref{linearsys-0}.  Since $r_n$ is bounded, we know $\bar\vp_\iy\not\equiv0$. Applying Lemma \ref{linear}, it holds
\be 	(\bar\vp_\iy,\bar\psi_\iy)=\sum_{j=0}^2c_j(\phi_j,\phi_j)+c_3(\phi_3,-\phi_3). \ee
The radial symmetry of $(\bar\vp_\iy,\bar\psi_\iy)$ yields $c_1=c_2=0$. Since \eqref{c3-7} says that $s_n(0)=1$ and $t_n(0)=0$, we see from $\vp(r_n)\to\infty$ that $\bar\vp_\iy(0)=\bar\psi_\iy(0)=0$, which yields $c_0=c_3=0$. So $\bar\vp_\iy\equiv0$, which is a contradiction with $\bar\vp_\iy\not\equiv0$. This finishes the proof of Lemma \ref{bound-0}.
\ep

As a consequence of Lemma \ref{bound-0}, we obtain
\bl\lab{tem-31}
It holds that up to a subsequence,
	$$(s_{n,i},t_{n,i})\to \f{1}{2}(\phi_0,\phi_0)+\f{1}{2}(\phi_3,-\phi_3),\quad\text{in}~\CR_{loc}^2(\R^2),$$
as $n\to\iy$, where $\phi_0$ and $\phi_3$ is given in Lemma \ref{linear}.
\el
\bp
By Lemmas \ref{limitsys-0}-\ref{bound-0}, it follows from 
the standard elliptic estimates that up to a subsequence, as $n\to\iy$, $(s_{n,i},t_{n,i})\to (s_{\iy,i},t_{\iy,i})$, and $(s_{\iy,i},t_{\iy,i})$ is a solution of the linearized system \eqref{linearsys-0}.
Again by Lemma \ref{linear}, we have
	$$(s_{\iy,i},t_{\iy,i})=\sum_{j=0}^2c_{i,j}(\phi_j,\phi_j)+c_{i,3}(\phi_3,-\phi_3).$$
Since \eqref{c3-7} says $s_{n,i}(0)=1$, $t_{n,i}(0)=0$ and $\nabla t_{n,i}(0)=0$, we have $s_{\iy,i}(0)=1$, $t_{\iy,i}(0)=0$ and $\nabla t_{\iy,i}(0)=0$, which implies $c_{i,0}=c_{i,3}=\f{1}{2}$ and $c_{i,1}=c_{i,2}=0$.
\ep

\bl\lab{tem-32}
For fixed $r\in(0,d_0)$, we have
\be\lab{est-1-1} 	
	\int_{B_{\f{r}{\mu_{n,i}}}(0)} \sbr{1+\f{z_{n,i}}{p_n}}^{p_n}\rd x=8\pi+\f{1}{\e_{n,i}}\int_{B_{\e_{n,i}}(0)} e^U t_{n,i}\rd x+o(\f{1}{\e_{n,i}}),
\ee
and
\be\lab{est-2-1} 	
	e^{\sigma_{n,i}}\int_{B_{\f{r}{\mu_{n,i}}}(0)} \sbr{1+\f{w_{n,i}}{p_n}}^{q_n}\rd x= 8\pi+\f{1}{\e_{n,i}}\int_{B_{\e_{n,i}}(0)} e^U s_{n,i}\rd x+o(\f{1}{ \e_{n,i} }).
\ee
\el
\bp
Note from \eqref{c3-07} that $\e_{n,i}<\f{r}{\mu_{n,i}}$ for $n$ large, then by Lemma \ref{decay-1} and $\e_{n,i}\to\infty$ we have
	$$\int_{B_{\f{r}{\mu_{n,i}}}(0)\setminus B_{\e_{n,i}}(0)} \sbr{1+\f{z_{n,i}}{p_n}}^{p_n}\rd x \le C\int_{B_{\f{r}{\mu_{n,i}}}(0)\setminus B_{\e_{n,i}}(0)} \f{1}{1+|x|^{4-\delta}}\rd x=O(\f{1}{\e_{n,i}^{2-\delta}})=o(\f{1}{\e_{n,i}}),$$
Thus
\be\lab{tem-22-1}   	
	\int_{B_{\f{r}{\mu_{n,i}}}(0)} \sbr{1+\f{z_{n,i}}{p_n}}^{p_n}\rd x=\int_{B_{\e_{n,i}}(0)} \sbr{1+\f{z_{n,i}}{p_n}}^{p_n}\rd x+o(\f{1}{\e_{n,i}}).
\ee

Note from \eqref{c3-7} that
\be\lab{expansion-w-1} 	w_{n,i}=U-\sigma_{n,i}+\f{s_{n,i}}{\e_{n,i}}, \quad z_{n,i}=U+\f{t_{n,i}}{\e_{n,i}}. \ee
For $x\in B_{\e_{n,i}}(0)$, Lemma \ref{decay-12} says that $z_{n,i}(x)=O(\ln \e_{n,i})$ and Lemma \ref{bound-0} says that $t_{n,i}(x)=O(\e_{n,i}^\tau)$, so by the Taylor's expansion and $\e_{n,i}/p_n\to 0$,
$$\begin{aligned}
	\sbr{1+\f{z_{n,i}}{p_n}}^{p_n}
	&=e^{p_n\ln(1+\f{z_{n,i}}{p_n})}=e^{z_{n,i}-\f{1}{2p_n}z_{n,i}^2+O(\f{z_{n,i}^3}{p_n^2})}\\
	&=e^{U+\f{t_{n,i}}{\e_{n,i}} -\f{1}{2p_n}\sbr{ U+\f{t_{n,i}}{\e_{n,i}} }^2 +o(\f{1}{\e_{n,i}}) }\\
	&=e^U\sbr{ 1+\f{t_{n,i}}{\e_{n,i}} -\f{1}{2p_n}U^2+o(\f{1}{\e_{n,i}}) }.
\end{aligned}$$
So
$$\begin{aligned}
	\int_{B_{\e_{n,i}}(0)} \sbr{1+\f{z_{n,i}}{p_n}}^{p_n}\rd x
	&=\int_{B_{\e_{n,i}}(0)}e^U\rd x + \f{1}{\e_{n,i}}\int_{B_{\e_{n,i}}(0)} e^U t_{n,i}\rd x+ O(\frac{1}{p_n})+o(\f{1}{\e_{n,i}})\\
	&=8\pi+\f{1}{\e_{n,i}}\int_{B_{\e_{n,i}}(0)} e^U t_{n,i}\rd x+o(\f{1}{\e_{n,i}}),
\end{aligned}$$
which together with \eqref{tem-22-1} gives \eqref{est-1-1}.

Similarly, for $n$ large we have
\be\lab{tem-23-1}   	
	\int_{B_{\f{r}{\mu_{n,i}}}(0)} \sbr{1+\f{w_{n,i}}{p_n}}^{q_n}\rd x  =\int_{B_{\e_{n,i}}(0)} \sbr{1+\f{w_{n,i}}{p_n}}^{q_n}\rd x+o(\f{1}{\e_{n,i}}).
\ee
Again in $B_{\e_{n,i}}(0)$, Lemma \ref{decay-12} says that $w_{n,i}=O(\ln \e_{n,i})$ and Lemma \ref{bound-0} says that $s_{n,i}(x)=O(\e_{n,i}^\tau)$. Since $\frac{q_n}{p_n}=1+\frac{\theta_n}{p_n}=1+O(\frac{1}{p_n})$, we have
$$\begin{aligned}
	\sbr{1+\f{w_{n,i}}{p_n}}^{q_n}
	&=e^{q_n\ln(1+\f{w_{n,i}}{p_n})}=e^{\f{q_n}{p_n}w_{n,i}-\f{q_n}{2p_n^2}w_{n,i}^2+O(\f{w_{n,i}^3}{p_n^2})}\\
	&=e^{\f{q_n}{p_n}\sbr{ U-\sigma_{n,i}+\f{s_{n,i}}{\e_{n,i}} } -\f{q_n}{2p_n^2}\sbr{ U-\sigma_{n,i}+\f{s_{n,i}}{\e_{n,i}} }^2 +o(\f{1}{\e_{n,i}}) }\\
	&=e^{U-\sigma_{n,i}}\sbr{ 1+\f{s_{n,i}}{\e_{n,i}}+\f{1}{p_n}\sbr{ -\f{1}{2}U^2+(\sigma_{n,i}+\q_n)U-(\q_n\sigma_{n,i}+\f{1}{2}\sigma_{n,i}^2) }+o(\f{1}{\e_{n,i}}) }.
\end{aligned}$$
So
	$$e^{\sigma_{n,i}}\int_{B_{\e_{n,i}}(0)} \sbr{1+\f{z_{n,i}}{p_n}}^{p_n}\rd x =8\pi+\f{1}{\e_{n,i}}\int_{B_{\e_{n,i}}(0)} e^U s_{n,i}\rd x+o(\f{1}{ \e_{n,i} }),$$
which together with \eqref{tem-23-1} gives \eqref{est-2-1}.
\ep

Now we are ready to achieve our goal of this section.
\bp[Proof of Proposition \ref{sharpest-0}]
By contradiction, we assume that \eqref{keyest-1} does not hold for some $i\in\lbr{1,\cdots,k}$, i.e. \eqref{tem-101} holds for such $i$ up to a subsequence. Then by the above discussions, we obtain
$$\begin{aligned}
	&\quad \int_{B_{\f{r}{\mu_{n,i}}}(0)} \sbr{1+\f{z_{n,i}}{p_n}}^{p_n}\rd x-e^{\sigma_{n,i}}\int_{B_{\f{r}{\mu_{n,i}}}(0)} \sbr{1+\f{w_{n,i}}{p_n}}^{q_n}\rd x\\
	&=\f{1}{\e_{n,i}}\int_{B_{\e_{n,i}}(0)} e^U (t_{n,i}-s_{n,i})\rd x+o(\f{1}{\e_{n,i}}).
\end{aligned}$$
Then we deduce from \eqref{tem-17} that
	$$\int_{B_{\e_{n,i}}(0)} e^U (t_{n,i}-s_{n,i})\rd x=o(1)+O\sbr{\f{\e_{n,i}}{p_n}}=o(1).$$
However, using Lemma \ref{bound-0}, the Dominated Convergence Theorem and Lemma \ref{tem-31}, we get
	$$\int_{B_{\e_{n,i}}(0)} e^U (t_{n,i}-s_{n,i})\rd x\to -\int_{\R^2} e^U\phi_3 \rd x<0,$$
a contradiction. This finishes the proof.
\ep

\section{Sharp estimates of \texorpdfstring{$\mu_{n,i}$}{} and \texorpdfstring{$v_n(x_{n,i})$}{} }
In this section, we give sharp estimates of $\mu_{n,i}$, $v_n(x_{n,i})$  and  $u_n(x_{n,i})$, and finish the proof of Theorem \ref{thm1}. To this aim, we need better estimates of $w_{n,i},z_{n,i}$ than the forms in \eqref{expansion-w-1}.  In this section, for simplicity, we use the same notations as in Section 3 without introducing ambiguity.

First of all, we estimate the errors between $(w_{n,i},z_{n,i})$ and its limit $(U-\f{\q}{2},U)$. Let us define
\be\lab{function-stn}
	s_{n,i}:=p_n(w_{n,i}+\sigma_{n,i}-U)\quad\text{and}\quad t_{n,i}:=p_n(z_{n,i}-U). 	
\ee
where $\sigma_{n,i}=\q_n\ln v_n(x_{n,i})\to \frac{\theta}{2}$. Note that the difference between \eqref{function-stn} and that one \eqref{c3-6} in Section 3 is only the replacement of $\e_{n,i}$ with $p_n$. Thanks to Proposition \ref{sharpest-0}, we see that
	\be\label{c4-0}|s_{n,i}(0)|\le C,\quad t_{n,i}(0)=0,\quad \nabla t_{n,i}(0)=0,\quad\text{for any}~n\ge1.\ee
First we have

\bl\lab{limitsys-1}
For fixed $r\in(0,d_0)$, we have
\be\label{c4-1}\bcs
	-\Delta s_{n,i}=t_{n,i}(e^U+g_{n,i})+h_{n,i},\quad\text{in}~B_{\f{r}{\mu_{n,i}}}(0),\\
	-\Delta t_{n,i}=s_{n,i}(e^U+\wt g_{n,i})+\wt h_{n,i},\quad\text{in}~B_{\f{r}{\mu_{n,i}}}(0),
\ecs\ee
where
\be \label{c4-2}	
	g_{n,i}=O\sbr{\f{|z_{n,i}-U|}{1+|x|^{4-\delta}}}, \quad \wt g_{n,i}=O\sbr{\f{|w_{n,i}+\sigma_{n,i}-U|}{1+|x|^{4-\delta}}},\quad\text{in }B_{\f{r}{\mu_{n,i}}}(0),
\ee
\be \label{c4-3}	
	\max\lbr{|h_{n,i}|,|\wt h_{n,i}|}=O\sbr{\f{1}{1+|x|^{4-\delta}}}, \quad\text{in }B_{\f{r}{\mu_{n,i}}}(0),	
\ee
for any given small $\delta>0$, and as $n\to\iy$
\be \label{c4-4}	
	h_{n,i}(x)\to-\f{U^2}{2}e^U,\quad \wt h_{n,i}(x)\to\sbr{-\f{U^2}{2}+\f{3}{2}\q U-\f{5}{8}\theta^2}e^U,\quad\text{in }\mathcal{C}_{loc}(\mathbb{R}^2).
\ee
\el
\bp
By the same argument as Lemma \ref{limitsys-0}, we have that \eqref{c4-1} holds with the corresponding
\be \label{c4-7}	
g_{n,i}:=\f{p_n\sbr{e^{z_{n,i}}-e^U}-t_{n,i}e^U}{t_{n,i}}, \quad\quad h_{n,i}:=p_n\sbr{ \sbr{1+\f{z_{n,i}}{p_n}}^{p_n}-e^{z_{n,i}} },
\ee
\be \label{c4-8}	\wt g_{n,i}:=\f{p_n\sbr{ e^{w_{n,i}+\sigma_{n,i}}-e^U }-s_{n,i}e^U}{s_{n,i}},\quad 	\wt h_{n,i}:=p_ne^{\sigma_{n,i}}\sbr{ \sbr{1+\f{w_{n,i}}{p_n}}^{q_n}-e^{w_{n,i}} }, \ee
satisfying \eqref{c4-2}-\eqref{c4-3}.

To prove \eqref{c4-4}, we fix any $R>0$. Since $z_{n,i}(x)\to U$ and  $w_{n,i}(x)\to U-\frac\theta2$ uniformly in $B_R(0)$, i.e. $|z_{n,i}(x)/p_n|\leq \frac12$ uniformly in $B_R(0)$ for $n$ large, we can use\footnote{Remark that $\ln(1+s)=s-\frac12 s^2+O(s^3)$ can not hold uniformly for $|s|<1$, otherwise letting $s\to -1$ leads to a contradiction. This is the reason why we can not use the arguments \eqref{tem-12-4}-\eqref{tem-20-4} in $B_{\frac{r}{\mu_{n,i}}}(0)$ to obtain \eqref{c4-3}.} $\ln(1+s)=s-\frac12 s^2+O(s^3)$ for $|s|\leq \frac{1}{2}$
to obtain
\be\lab{tem-12-4} 	\begin{aligned}
h_{n,i}
&=p_ne^{z_{n,i}}\sbr{ e^{p_n\ln(1+\f{z_{n,i}}{p_n})-z_{n,i}} -1 }
	=p_ne^{z_{n,i}}\sbr{ e^{ -\f{z_{n,i}^2}{2p_n}+O(\f{1}{p_n^2}) } -1 }\\
&=p_ne^{z_{n,i}}\sbr{ -\f{z_{n,i}^2}{2p_n}+O(\f{1}{p_n^2})  }
	\to -\f{U^2}{2}e^U\quad\text{uniformly in $B_R(0)$}.
\end{aligned}\ee
Similarly,\be\lab{tem-20-4} 	
\begin{aligned}
	\wt h_{n,i}&=e^{w_{n,i}+\sigma_{n,i}}\sbr{ \q_nw_{n,i}-\f{1}{2}w_{n,i}^2+O(\f{1}{p_n})}
\to \sbr{ \q\sbr{U-\f{\q}{2}}-\f{1}{2}\sbr{U-\f{\q}{2}}^2 }e^U\\&=\sbr{-\f{1}{2}U^2+\f{3}{2}\q U-\f{5}{8}\q^2}e^U\quad\text{uniformly in $B_R(0)$}.
\end{aligned}\ee
The proof is complete.
\ep

\bl\lab{bound-1}
For fixed $r\in(0,d_0)$, $i=1,\cdots,k$ and any $\tau\in(0,1)$, there exists a $C_\tau>0$ such that
\be
	\max\lbr{|s_{n,i}(x)|,~|t_{n,i}(x)|}\le C_\tau(1+|x|)^\tau\quad\text{for any}~x\in B_{\f{r}{\mu_{n,i}}}(0).
\ee
\el
\bp
The proof is the same as Lemma \ref{bound-0} by just replacing $\e_{n,i}$ with $p_n$.
\ep

A direct consequence of Lemmas \ref{limitsys-1}-\ref{bound-1} is

\bl
For any $i=1,\cdots,k$, it holds
	$$(s_{n,i},t_{n,i})\to(s_{\iy,i},t_{\iy,i}),\quad\text{in}~\CR_{loc}^2(\R^2),$$
as $n\to\iy$, where $(s_{\iy,i},t_{\iy,i})$ is a solution of the system
\be\lab{tem-13} \bcs
	-\Delta s_{\iy,i}=t_{\iy,i}e^U-\f{1}{2}U^2e^U,\quad\text{in}~\R^2\\
	-\Delta t_{\iy,i}=s_{\iy,i}e^U+\sbr{-\f{U^2}{2}+\f{3}{2}\q U-\f{5}{8}\q^2}e^U,\quad\text{in}~\R^2,
\ecs \ee
and for any $\tau\in(0,1)$, there exists a $C_\tau>0$ such that
\be
	\max\lbr{|s_{\iy,i}(x)|,~|t_{\iy,i}(x)|}\le C_\tau(1+|x|)^\tau.
\ee
\el

Moreover, we can compute the exact formulas for the limit functions $(s_{\iy,i},t_{\iy,i})$.
\bl\lab{formula-st}
For any $i=1,\cdots,k$, there exist constants $c_{i,j}\in\R$, $j=0,3$ such that
\be\lab{function-s}
	s_{\iy,i}=-\f{3}{2}\q U+\f{5}{8}\q^2+\psi_0+c_{i,0}\phi_0 +c_{i,3}\phi_3,
\ee
and
\be\lab{function-t}
	t_{\iy,i}=-\f{3}{2}\q +\psi_0+c_{i,0}\phi_0 -c_{i,3}\phi_3,
\ee
where $\psi_0$ is a smooth radial solution of $-\Delta \psi-e^U\psi+\f{1}{2}U^2e^U=0$ such that $\psi_0(0)=\nabla\psi_0(0)=0$ and as $|x|\to\iy$,
\be \label{c4-16}	\psi_0(x)=12\ln|x|+O(1),\quad \nabla\psi_0(x)=12\f{x}{|x|^2}+O(\f{1}{|x|^2}). 	\ee
\el
\br
	From $t_{\iy,i}(0)=0$, we see that $c_{i,0}=c_{i,3}+\f{3}{2}\q$. In Remark \ref{remark-ci} below, we will give the exact values of $c_{i,0}$ and $c_{i,3}$.
\er
\bp[Proof of Lemma \ref{formula-st}]
Let
	$$(\bar s_{\iy,i},~\bar t_{\iy,i})=\sbr{ s_{\iy,i}+\f{3}{2}\q U-\f{5}{8}\q^2,~t_{\iy,i}+\f{3}{2}\q },$$
then from \eqref{tem-13} we see that
\be\lab{tem-14} \bcs
	-\Delta \bar s_{\iy,i}=\bar t_{\iy,i}e^U-\f{1}{2}U^2e^U,\quad\text{in}~\R^2,\\
	-\Delta \bar t_{\iy,i}=\bar s_{\iy,i}e^U-\f{1}{2}U^2e^U,\quad\text{in}~\R^2.
\ecs \ee
Setting $\vp_i:=\f{\bar s_{\iy,i}+\bar t_{\iy,i}}{2}$ and $\psi_i:=\f{\bar s_{\iy,i}-\bar t_{\iy,i}}{2}$, i.e.
	$$\bar s_{\iy,i}=\vp_i+\psi_i \quad\text{and}\quad \bar t_{\iy,i}=\vp_i-\psi_i.$$
From \eqref{tem-14}, we obtain $-\Delta \psi_i+e^U\psi_i=0$ in $\R^2$. Since $|\psi_i|\le C_\tau(1+|x|)^\tau$, we deduce from Lemma \ref{linear-1} that $\psi_i=c_{i,3}\phi_3$. Now we turn our attention to $\vp_i$. From \eqref{tem-14}, we obtain
\be\lab{tem-15}  	-\Delta \vp_i-e^U\vp_i=-\f{1}{2}U^2e^U\quad\text{in}~\R^2.  \ee
Recall Lemma \ref{linear-scalar} that $\phi_0(x)=\frac{8-|x|^2}{8+|x|^2}$ is a solution of $-\Delta u- e^U u=0$, i.e.
$\phi_0(r)=\f{8-r^2}{8+r^2}$ is a solution of $-u''-\f{1}{r}u'-\f{1}{(1+\f{1}{8}r^2)^2}u=0$. By ODE methods (see \cite[Section 2]{Chae=CMP} for detailed discussions), we obtain a radial special solution $\psi_0$ of \eqref{tem-15} as follows
\be
	\psi_0(x)=\f{8-|x|^2}{8+|x|^2} \sbr{ \int_0^{\f{|x|}{\sqrt8}}\f{\wt\psi_0(s)-\wt\psi_0(1)}{(1-s)^2}\rd s+\wt\psi_0(1)\f{|x|}{\sqrt8-|x|} },
\ee
where
	$$\wt\psi_0(r)=\f{16(1+r^2)^2}{r(1+r)^2}\int_0^r \f{s(1-s^2)}{(1+s^2)^3}\ln ^2(1+s^2)\rd s.$$
Obviously $\psi_0(0)=\nabla\psi_0(0)=0$. By direct computations, we obtain as $|x|\to\iy$,
	$$\psi_0(x)=-16\ln |x|\int_0^\iy\f{t(1-t^2)}{(1+t^2)^3}\ln^2(1+t^2)\rd t +O(1)=12\ln |x|+O(1),$$
and
	$$\nabla \psi_0(x)=-16\f{x}{|x|^2}\int_0^\iy\f{t(1-t^2)}{(1+t^2)^3}\ln^2(1+t^2)\rd t +O(\f{1}{|x|^2})=12\f{x}{|x|^2}+O(\f{1}{|x|^2}).$$
Since $|\vp_i|\le C_\tau(1+|x|)^\tau$, we deduce from Lemma \ref{linear-scalar} that
\be  \vp_i-\psi_0=\sum_{j=0}^2 c_{i,j}\phi_j,  \ee
where $c_{i,j}\in\R$ are constants. Therefore
	$$\bar s_{\iy,i}=\vp_i+\psi_i=\psi_0+\sum_{j=0}^2 c_{i,j}\phi_j+c_{i,3}\phi_3, $$
	$$\bar t_{\iy,i}=\vp_i-\psi_i=\psi_0+\sum_{j=0}^2 c_{i,j}\phi_j-c_{i,3}\phi_3, $$
which give \eqref{function-s} and \eqref{function-t}. Finally, $\nabla t_{\iy,i}(0)=0$ implies $c_{i,1}=c_{i,2}=0$.
\ep

Now we estimate the errors between $(s_{n,i},t_{n,i})$ and its limit $(s_{\iy,i},t_{\iy,i})$. For any $i=1\cdots,k$, let us define

\be\lab{function-st*} \bcs
	s_{n,i}^*:=\psi_0+l_{n,i}\phi_0+m_{n,i}\phi_3-(\q_n+\sigma_{n,i})U+\sigma_{n,i}\q_n+\f{1}{2}\sigma_{n,i}^2,\\
	t_{n,i}^*:=\psi_0+l_{n,i}\phi_0-m_{n,i}\phi_3-(\q_n+\sigma_{n,i}),
\ecs \ee
where
\be\lab{constant-lm} \bcs
	l_{n,i}:=\f{1}{2}(s_{n,i}(0)-(\sigma_{n,i}\q_n+\f{1}{2}\sigma_{n,i}^2)+\q_n+\sigma_{n,i}),\\
	m_{n,i}:=\f{1}{2}(s_{n,i}(0)-(\sigma_{n,i}\q_n+\f{1}{2}\sigma_{n,i}^2)-(\q_n+\sigma_{n,i})).
\ecs \ee
In Remark \ref{remark-ci}, we will give exact values of $l_{n,i}$ and $m_{n,i}$. By direct computations, we see that $l_{n,i}\to c_{i,0}=\frac12(s_{\infty, i}(0)-\frac58\theta^2+\frac32\theta)$ and $m_{n,i}\to c_{i,3}=c_{i,0}-\frac32\theta$. So $(s_{n,i}^*,t_{n,i}^*)\to(s_{\iy,i},t_{\iy,i})$. Furthermore, we have
\be\label{c4-23} \bcs
	-\Delta s_{n,i}^*=e^U t_{n,i}^*-\f{1}{2}U^2e^U,\quad\text{in}~\R^2,\\
	-\Delta t_{n,i}^*=e^U s_{n,i}^*+\sbr{ -\f{1}{2}(U-\sigma_{n,i})^2+\q_n(U-\sigma_{n,i}) }e^U,\quad\text{in}~\R^2.
\ecs \ee
Define
\be\label{c4-5}  	
	\al_{n,i}:= p_n(s_{n,i}-s_{n,i}^*)\quad\text{and}\quad	\beta_{n,i}:= p_n(t_{n,i}-t_{n,i}^*).
\ee
Then it follows from \eqref{constant-lm} and \eqref{c4-0} that
\be\al_{n,i}(0)=0,\quad \beta_{n,i}(0)=0,\quad \nabla\beta_{n,i}(0)=0 \quad\text{for any}~n\ge1.\ee
\begin{remark}\label{rmk-42}
Due to the complexity of the Lane-Emden system, one may see that we can not study directly the error terms $p_n(s_{n,i}-s_{\infty,i})$ and $p_n(t_{n,i}-t_{\infty,i})$, but have to  modify them with \eqref{c4-5}, i.e. we have to replace the limiting functions $(s_{\infty,i}, t_{\infty,i})$ with $(s_{n,i}^*, t_{n,i}^*)$ by choosing those approximate parameters $l_{n,i}, m_{n,i}$ carefully.
\end{remark}

\bl
	For fixed $r\in(0,d_0)$, we have
	\be\bcs
		-\Delta \al_{n,i}=\beta_{n,i}e^U+k_{n,i},\quad\text{in}~B_{\f{r}{\mu_{n,i}}}(0),\\
		-\Delta \beta_{n,i}=\al_{n,i}e^U+\wt k_{n,i},\quad\text{in}~B_{\f{r}{\mu_{n,i}}}(0),
	\ecs\ee
	where in $B_{\f{r}{\mu_{n,i}}}(0)$ it holds
	\be 	\max\lbr{|k_{n,i}|,|\wt k_{n,i}|}=O\sbr{ \f{1}{1+|x|^{4-\delta}} }, 	\ee
	for any given small $\delta>0$.
\el
\bp
By \eqref{c4-1} and \eqref{c4-23}, we have
\be\begin{aligned}
	-\Delta \al_{n,i}
	&=\beta_{n,i}e^U+p_n(t_{n,i}g_{n,i}+h_{n,i}+\frac12 U^2 e^U)\\
	&=:\beta_{n,i}e^U+k_{n,i}.
\end{aligned}\ee
Using \eqref{function-stn}, \eqref{c4-2} and Lemma \ref{bound-1}, we have
	$$ p_nt_{n,i}g_{n,i}=O\sbr{\frac{t_{n,i}^2}{1+|x|^{4-\frac{\delta}{2}}}} =O\sbr{ \f{1}{1+|x|^{4-\delta}} }.$$
To estimate the term $p_n(h_{n,i}+\frac12 U^2 e^U)$, we use
\eqref{c4-7} and \eqref{c3-15}-\eqref{c3-17} to deduce that
\begin{align*}
h_{n,i}&=p_n\sbr{ \sbr{1+\f{z_{n,i}}{p_n}}^{p_n}-e^{z_{n,i}} }
=p_n\sbr{ \sbr{1+\f{z_{n,i}}{p_n}}^{p_n}-\sbr{1+\f{y_{n,i}}{p_n}}^{p_n} }\\
&=p_n\sbr{1+\f{y_{n,i}}{p_n}}^{p_n-1}(z_{n,i}-y_{n,i})+\frac{p_n-1}{2}\sbr{1+\f{\xi}{p_n}}^{p_n-2}(z_n-y_n)^2\\
&=-\frac{1}{2}z_{n,i}^2e^{\frac{p_n-1}{p_n}z_{n,i}}+O(\frac{z_{n,i}^3}{p_n}e^{\frac{p_n-1}{p_n}z_{n,i}})+\sbr{1+\f{\xi}{p_n}}^{p_n-2}O(\frac{z_{n,i}^4}{p_n}),
\end{align*}
where $\xi$ is between $z_{n,i}$ and $y_{n,i}$. Since $(1+\frac{s}{p_n})^{p_n-2}$ is increasing for $s>-p_n$, we have
\[\sbr{1+\f{\xi}{p_n}}^{p_n-2}\leq\max\left\{\sbr{1+\f{z_{n,i}}{p_n}}^{p_n-2},e^{\frac{p_n-2}{p_n}z_{n,i}}\right\}
\leq C\left(\frac{C}{1+|x|^{4-\frac{\delta}{2}}}\right)^{\frac{p_n-2}{p_n}},\]
so we see from Lemma \ref{decay-12} that
\[O(\frac{z_{n,i}^3}{p_n}e^{\frac{p_n-1}{p_n}z_{n,i}})+\sbr{1+\f{\xi}{p_n}}^{p_n-2}O(\frac{z_{n,i}^4}{p_n})
=O\sbr{ \frac{1}{p_n}\f{1}{1+|x|^{4-\delta}} }.\]
Similarly, we have
\begin{align*}
U^2e^{U}-z_{n,i}^2e^{\frac{p_n-1}{p_n}z_{n,i}}
=&U^2e^{U}-z_{n,i}^2e^{z_{n,i}}+z_{n,i}^2e^{z_{n,i}}(1-e^{-\frac{z_{n,i}}{p_n}})\\
=&(\xi_1^2+2\xi_1)e^{\xi_1}(U-z_{n,i})+O(\frac{z_{n,i}^3}{p_n}e^{z_{n,i}})\\
=&\frac{-1}{p_n}(\xi_1^2+2\xi_1)e^{\xi_1} t_{n,i}+O(\frac{z_{n,i}^3}{p_n}e^{z_{n,i}})\\
=&O\sbr{ \frac{1}{p_n}\f{1}{1+|x|^{4-\delta}} },
\end{align*}
where $\xi_1$ is between $z_{n,i}$ and $U$, and Lemma \ref{bound-1} is used.
Combining the above estimates, we obtain
$p_n(h_{n,i}+\frac12 U^2 e^U)=O\sbr{\f{1}{1+|x|^{4-\delta}} }$
and so $|k_{n,i}|=O\sbr{ \f{1}{1+|x|^{4-\delta}} }$ for any $|x|\leq r/{\mu_{n,i}}$. Remark that this also improves \eqref{tem-12-4} as follows:
\begin{equation}\label{tem-12-4-1}
h_{n,i}\to -\frac12 U^2 e^U \quad\text{uniformly in }\; B_{\frac{r}{\mu_{n,i}}}(0).
\end{equation}

Again by \eqref{c4-1} and \eqref{c4-23}, we have
\be\begin{aligned}
	-\Delta \beta_{n,i}
	&=\alpha_{n,i}e^U
		+p_n\sbr{s_{n,i}\wt g_{n,i}+ \wt h_{n,i}+(\f{1}{2}(U-\sigma_{n,i})^2-\q_n(U-\sigma_{n,i}))e^U }\\
&=:\alpha_{n,i}e^U+\wt k_{n,i}.
\end{aligned}\ee
where $\wt h_{n,i}$ is given in Lemma \ref{limitsys-1}. By similar arguments as above, we can prove that
\[p_ns_{n,i}\wt g_{n,i}=O\sbr{ \f{1}{1+|x|^{4-\delta}} },\]
\[p_n\sbr{\wt h_{n,i}+(\f{1}{2}(U-\sigma_{n,i})^2-\q_n(U-\sigma_{n,i}))e^U }=O\sbr{ \f{1}{1+|x|^{4-\delta}} },\]
so $|\wt k_{n,i}|=O\sbr{ \f{1}{1+|x|^{4-\delta}} }$ for any $|x|\leq r/{\mu_{n,i}}$. The details are omitted here.
\ep

\bl\lab{bound-2}
For fixed $r\in(0,d_0)$, $i=1,\cdots,k$ and any $\tau\in(0,1)$, there exists a $C_\tau>0$ such that
	$$\max\lbr{|\al_{n,i}(x)|,~|\beta_{n,i}(x)|}\le C_\tau(1+|x|)^\tau\quad\text{for any}~x\in B_{\f{r}{\mu_{n,i}}}(0).$$
\el
\bp
The proof is similar to Lemma \ref{bound-0}, but we have to make some minor modifications. Again to simplify the notations, we omit the subscript $i$. Let
\be 	
	M_n:=\max_{|x|\le\f{r}{\mu_n}}\lbr{ \f{|\al_n(x)|}{(1+|x|)^\tau},~\f{|\beta_n(x)|}{(1+|x|)^\tau} },
\ee
and assume by contradiction that $M_n\to\iy$ as $n\to\iy$. We need to claim that
\be\label{c4-32}
	 M_n\to\iy\quad \text{implies} \quad\quad M_n^*=o(1)M_n,
\ee
where
	$$M_n^*:=\max_{|x|\le\f{r}{\mu_n}}\max_{|x'|=|x|}\lbr{ \f{|\al_n(x)-\al_n(x')|}{(1+|x|)^\tau},~\f{|\beta_n(x)-\beta_n(x')|}{(1+|x|)^\tau} }. $$
Again we prove this claim by contradiction, i.e. we suppose $M_n\to\iy$ and $M_n^*\ge CM_n$ for some $C>0$. Without loss of generality, we take $x_n'$ and $x_n''$ such that $|x_n'|=|x_n''|\le\f{r}{\mu_n}$ and
	$$M_n^*=\f{|\al_n(x_n')-\al_n(x_n'')|}{(1+|x_n'|)^\tau}.$$
By rotation, we assume that $x_n'=(x_{n,1}', x_{n,2}')$ with $x'_{n,2}>0$ and $x_n''=x_n'^-$, where $x^-=(x_1,x_2)^-=(x_1,-x_2)$. Set
\be 	
	\delta_n^*(x):=\al_n(x)-\al_n(x^-),\quad\quad \ga_n^*(x):=\beta_n(x)-\beta_n(x^-),
\ee
\[\delta_n(x):=\f{\delta_n^*(x)}{(1+x_2)^\tau}\quad\quad\ga_n(x):=\f{\ga_n^*(x)}{(1+x_2)^\tau},\]
\be 	
	M_n^{**}:=\max_{|x|\le\f{r}{\mu_n},x_2\ge0}\lbr{ |\delta_n(x)|,~|\ga_n(x)| },
\ee
and take $x_n^*\in \lbr{x\in\R^2:|x|\le\f{r}{\mu_n},x_2>0}$ such that $M_n^{**}=|\delta_n(x_n^*)|$. Then
	$$M_n^{**}\ge \f{|\delta_n^*(x_n')|}{(1+x_{n,2}')^\tau}\geq \f{|\delta_n^*(x_n')|}{(1+|x_n'|)^\tau}=M_n^*\ge C M_n\to +\infty.$$
Following the strategy of Lemma \ref{claim-1}, first we need to show
$|x_n^*|\le \f{r}{2\mu_n}$. Again assume by contradiction that $\f{r}{2\mu_n}\le |x_n^*|\le \f{r}{\mu_n}$. Here the key point is that since  \eqref{function-st*} says that $s_n^*(x)$ is radially symmetric, i.e. $s_n^*(x_n^*)=s_n^*(x_n^{*-})$, so
\begin{align*}\delta_n^*(x_n^*)=\al_n(x_n^*)-\al_n(x_n^{*-})=p_n\sbr{s_n(x_n^*)-s_n(x_n^{*-}) }=p_n^2\sbr{w_n(x_n^*)-w_n(x_n^{*-}) }.\end{align*}
Then
similarly as \eqref{tem-28}-\eqref{c3-28}, we have
\[
	\delta_n^*(x_n^*)
	=O(p_n^2\mu_nx_{n,2}^*)+o(p_n^2\mu_n)
\]
and so
\be 	M_n^{**}=\f{|\delta_n^*(x_n^*)|}{(1+x_{n,2}^*)^\tau}=O(p_n^{2}\mu_n(x_{n,2}^*)^{1-\tau})+o(p_n^{2}\mu_n)=o(1), \ee
which is a contradiction with $M_n^{**}\to\iy$. This proves $|x_n^*|\le \f{r}{2\mu_n}$. The a similar argument as Lemma \ref{claim-1} implies $M_n^*=o(1)M_n$, i.e. the claim \eqref{c4-32} holds.
Using this, one can prove Lemma \ref{bound-2} following the strategy of Lemma \ref{bound-0}, and we omit the details here.
\ep

Note from \eqref{function-stn} and \eqref{c4-5} that
\be\lab{expansion-w} 	w_{n,i}=U-\sigma_{n,i}+\f{s_{n,i}^*}{p_n}+\f{\al_{n,i}}{p_n^2}, \quad	z_{n,i}=U+\f{t_{n,i}^*}{p_n}+\f{\beta_{n,i}}{p_n^2}. \ee
Then we can give the following sharp estimates.
\bl
For fixed $r\in(0,d_0)$ and $i=1,\cdots,k$, it holds
\be\lab{est-1} 	
	\int_{B_{\f{r}{\mu_{n,i}}}(0)} \sbr{1+\f{z_{n,i}}{p_n}}^{p_n}\rd x =8\pi-\f{1}{p_n}\int_{\R^2}\Delta s_{n,i}^*\rd x+O(\f{1}{p_n^{2-\delta}}),
\ee
and
\be\lab{est-2} 	
	e^{\sigma_{n,i}}\int_{B_{\f{r}{\mu_{n,i}}}(0)} \sbr{1+\f{w_{n,i}}{p_n}}^{q_n}\rd x =8\pi-\f{1}{p_n}\int_{\R^2}\Delta t_{n,i}^*\rd x+O(\f{1}{p_n^{2-\delta}}),
\ee
for any given small $\delta>0$.
\el
\bp
The proof is similar as Lemma \ref{tem-32}. Noting that $p_n<\frac{r}{\mu_{n,i}}$ for $n$ large, we have
\be\lab{tem-22}
	\begin{aligned} 	
	\int_{B_{\f{r}{\mu_{n,i}}}(0)} \sbr{1+\f{z_{n,i}}{p_n}}^{p_n}\rd x
	&=\sbr{\int_{B_{p_n}(0)} + \int_{B_{\f{r}{\mu_{n,i}}}(0)\setminus B_{p_n}(0)} }\sbr{1+\f{z_{n,i}}{p_n}}^{p_n}\rd x\\
	&=\int_{B_{p_n}(0)} \sbr{1+\f{z_{n,i}}{p_n}}^{p_n}\rd x+O(\f{1}{p_n^{2-\delta}}).
	\end{aligned}
\ee
While in $B_{p_n}(0)$, we obtain from Lemma \ref{decay-12} that $z_{n,i}=O(\ln p_n)$, from \eqref{function-st*} that $t_{n,i}^*=O(\ln p_n)$ and from Lemma \ref{bound-2} that $\beta_{n,i}=O(p_n^{\tau})$, so we deduce from \eqref{expansion-w} that
\be\lab{tem-49}\begin{aligned}
	\sbr{1+\f{z_{n,i}}{p_n}}^{p_n}
	&=e^{p_n\ln(1+\f{z_{n,i}}{p_n})}=e^{z_{n,i}-\f{1}{2p_n}z_{n,i}^2+O(\f{z_{n,i}^3}{p_n^2})}\\
	&=e^{U+\f{t_{n,i}^*}{p_n}+\f{\beta_{n,i}}{p_n^2} -\f{1}{2p_n}\sbr{ U+\f{t_{n,i}^*}{p_n}+\f{\beta_{n,i}}{p_n^2} }^2 +O(\f{1}{p_n^{2-\delta}}) }\\
	&=e^U\sbr{ 1+\f{1}{p_n}(t_{n,i}^*-\f{1}{2}U^2)+O(\f{1}{p_n^{2-\delta}}) }.
\end{aligned}\ee
Since
	$$\int_{\R^2\setminus B_{p_n}(0)} e^U\sbr{ 1+\f{1}{p_n}(t_{n,i}^*-\f{1}{2}U^2) }\rd x=O\left(\int_{\R^2\setminus B_{p_n}(0)} \frac{(\ln|x|)^2}{|x|^4}\rd x\right)=O(\f{1}{p_n^{2-\delta}}),$$
we obtain that
\be  \begin{aligned}
\int_{B_{p_n}(0)} \sbr{1+\f{z_{n,i}}{p_n}}^{p_n}\rd x
&=\int_{\R^2}e^U\sbr{ 1+\f{1}{p_n}(t_{n,i}^*-\f{1}{2}U^2) }\rd x+O(\f{1}{p_n^{2-\delta}})\\
&=8\pi-\f{1}{p_n}\int_{\R^2}\Delta s_{n,i}^*\rd x+O(\f{1}{p_n^{2-\delta}}),
\end{aligned} \ee
which together with \eqref{tem-22} gives \eqref{est-1}.

Similarly, for $n$ large we have
\be\lab{tem-23}   	
	\int_{B_{\f{r}{\mu_{n,i}}}(0)} \sbr{1+\f{w_{n,i}}{p_n}}^{q_n}\rd x=\int_{B_{p_n}(0)} \sbr{1+\f{w_{n,i}}{p_n}}^{q_n}\rd x+O(\f{1}{p_n^{2-\delta}}).
\ee
Again in $B_{p_n}(0)$, we obtain from Lemma \ref{decay-12} that $w_{n,i}=O(\ln p_n)$, from \eqref{function-st*} that $s_{n,i}^*=O(\ln p_n)$ and from Lemma \ref{bound-2} that $\al_{n,i}=O(p_n^{\tau})$, so we deduce from \eqref{expansion-w} and $\frac{q_n}{p_n}=1+\frac{\theta_n}{p_n}=1+O(\frac{1}{p_n})$ that
$$ \begin{aligned}
	\sbr{1+\f{w_{n,i}}{p_n}}^{q_n} &=e^{q_n\ln(1+\f{w_{n,i}}{p_n})}=e^{\f{q_n}{p_n}w_{n,i}-\f{q_n}{2p_n^2}w_{n,i}^2+O(\f{w_{n,i}^3}{p_n^2})}\\
	&=e^{\f{q_n}{p_n}\sbr{ U-\sigma_{n,i}+\f{s_{n,i}^*}{p_n}+\f{\al_{n,i}}{p_n^2} }
		-\f{q_n}{2p_n^2}\sbr{ U-\sigma_{n,i}+\f{s_{n,i}^*}{p_n}+\f{\al_{n,i}}{p_n^2} }^2 +O(\f{1}{p_n^{2-\delta}}) }\\
	&=e^{U-\sigma_{n,i}}\sbr{ 1+\f{1}{p_n}\sbr{ s_{n,i}^*+\q_n(U-\sigma_{n,i})-\f{1}{2}(U-\sigma_{n,i})^2 }+O(\f{1}{p_n^{2-\delta}}) }.
\end{aligned}$$
Thus a similar argument implies
	$$e^{\sigma_{n,i}}\int_{B_{p_n}(0)} \sbr{1+\f{w_{n,i}}{p_n}}^{q_n}\rd x = 8\pi-\f{1}{p_n}\int_{\R^2}\Delta t_{n,i}^*\rd x+O(\f{1}{p_n^{2-\delta}}),$$
which together with \eqref{tem-23} gives \eqref{est-2}.
\ep

\bl\lab{tem-33}
For any $i=1,\cdots,k$, we have
	$$ v_n(x_{n,i})-\sqrt e=O(\f{\ln p_n}{p_n}). $$
\el
\bp
By \eqref{est-1} we have $\int_{B_{\f{r}{\mu_{n,i}}}(0)} \sbr{1+\f{z_{n,i}(y)}{p_n}}^{p_n}\rd y=8\pi+O(\frac{1}{p_n})$. From here and
 \eqref{tem-16}, we obtain
	$$1+O(\f{1}{p_n})=\f{1}{4\pi}\sbr{\f{\ln p_n}{p_n}+\f{p_n-1}{p_n}\ln v_n(x_{n,i})}(8\pi+O(\frac{1}{p_n})),$$
which implies $v_n(x_{n,i})=\sqrt e+O(\f{\ln p_n}{p_n})$.
\ep

Now we can prove the sharp estimates of $\mu_{n,i}$ and $v_n(x_{n,i})$.
\bl
For $i=1,\cdots,k$, we have
\be\lab{est-mu} 	
	\mu_{n,i}=e^{-\f{p_n}{4}}\sbr{ e^{-\sbr{ 2\pi\Phi_{k,i}(\boldsymbol{x}_n)+\f{3}{2}\ln 2+\f{3}{4}+\f{1}{8}\q_n }} +O(\f{1}{p_n^{1-\delta}}) },
\ee
and
\be\lab{est-v} 	
	v_n(x_{n,i})=\sqrt e\sbr{ 1-\f{\ln p_n}{p_n-1}+\f{1}{p_n}(4\pi\Phi_{k,i}(\boldsymbol{x}_n)+3\ln 2+2+\f{1}{4}\q_n) +O(\f{1}{p_n^{2-\delta}}) },
\ee
where $\Phi_{k,i}$ is defined by \eqref{function-Phik}, $\boldsymbol{x_n}=(x_{n,1},\cdots,x_{n,k})$ and $\delta>0$ is any fixed small constant.
\el
\bp
By taking $r\in(0,d_0)$ and using \eqref{c2-30}-\eqref{c2-31}, we have
\be\lab{tem-24} \begin{aligned}
	u_n(x_{n,i})
	&=\int_\Omega G(x_{n,i},y)v_n^{p_n}(y)\rd y=\sum_{j=1}^k \int_{B_r(x_{n,j})} G(x_{n,i},y)v_n^{p_n}(y)\rd y +O(\f{C^{p_n}}{p_n^{p_n}}).
\end{aligned}\ee
We compute each term in the right-hand side of \eqref{tem-24}. For $j=i$, we deduce from \eqref{c1-g} that
\be\begin{aligned}
	&\int_{B_r(x_{n,i})} G(x_{n,i},y)v_n^{p_n}(y)\rd y\\
	&=\f{v_n(x_{n,i})}{p_n}\int_{B_{\f{r}{\mu_{n,i}}}(0)}G(x_{n,i},x_{n,i}+\mu_{n,i}y)\sbr{1+\f{z_{n,i}(y)}{p_n}}^{p_n}\rd y\\
	&=-\f{v_n(x_{n,i})\ln \mu_{n,i}}{2\pi p_n}\int_{B_{\f{r}{\mu_{n,i}}}(0)}\sbr{1+\f{z_{n,i}(y)}{p_n}}^{p_n}\rd y\\
		&\quad -\f{v_n(x_{n,i})}{p_n}\int_{B_{\f{r}{\mu_{n,i}}}(0)}\sbr{ H(x_{n,i},x_{n,i}+\mu_{n,i}y)+\f{\ln|y|}{2\pi} }\sbr{1+\f{z_{n,i}(y)}{p_n}}^{p_n}\rd y\\
	&=:I_{n,1}+I_{n,2}.
\end{aligned}\ee
Using \eqref{est-1}, we obtain
	$$I_{n,1}=-\f{v_n(x_{n,i})\ln\mu_{n,i}}{2\pi p_n}\sbr{ 8\pi-\f{1}{p_n}\int_{\R^2}\Delta s_{n,i}^*\rd x+O(\f{1}{p_n^{2-\delta}}) }.$$
By the Taylor's expansion, \eqref{est-1} and Lemma \ref{decay-1}, we have
\be\label{c4-47}\begin{aligned}
	&\int_{B_{\f{r}{\mu_{n,i}}}(0)} H(x_{n,i},x_{n,i}+\mu_{n,i}y) \sbr{1+\f{z_{n,i}(y)}{p_n}}^{p_n}\rd y\\
	&=H(x_{n,i},x_{n,i})\int_{B_{\f{r}{\mu_{n,i}}}(0)}\sbr{1+\f{z_{n,i}(y)}{p_n}}^{p_n}\rd y+O(\mu_{n,i})\int_{B_{\f{r}{\mu_{n,i}}}(0)}|y|\sbr{1+\f{z_{n,i}(y)}{p_n}}^{p_n}\rd y\\
	&=H(x_{n,i},x_{n,i})\sbr{8\pi+O(\f{1}{p_n})}+O(\mu_{n,i})\\
	&=8\pi  H(x_{n,i},x_{n,i})+O(\f{1}{p_n}).
\end{aligned}\ee
Furthermore, similarly to the proof of \eqref{est-1}, we have
\be\label{c4-48}\begin{aligned}
	&\f{1}{2\pi}\int_{B_{\f{r}{\mu_{n,i}}}(0)}\ln|y|\sbr{1+\f{z_{n,i}(y)}{p_n}}^{p_n}\rd y\\
	&=\f{1}{2\pi}\int_{B_{p_n}(0)}\ln|y|\sbr{1+\f{z_{n,i}(y)}{p_n}}^{p_n}\rd y+O(\f{1}{p_n^{2-\delta}})\\
	&=\f{1}{2\pi}\int_{\R^2}\ln|y|e^U\rd y+\f{1}{2\pi p_n}\int_{\R^2}\ln|y|e^U(t_{n,i}^*-\f{1}{2}U^2)\rd y+O(\f{1}{p_n^{2-\delta}})\\
	&=\int_{0}^{\infty}\ln r\frac{r}{(1+\frac18 r^2)^2}dr+O(\f{1}{p_n})=6\ln 2+O(\f{1}{p_n}).
\end{aligned}\ee
Inserting these estimates into $I_{n,2}$ leads to
\[ 	I_{n,2}=-\f{v_n(x_{n,i})}{p_n}\sbr{ 8\pi  H(x_{n,i},x_{n,i})+6\ln 2+O(\f{1}{p_n}) },\]
and so
\be\lab{tem-110} \begin{aligned}	
	\int_{B_r(x_{n,i})} G(x_{n,i},y)v_n^{p_n}(y)\rd y
	&=-\f{v_n(x_{n,i})\ln\mu_{n,i}}{2\pi p_n}\sbr{ 8\pi-\f{1}{p_n}\int_{\R^2}\Delta s_{n,i}^*\rd x+O(\f{1}{p_n^{2-\delta}}) }\\
		&\quad -\f{v_n(x_{n,i})}{p_n}\sbr{ 8\pi  H(x_{n,i},x_{n,i})+6\ln 2+O(\f{1}{p_n}) }.
\end{aligned}\ee
While for $j\neq i$, note from Lemma \ref{tem-33} that $v_n(x_{n,j})=v_n(x_{n,i})+O(\frac{\ln p_n}{p_n})$, so it follows from \eqref{est-1} that
\be\lab{tem-111} \begin{aligned}
	&\int_{B_r(x_{n,j})} G(x_{n,i},y)v_n^{p_n}(y)\rd y\\
&=\f{v_n(x_{n,j})}{p_n}\int_{B_{\f{r}{\mu_{n,i}}}(0)}G(x_{n,i},x_{n,j}+\mu_{n,j}y)\sbr{1+\f{z_{n,j}(y)}{p_n}}^{p_n}\rd y\\
	&=\f{v_n(x_{n,j})}{p_n}G(x_{n,i},x_{n,j})\int_{B_{\f{r}{\mu_{n,i}}}(0)}\sbr{1+\f{z_{n,j}(y)}{p_n}}^{p_n}\rd y\\
&\quad+\f{v_n(x_{n,j})}{p_n}O(\mu_{n,j})\int_{B_{\f{r}{\mu_{n,i}}}(0)}|y|\sbr{1+\f{z_{n,j}(y)}{p_n}}^{p_n}\rd y\\ &=G(x_{n,i},x_{n,j})\f{v_n(x_{n,j})}{p_n}\sbr{8\pi+O(\f{1}{p_n})}+O(\f{\mu_{n,j}}{p_n})\\
	&= \f{v_n(x_{n,i})}{p_n} 8\pi G(x_{n,i},x_{n,j})+O(\f{1}{p_n^{2-\delta}}).
\end{aligned}\ee
Therefore, substituting \eqref{tem-110}-\eqref{tem-111} into \eqref{tem-24} leads to
\be\label{c4-50}\begin{aligned}
	\f{u_n(x_{n,i})}{v_n(x_{n,i})}
	&=-\f{\ln \mu_{n,i}}{2\pi p_n}\sbr{ 8\pi-\f{1}{p_n}\int_{\R^2}\Delta s_{n,i}^*\rd x+O(\f{1}{p_n^{2-\delta}}) }\\
		&\quad -\f{1}{p_n}\sbr{ 8\pi\Phi_{k,i}(\boldsymbol{x}_n)+6\ln 2 }+O(\f{1}{p_n^{2-\delta}}).
\end{aligned}\ee
On the other hand, Proposition \ref{sharpest-0} says $\f{u_n(x_{n,i})}{v_n(x_{n,i})}=1-\frac{1}{p_n}\sigma_{n,i}+O(\frac{1}{p_n^2})$. Inserting this into \eqref{c4-50} gives
\be\lab{tem-25} \begin{aligned}
	\ln \mu_{n,i}
	&=\f{ 1+\f{1}{p_n}\sbr{ 8\pi\Phi_{k,i}(\boldsymbol{x}_n)+6\ln 2-\sigma_{n,i} }+O(\f{1}{p_n^{2-\delta}}) }
		{ -\f{1}{2\pi p_n}\sbr{ 8\pi-\f{1}{p_n}\int_{\R^2}\Delta s_{n,i}^*\rd x+O(\f{1}{p_n^{2-\delta}}) }}\\
	&=-\f{p_n}{4}-\f{1}{4}\sbr{ 8\pi\Phi_{k,i}(\boldsymbol{x}_n)+6\ln 2-\sigma_{n,i}+\f{1}{8\pi}\int_{\R^2}\Delta s_{n,i}^*\rd x } +O(\f{1}{p_n^{1-\delta}}).
\end{aligned}\ee

Similarly, by the Green's representation formula we also have
\be\label{c4-53} \begin{aligned}
	v_n(x_{n,i})
	=\int_\Omega G(x_{n,i},y)u_n^{q_n}(y)\rd y
	=\sum_{j=1}^k \int_{B_r(x_{n,j})} G(x_{n,i},y)u_n^{q_n}(y)\rd y +O(\f{C^{q_n}}{p_n^{q_n}}).
\end{aligned}\ee
For $j=i$, we have
\be\begin{aligned}
	&\int_{B_r(x_{n,i})} G(x_{n,i},y)u_n^{q_n}(y)\rd y\\
	&=\f{v_n^{\q_n+1}(x_{n,i})}{p_n}\int_{B_{\f{r}{\mu_{n,i}}}(0)}G(x_{n,i},x_{n,i}+\mu_{n,i}y)\sbr{1+\f{w_{n,i}(y)}{p_n}}^{q_n}\rd y\\
	&=-\f{v_n^{\q_n+1}(x_{n,i})\ln \mu_{n,i}}{2\pi p_n}\int_{B_{\f{r}{\mu_{n,i}}}(0)}\sbr{1+\f{w_{n,i}(y)}{p_n}}^{q_n}\rd y\\
		&\quad -\f{v_n^{\q_n+1}(x_{n,i})}{p_n}\int_{B_{\f{r}{\mu_{n,i}}}(0)}\sbr{ H(x_{n,i},x_{n,i}+\mu_{n,i}y)+\f{\ln|y|}{2\pi} }\sbr{1+\f{w_{n,i}(y)}{p_n}}^{q_n}\rd y\\
	&=:I_{n,3}+I_{n,4}.
\end{aligned}\ee
Using \eqref{est-2} and $v_n^{\q_n}(x_{n,i})=e^{\sigma_{n,i}}$, we obtain
	$$I_{n,3}=-\f{v_n(x_{n,i})\ln\mu_{n,i}}{2\pi p_n}\sbr{ 8\pi-\f{1}{p_n}\int_{\R^2}\Delta t_{n,i}^*\rd x+O(\f{1}{p_n^{2-\delta}}) }.$$
Similarly to \eqref{c4-47}-\eqref{c4-48}, we have
\be\begin{aligned}
	&v_n^{\q_n}(x_{n,i})\int_{B_{\f{r}{\mu_{n,i}}}(0)}\sbr{ H(x_{n,i},x_{n,i}+\mu_{n,i}y)+\f{\ln|y|}{2\pi} }\sbr{1+\f{w_{n,i}(y)}{p_n}}^{q_n}\rd y\\	&=H(x_{n,i},x_{n,i})\sbr{8\pi+O(\f{1}{p_n})}+O(\mu_{n,i})+\f{e^{\sigma_{n,i}}}{2\pi}\int_{B_{\f{r}{\mu_{n,i}}}(0)}\ln|y|\sbr{1+\f{w_{n,i}(y)}{p_n}}^{q_n}\rd y\\
	&=8\pi  H(x_{n,i},x_{n,i})+6\ln 2+O(\f{1}{p_n}),
\end{aligned}\ee
thus
\be\label{c4-56} \begin{aligned}	
	\int_{B_r(x_{n,i})} G(x_{n,i},y)u_n^{q_n}(y)\rd y
	&=-\f{v_n(x_{n,i})\ln\mu_{n,i}}{2\pi p_n}\sbr{ 8\pi-\f{1}{p_n}\int_{\R^2}\Delta t_{n,i}^*\rd x+O(\f{1}{p_n^{2-\delta}}) }\\
		&\quad -\f{v_n(x_{n,i})}{p_n}\sbr{ 8\pi  H(x_{n,i},x_{n,i})+6\ln 2+O(\f{1}{p_n}) }.
\end{aligned}\ee For $j\neq i$, similarly to \eqref{tem-111}, we have
\be \label{c4-57}	\int_{B_r(x_{n,j})} G(x_{n,i},y)u_n^{q_n}(y)\rd y= \f{v_n(x_{n,i})}{p_n} 8\pi G(x_{n,i},x_{n,j})+O(\f{1}{p_n^{2-\delta}}).\ee
Inserting \eqref{c4-56}-\eqref{c4-57} into \eqref{c4-53} leads to
\be\begin{aligned}
	1&=-\f{\ln \mu_{n,i}}{2\pi p_n}\sbr{ 8\pi-\f{1}{p_n}\int_{\R^2}\Delta t_{n,i}^*\rd x+O(\f{1}{p_n^{2-\delta}}) }\\
		&\quad -\f{1}{p_n}\sbr{ 8\pi\Phi_{k,i}(\boldsymbol{x}_n)+6\ln 2 }+O(\f{1}{p_n^{2-\delta}}),
\end{aligned}\ee
which gives
\be\lab{tem-26} \begin{aligned}
	\ln \mu_{n,i}
	&=\f{ 1+\f{1}{p_n}\sbr{ 8\pi\Phi_{k,i}(\boldsymbol{x}_n)+6\ln 2 }+O(\f{1}{p_n^{2-\delta}}) }
		{ -\f{1}{2\pi p_n}\sbr{ 8\pi-\f{1}{p_n}\int_{\R^2}\Delta t_{n,i}^*\rd x+O(\f{1}{p_n^{2-\delta}}) }}\\
	&=-\f{p_n}{4}-\f{1}{4}\sbr{ 8\pi\Phi_{k,i}(\boldsymbol{x}_n)+6\ln 2+\f{1}{8\pi}\int_{\R^2}\Delta t_{n,i}^*\rd x } +O(\f{1}{p_n^{1-\delta}}).
\end{aligned}\ee

Comparing \eqref{tem-25} and \eqref{tem-26}, we obtain
	\be\label{c4-59}-\sigma_{n,i}+\f{1}{8\pi}\int_{\R^2}\Delta s_{n,i}^*\rd x=\f{1}{8\pi}\int_{\R^2}\Delta t_{n,i}^*\rd x+O(\f{1}{p_n^{1-\delta}}) .\ee
On the other hand, by \eqref{function-st*} we see that
	$$s_{n,i}^*-t_{n,i}^*=2m_{n,i}\phi_3-(\sigma_{n,i}+\q_n)(U-1)+(\sigma_{n,i}\q_n+\f{1}{2}\sigma_{n,i}^2). $$
Then as $|x|\to\iy$, it follows from Lemma \ref{linear} that
$$\begin{aligned}
	\nabla(s_{n,i}^*-t_{n,i}^*)
	&=2m_{n,i}C_0(\f{x}{|x|^2}+o(\frac{1}{|x|}))+(\sigma_{n,i}+\q_n)\f{4x}{8+|x|^2}\\
	&=(2m_{n,i}C_0+4(\sigma_{n,i}+\q_n))\f{x}{|x|^2}+o(\f{1}{|x|}),
\end{aligned}$$
where $C_0:=\pi^{-1}(e^{\f{\sqrt7}{2}\pi}+e^{-\f{\sqrt7}{2}\pi})$. So \eqref{c4-59} implies
$$\begin{aligned}
	8\pi\sigma_{n,i}+O(\f{1}{p_n^{1-\delta}})
	&=\int_{\R^2}\Delta (s_{n,i}^*-t_{n,i}^*)\rd x
		=\lim_{R\to\iy}\int_{\pa B_R(0)}\abr{ \nabla(s_{n,i}^*-t_{n,i}^*),\f{x}{R} }\rd \sigma_x\\
	&=2\pi(2m_{n,i}C_0+4(\sigma_{n,i}+\q_n)),
\end{aligned}$$
which means \be\label{c4-61} m_{n,i}=-2\q_nC_0^{-1}+O(\f{1}{p_n^{1-\delta}}).\ee From here, \eqref{function-st*}, \eqref{c4-16} and Lemma \ref{linear}, we get
	$$\nabla t_{n,i}^*=\nabla \psi_0+l_{n,i}\nabla \phi_0-m_{n,i}\nabla\phi_3 =\sbr{12+2\q_n+O(\f{1}{p_n^{1-\delta}})}\f{x}{|x|^2}+o(\f{1}{|x|}),\;\text{as } |x|\to\infty,$$
so
	\be\label{c4-80}\int_{\R^2}\Delta t_{n,i}^*\rd x=\lim_{R\to\iy}\int_{\pa B_R(0)}\abr{ \nabla t_{n,i}^*,\f{x}{R} }\rd \sigma_x=4\pi(6+\q_n)+O(\f{1}{p_n^{1-\delta}}).\ee
It follows from here and \eqref{tem-26} that
	$$\ln \mu_{n,i}=-\f{p_n}{4}-\sbr{ 2\pi\Phi_{k,i}(\boldsymbol{x}_n)+\f{3}{2}\ln 2+\f{3}{4}+\f{1}{8}\q_n } +O(\f{1}{p_n^{1-\delta}}),$$
which gives \eqref{est-mu}. Finally, since
	$$\ln\mu_{n,i}=-\f{1}{2}\ln p_n-\f{p_n-1}{2}\ln v_n(x_{n,i}),$$
we have
\[ \ln v_n(x_{n,i})=\f{1}{2}-\f{\ln p_n}{p_n-1}+\f{1}{p_n}\sbr{ 4\pi\Phi_{k,i}(\boldsymbol{x}_n)+3\ln 2+2+\f{1}{4}\q_n} +O(\f{1}{p_n^{2-\delta}}) \]
which gives \eqref{est-v}. The proof is complete.
\ep
\br\lab{remark-ci}
	In the above proof,
we see from \eqref{c4-59} and \eqref{c4-80} that
\be\label{c4-81}\int_{\R^2}\Delta s_{n,i}^*\rd x=8\pi (3+\theta_n)+O(\frac1{p_n^{1-\delta}}).\ee
Besides, we see from \eqref{c4-61} that	$$m_{n,i}=-2\q_nC_0^{-1}+O(\f{1}{p_n^{1-\delta}})=\f{-2\pi\q_n}{e^{\f{\sqrt7}{2}\pi}+e^{-\f{\sqrt7}{2}\pi}}+O(\f{1}{p_n^{1-\delta}}).$$
	Then we obtain from \eqref{constant-lm} that $$l_{n,i}=m_{n,i}+\q_n+\sigma_{n,i}=\f{-2\pi\q_n}{e^{\f{\sqrt7}{2}\pi}+e^{-\f{\sqrt7}{2}\pi}}+\q_n+\sigma_{n,i}+O(\f{1}{p_n^{1-\delta}}).$$
	By letting $n\to\iy$, we get
		$$c_{i,0}=\lim_{n\to\iy}l_{n,i} =\sbr{ \f{-2\pi}{e^{\f{\sqrt7}{2}\pi}+e^{-\f{\sqrt7}{2}\pi}}+\f{3}{2} }\q,\quad c_{i,3} =\lim_{n\to\iy}m_{n,i}=\f{-2\pi}{e^{\f{\sqrt7}{2}\pi}+e^{-\f{\sqrt7}{2}\pi}}\q.$$
\er

\vskip0.2in

\bp[Proof of Theorem \ref{thm1}] Note from \eqref{expansion-w} that the assertion (c) holds, so
it remains to compute the sharp estimate of $u_n(x_{n,i})$. By Proposition \ref{sharpest-0}, we have
$$\begin{aligned}
	u_n(x_{n,i})
	&=v_n(x_{n,i})-\f{1}{p_n}\theta_nv_n(x_{n,i})\ln v_n(x_{n,i})+O(\f{1}{p_n^2})\\
	&=\sqrt e\mbr{ 1 -\f{\ln p_n}{p_n-1} +\f{1}{p_n}\sbr{ 4\pi\Phi_{k,i}(\boldsymbol{x}_n)+3\ln2+2-\f{1}{4}\q_n }+O(\f{1}{p_n^{2-\delta}}) }.
\end{aligned}$$
Then we finish the proof.
\ep

\section{Nondegeneracy of \texorpdfstring{$1$}{}-bubble solutions}
In this section, we prove Theorem \ref{thm2} by contradiction. Let $(u_n,v_n)$ be a $1$-bubble solution concentrating at $x_\iy$, i.e. $k=i=1$ in the previous sections, so from now on we omit the subscript $i$. Assume that $x_\iy$ is a nondegenerate critical point of the Robin function $R(x)$. To prove that $(u_n, v_n)$ is nondegenerate for large $n$, we suppose by contradiction that there exists a sequence $(\xi_n,\eta_n)\in H_0^1(\Omega)$ of solutions to the linearized system \eqref{linear-0}, satisfying
	$$\nm{\xi_n}_{L^\iy(\Omega)}=\max\lbr{\nm{\xi_n}_{L^\iy(\Omega)},~\nm{\eta_n}_{L^\iy(\Omega)}}=1,\quad\text{for any}~n\ge1.$$
Let $d_0$ be defined by \eqref{constant-d0}, $x_n$ be the local maximum points $v_n(x_n)=\max_{B_{d_0}(x_\iy)} v_n$,
the scaling parameter be $\mu_n=\sbr{p_n v_n^{p_n-1}(x_n)}^{-\f{1}{2}}$ and the corresponding scaling functions $(w_n,z_n)$ be defined as \eqref{function-wz}.
Define
\be\lab{function-wtxi}
	\wt\xi_n (x):=\xi_n(x_n+\mu_nx)\quad\text{and}\quad \wt\eta_n (x):=\eta_n(x_n+\mu_nx) \quad\text{for}\quad x\in\Omega_n:=\f{\Omega-x_n}{\mu_n}.
\ee
By direct computations, we know that
\be\label{c5-2} \bcs
	-\Delta\wt\xi_n=\sbr{1+\f{z_n}{p_n}}^{p_n-1}\wt\eta_n,\quad\text{in}~\Omega_n,\\
	-\Delta\wt\eta_n=\f{q_n}{p_n}v_n^{\q_n}(x_n)\sbr{1+\f{w_n}{p_n}}^{q_n-1}\wt\xi_n,\quad\text{in}~\Omega_n,\\
	\wt\xi_n=\wt\eta_n=0,\quad\text{on}~\pa\Omega_n\\
\|\wt \eta_n\|_{L^\infty}\leq \|\wt \xi_n\|_{L^\infty}=1.
\ecs \ee
Then by standard elliptic estimates, we obtain that up to a subsequence, $(\wt\xi_n,\wt\eta_n)\to(\wt\xi_\iy,\wt\eta_\iy)$ in $\CR_{loc}^2(\R^2)$, and $(\wt\xi_\iy,\wt\eta_\iy)$ solves the linearized system \eqref{linearsys-0}. Then by Lemma \ref{linear} and $\|\wt\xi_\iy\|_{L^\infty}, \|\wt \eta_\iy\|_{L^\infty}\leq 1$, it holds
\be\lab{tem-51}
	(\wt\xi_\iy,\wt\eta_\iy)=\sum_{j=0}^2 a_j(\phi_j,\phi_j),
\ee
where $a_j\in\R$, $j=0,1,2$ are constants. We want to prove that $a_j=0$ for all $j$ and obtain a contradiction. First, we need the following lemma.
\bl\lab{lemma-out}
	For fixed $r\in(0,d_0)$, it holds
	\be\label{c5-4} \bcs
		\xi_n(x)=A_{n,0}G(x_n,x)+\sum_{j=1}^2A_{n,j}\pa_jG(x_n,x)+o(\mu_n),\\
		\eta_n(x)=\wt A_{n,0}G(x_n,x)+\sum_{j=1}^2\wt A_{n,j}\pa_jG(x_n,x)+o(\mu_n),\\
	\ecs \ee
	in $\CR^1(\Omega\setminus B_{2r}(x_n))$. Moreover,
		$$A_{n,0},~\wt A_{n,0}=\f{8\pi}{p_n}a_0+o(\f{1}{p_n}),$$
		\be\label{c5-5} A_{n,j},~\wt A_{n,j}=2\pi a_j\mu_n+o(\mu_n),\quad j=1,2.\ee
\el
\bp
For any $x\in\Omega\setminus B_{2r}(x_n)$, by the Green's representation formula and \eqref{c2-31} that $v_n=O(\f{1}{p_n})$ in $\Omega\setminus B_r(x_n)$, we have
$$\begin{aligned}
	\xi_n(x)
	&=\int_\Omega G(y,x)p_nv_n^{p_n-1}(y)\eta_n(y)\rd y\\
	&=\sbr{ \int_{B_r(x_n)}+ \int_{\Omega\setminus B_r(x_n)} } G(y,x)p_nv_n^{p_n-1}(y)\eta_n(y)\rd y\\
&=\int_{B_r(x_n)} G(y,x)p_nv_n^{p_n-1}(y)\eta_n(y)\rd y +O\left(p_n(\frac{C}{p_n})^{p_n-1}\right)\\
	&=\int_{B_r(x_n)} G(y,x)p_nv_n^{p_n-1}(y)\eta_n(y)\rd y +o(\mu_n).
\end{aligned}$$
Let
\be\label{c05-5}
	A_{n,0}:=\int_{B_r(x_n)}p_nv_n^{p_n-1}\eta_n\rd y \quad\text{and}\quad  A_{n,j}:=\int_{B_r(x_n)}(y-x_n)_jp_nv_n^{p_n-1}\eta_n\rd y,
\ee
for $j=1,2$. Since $\nabla^2 G(y,x)=O(1)$ for $x\in\Omega\setminus B_{2r}(x_n)$ and $y\in B_r(x_n)$,  we get by the Taylor's expansion that
$$\begin{aligned}
	\int_{B_r(x_n)} G(y,x)p_nv_n^{p_n-1}(y)\eta_n(y)\rd y
	&=A_{n,0}G(x_n,x)+\sum_{j=1}^2A_{n,j}\pa_jG(x_n,x)\\
	&\quad+\int_{B_r(x_n)} O(|y-x_n|^2)p_nv_n^{p_n-1}(y)\eta_n(y)\rd y.
\end{aligned}$$
Since Remark \ref{decay-10} yields
	$$\int_{B_{r}(x_{n})} |y-x_{n}|^2p_nv_n^{p_n-1}(y)\eta_n(y)\rd y=O(\mu_n^2)\int_{B_{\f{r}{\mu_{n}}}(0)}|y|^2\sbr{1+\f{z_{n,i}(y)}{p_n}}^{p_n-1}\rd y=o(\mu_n),$$
we obtain
\be\lab{tem-36}
	\xi_n(x)=A_{n,0}G(x_n,x)+\sum_{j=1}^2A_{n,j}\pa_jG(x_n,x)+o(\mu_n),
\ee
for $x\in\Omega\setminus B_{2r}(x_n)$. Similarly, for $h=1,2$, it holds
$$\begin{aligned}
	\pa_h\xi_n(x)
	&=\int_\Omega D_hG(y,x)p_nv_n^{p_n-1}(y)\eta_n(y)\rd y\\
	&=\int_{B_r(x_n)} D_hG(y,x)p_nv_n^{p_n-1}(y)\eta_n(y)\rd y +o(\mu_n).
\end{aligned}$$
then one can get
	$$\pa_h\xi_n(x)=A_{n,0}D_hG(x_n,x)+\sum_{j=1}^2A_{n,j}\pa_jD_hG(x_n,x)+o(\mu_n),$$
which means that \eqref{tem-36} holds in $\CR^1(\Omega\setminus B_{2r}(x_n))$.

Similarly, by letting
\be\label{c05-7}
	\wt A_{n,0}:=\int_{B_r(x_n)}q_nu_n^{q_n-1}\xi_n\rd y \quad\text{and}\quad \wt A_{n,j}:=\int_{B_r(x_n)}(y-x_n)_jq_nu_n^{q_n-1}\xi_n\rd y,
\ee
for $j=1,2$, we can also obtain
\be
	\eta_n(x)=\wt A_{n,0}G(x_n,x)+\sum_{j=1}^2\wt A_{n,j}\pa_jG(x_n,x)+o(\mu_n),
\ee
in $\CR^1(\Omega\setminus B_{2r}(x_n))$. This proves the first part of Lemma \ref{lemma-out}.

Now we compute the estimates of $A_{n,j}$ and $\wt A_{n,j}$ for $j=1,2$. By \eqref{c05-5}, Remark \ref{decay-10}, the Dominated Convergence Theorem and \eqref{tem-51} we see that
\be\label{c5-91}\begin{aligned}
	A_{n,j}
	&=\mu_n\int_{B_{\f{r}{\mu_n}}(0)} y_j\sbr{1+\f{z_n(y)}{p_n}}^{p_n-1}\wt\eta_n(y)\rd y\\
	&=\mu_n\int_{\R^2}y_je^{U(y)}\wt\eta_\iy(y)\rd y +o(\mu_n)\\
&=\mu_n\sum_{i=0}^2a_i\int_{\R^2}y_je^{U(y)}\phi_i(y)\rd y +o(\mu_n)\\
	&=\mu_na_j\int_{\R^2}y_je^{U(y)}\phi_j(y)\rd y +o(\mu_n)\\
	&=2\pi a_j\mu_n+o(\mu_n),\qquad j=1,2,
\end{aligned}\ee
where we have used $\int_{\R^2}y_je^{U(y)}\phi_i(y)\rd y=0$ for $i\neq j$ by the symmetry
and
\[\int_{\R^2}y_je^{U(y)}\phi_j(y)\rd y=\int_{\R^2}\frac{64y_j^2}{(8+|y|^2)^3}\rd y=\int_{\R^2}\frac{32|y|^2}{(8+|y|^2)^3}\rd y=2\pi.\]
Similarly,
$$\begin{aligned}
	\wt A_{n,j}
	&=\f{q_n}{p_n}\mu_nv_n^{\q_n}(x_n)\int_{B_{\f{r}{\mu_n}}(0)} y_j\sbr{1+\f{w_n(y)}{p_n}}^{q_n-1}\wt\xi_n(y)\rd y\\
	&=\mu_ne^{\f{\q}{2}}\int_{\R^2}y_je^{U(y)-\f{\q}{2}}\wt\xi_\iy(y)\rd y +o(\mu_n)\\
	&=\mu_n\sum_{i=0}^2a_i\int_{\R^2}y_je^{U(y)}\phi_i(y)\rd y +o(\mu_n)\\
	&=2\pi a_j\mu_n+o(\mu_n).
\end{aligned}$$

Finally, we compute the estimates of $A_{n,0}$ and $\wt A_{n,0}$. Observe from \eqref{c05-5} and \eqref{c05-7} that
\be\lab{tem-40}\begin{aligned}
	A_{n,0}
	&=\int_{B_r(x_n)}\sbr{1-\f{v_n}{v_n(x_n)}}p_nv_n^{p_n-1}\eta_n\rd x +\f{1}{v_n(x_n)}\int_{B_r(x_n)} p_nv_n^{p_n}\eta_n\rd x,
\end{aligned}\ee
and
\be\lab{tem-41}\begin{aligned}
	\wt A_{n,0}
	&=\int_{B_r(x_n)}\sbr{1-\f{u_n}{v_n(x_n)}}q_nu_n^{q_n-1}\xi_n\rd x +\f{1}{v_n(x_n)}\int_{B_r(x_n)} q_nu_n^{q_n}\xi_n\rd x.
\end{aligned}\ee
On the one hand, by the Dominated Convergence Theorem we see that
\be\lab{tem-42}\begin{aligned}
	&\quad \int_{B_r(x_n)}\sbr{1-\f{v_n}{v_n(x_n)}}p_nv_n^{p_n-1}\eta_n\rd x \\
	&=\f{1}{p_n}\int_{B_{\f{r}{\mu_n}}(0)}-z_n(y)\sbr{1+\f{z_n(y)}{p_n}}^{p_n-1}\wt\eta_n(y)\rd y\\
&=-\f{1}{p_n}\sum_{i=0}^2a_i\int_{\R^2}U(y)e^{U(y)}\phi_i(y)\rd y +o(\f{1}{p_n})\\
	&=\f{8\pi}{p_n}a_0+o(\f{1}{p_n}),
\end{aligned}\ee
where we have used $\int_{\R^2}U(y)e^{U(y)}\phi_i(y)\rd y=0$ for $i\in \{1,2\}$ by the symmetry
and
\[\int_{\R^2}U(y)e^{U(y)}\phi_0(y)\rd y=\int_{\R^2}-2\ln(1+\frac{|y|^2}{8})\frac{1}{(1+\frac{|y|^2}{8})^2}\frac{8-|y|^2}{8+|y|^2}\rd y=-8\pi.\]
Similarly,
\be\lab{tem-43}\begin{aligned}
	&\quad \int_{B_r(x_n)}\sbr{1-\f{u_n}{v_n(x_n)}}q_nu_n^{q_n-1}(x)\xi_n\rd x \\
	&=\f{q_n}{p_n^2}v_n^{\q_n}(x_n)\int_{B_{\f{r}{\mu_n}}(0)}-w_n(y)\sbr{1+\f{w_n(y)}{p_n}}^{q_n-1}\wt\xi_n(y)\rd y\\
	&=-\f{1}{p_n}\sum_{i=0}^2a_i\int_{\R^2}(U(y)-\f{\q}{2})e^{U(y)}\phi_i(y)\rd y +o(\f{1}{p_n})\\
	&=\f{8\pi}{p_n}a_0+o(\f{1}{p_n}),
\end{aligned}\ee
where we have used
\be\label{c5-93}\int_{\R^2}e^{U(y)}\phi_0(y)\rd y=\int_{\R^2}\frac{1}{(1+\frac{|y|^2}{8})^2}\frac{8-|y|^2}{8+|y|^2}\rd y=0.\ee

Substituting \eqref{tem-42}-\eqref{tem-43} into \eqref{tem-40}-\eqref{tem-41} and using \eqref{c2-25}-\eqref{c2-26}, we get
\be\lab{tem-44}
	A_{n,0}=\f{8\pi}{p_n}a_0+o(\f{1}{p_n})+\f{p_n}{v_n(x_n)}\int_{B_r(x_n)} v_n^{p_n}\eta_n\rd x=O(1),
\ee
\be\lab{tem-45}
	\wt A_{n,0}=\f{8\pi}{p_n}a_0+o(\f{1}{p_n})+\f{q_n}{v_n(x_n)}\int_{B_r(x_n)} u_n^{q_n}\xi_n\rd x=O(1).
\ee
On the other hand, to give sharper estimates of $\int_{B_r(x_n)} v_n^{p_n}\eta_n\rd x$ and $\int_{B_r(x_n)} u_n^{q_n}\xi_n\rd x$, we use the equations satisfied by $(u_n,v_n)$ and $(\xi_n,\eta_n)$ and obtain
\be\lab{tem-74}\begin{aligned}
	&\quad p_n\int_{B_{2r}(x_n)} v_n^{p_n}\eta_n\rd x=\int_{B_{2r}(x_n)} -\Delta\xi_nv_n\rd x\\
	&=\int_{B_{2r}(x_n)} \Delta v_n\xi_n-\Delta\xi_nv_n\rd x+\int_{B_{2r}(x_n)} u_n^{q_n}\xi_n\rd x\\
	&=\int_{\pa B_{2r}(x_n)}\abr{\nabla v_n,\nu}\xi_n-\abr{\nabla\xi_n,\nu}v_n\rd\sigma_x+\int_{B_{2r}(x_n)} u_n^{q_n}\xi_n\rd x,
\end{aligned}\ee
and
\be\lab{tem-75}\begin{aligned}
	&\quad q_n\int_{B_{2r}(x_n)} u_n^{q_n}\xi_n\rd x=\int_{B_{2r}(x_n)} -\Delta\eta_nu_n\rd x\\
	&=\int_{B_{2r}(x_n)} \Delta u_n\eta_n-\Delta\eta_nu_n\rd x+\int_{B_{2r}(x_n)} v_n^{p_n}\eta_n\rd x\\
	&=\int_{\pa B_{2r}(x_n)}\abr{\nabla u_n,\nu}\eta_n-\abr{\nabla\eta_n,\nu}u_n\rd\sigma_x+\int_{B_{2r}(x_n)} v_n^{p_n}\eta_n\rd x.
\end{aligned}\ee
Since Lemma \ref{expansion-1} and \eqref{2-27}-\eqref{2-28} say that
\[\bcs
		u_n(x)=C_{n}G(x_{n},x)+o(\f{\mu_n}{p_n}),\\
		v_n(x)=\wt C_{n} G(x_{n},x)+o(\f{\mu_n}{p_n}),
	\ecs\qquad \text{in}~\CR^1(\Omega\setminus\cup_{i=1}^kB_{2r}(x_{n,i})),
	\]
with $C_n, \wt C_{n}=O(\frac{1}{p_n})$, we deduce from \eqref{c5-4}, \eqref{c5-5} and \eqref{tem-44}-\eqref{tem-45} that
\[
	\int_{\pa B_{2r}(x_n)}\abr{\nabla v_n,\nu}\xi_n-\abr{\nabla\xi_n,\nu}v_n\rd\sigma_x=O\sbr{\f{ A_{n,1}+ A_{n,2}}{p_n}}+o\sbr{\f{ A_{n,0}}{p_n}\mu_n}=O(\f{\mu_n}{p_n}),
\]
\[
	\int_{\pa B_{2r}(x_n)}\abr{\nabla u_n,\nu}\eta_n-\abr{\nabla\eta_n,\nu}u_n\rd\sigma_x=O\sbr{\f{\wt A_{n,1}+\wt A_{n,2}}{p_n}}+o\sbr{\f{\wt A_{n,0}}{p_n}\mu_n}=O(\f{\mu_n}{p_n}).
\]
Inserting these into \eqref{tem-74}-\eqref{tem-75} leads to
	$$p_n\int_{B_{2r}(x_n)} v_n^{p_n}\eta_n\rd x=\int_{B_{2r}(x_n)} u_n^{q_n}\xi_n\rd x+O(\f{\mu_n}{p_n}),$$
	$$q_n\int_{B_{2r}(x_n)} u_n^{q_n}\xi_n\rd x=\int_{B_{2r}(x_n)} v_n^{p_n}\eta_n\rd x+O(\f{\mu_n}{p_n}),$$
which means
	$$\int_{B_{2r}(x_n)} v_n^{p_n}\eta_n\rd x=O(\f{\mu_n}{p_n^2}) \quad\text{and}\quad \int_{B_{2r}(x_n)} u_n^{q_n}\xi_n\rd x=O(\f{\mu_n}{p_n^2}).$$
Since $u_n,v_n=O(\f{1}{p_n})$ in $B_{2r}(x_n)\setminus B_r(x_n)$, we immediately obtain
\be\lab{tem-46}
	\int_{B_{r}(x_n)} v_n^{p_n}\eta_n\rd x=O(\f{\mu_n}{p_n^2})=o(\f{1}{p_n^2}) \quad\text{and}\quad \int_{B_{r}(x_n)} u_n^{q_n}\xi_n\rd x=O(\f{\mu_n}{p_n^2})=o(\f{1}{p_n^2}).
\ee
Inserting \eqref{tem-46} into \eqref{tem-44}-\eqref{tem-45}, we obtain
	$A_{n,0},~\wt A_{n,0}=\f{8\pi}{p_n}a_0+o(\f{1}{p_n}).$
This completes the proof.
\ep

Lemma \ref{lemma-out} implies that $\xi_n,\eta_n$ are actually small at points far away from the blow up point $x_n$.
\bc\lab{out-1}
	For fixed $r\in(0,d_0)$, it holds
	$$\nm{\xi_n}_{L^\iy(\Omega\setminus B_r(x_n))}=o(1),\quad \nm{\eta_n}_{L^\iy(\Omega\setminus B_r(x_n))}=o(1).$$
\ec

Further, we have the following result. Note that $\mu_np_n\to 0$.
\bl\lab{out-2}
	For $n$ large, it holds
	\be\lab{tem-54}
		\nm{\xi_n}_{L^\iy(\Omega\setminus B_{2\mu_np_n}(x_n))}=O\sbr{p_nA_{n,0}}+O(\f{1}{p_n^{1-\delta}}),
	\ee
	\be\lab{tem-55}
		\nm{\eta_n}_{L^\iy(\Omega\setminus B_{2\mu_np_n}(x_n))}=O\sbr{p_n\wt A_{n,0}}+O(\f{1}{p_n^{1-\delta}}),
	\ee
	where $\delta>0$ is small.
\el
\bp
For any $x\in\Omega\setminus B_{2\mu_np_n}(x_n)$, we have
\be\lab{tem-102}\begin{aligned}
	\xi_n(x)
	&=\int_\Omega G(x,y)p_nv_n^{p_n-1}(y)\eta_n(y)\rd y\\
	&=\int_{B_{\mu_np_n}(x_n)}G(x,y)p_nv_n^{p_n-1}(y)\eta_n(y)\rd y\\
		&\quad +\int_{B_r(x_n)\setminus B_{\mu_np_n}(x_n)}G(x,y)p_nv_n^{p_n-1}(y)\eta_n(y)\rd y+O(\f{C^{p_n-1}}{p_n^{p_n-2}})\\
	&=:I_{n,1}+I_{n,2}+O(\f{C^{p_n-1}}{p_n^{p_n-2}}),
\end{aligned}\ee
where $r\in(0,d_0)$. For any $x\in\Omega\setminus B_{2\mu_np_n}(x_n)$ and $y\in B_{\mu_np_n}(x_n)$, it holds
	$$|D_yG(x,y)|=O(\f{1}{|x-y|})=O(\f{1}{\mu_np_n}),$$
then by the mean value theorem and $\|\eta_n\|_{L^\infty}\leq 1$ we get
$$\begin{aligned}
	I_{n,1}
	&=G(x,x_n)\int_{B_{\mu_np_n}(x_n)} p_nv_n^{p_n-1}(y)\eta_n(y)\rd y\\
		&\quad +O(\f{1}{\mu_np_n})\int_{B_{\mu_np_n}(x_n)} p_n|y-x_n|v_n^{p_n-1}(y) \rd y.
\end{aligned}$$
Form the definition \eqref{c05-5} of $A_{n,0}$ and Remark \ref{decay-10}, we see that
$$\begin{aligned}
	\int_{B_{\mu_np_n}(x_n)} p_nv_n^{p_n-1}(y)\eta_n(y)\rd y
	&=A_{n,0}-\int_{B_r(x_n)\setminus B_{\mu_np_n}(x_n)} p_nv_n^{p_n-1}(y)\eta_n(y)\rd y\\
&=A_{n,0}-\int_{B_{\frac{r}{\mu_n}}(0)\setminus B_{p_n}(0)} \left(1+\frac{z_n(y)}{p_n}\right)^{p_n-1}\widetilde{\eta}_n(y)\rd y\\
	&=A_{n,0}+O(\f{1}{p_n^{2-\delta}}),
\end{aligned}$$
and
$$\int_{B_{\mu_np_n}(x_n)} p_n|y-x_n|v_n^{p_n-1}(y) \rd y=\mu_n\int_{B_{p_n}(0)}|y|\sbr{1+\f{z_n(y)}{p_n}}^{p_n-1}\rd y=O(\mu_n).$$
Since $|x-x_n|\geq 2\mu_np_n$ implies
	$G(x,x_n)=O(|\ln(\mu_np_n)|)=O(p_n),$
we obtain
	$$I_{n,1}=O\sbr{p_nA_{n,0}}+O(\f{1}{p_n^{1-\delta}}).$$
Again by Remark \ref{decay-10} and $\ln|\mu_n|=O(p_n)$,
$$\begin{aligned}
	|I_{n,2}|
	&\le C\int_{B_r(x_n)\setminus B_{\mu_np_n}(x_n)}|\ln|x-y||p_nv_n^{p_n-1}(y)\rd y\\
&=C \int_{B_{\frac{r}{\mu_n}}(0)\setminus B_{p_n}(0)} |\ln|\mu_n(y_n-y)||\left(1+\frac{z_n(y)}{p_n}\right)^{p_n-1}\rd y\\
	&\le C\int_{B_{\f{r}{\mu_n}}(0)\setminus B_{p_n}(0)}\f{|\ln\mu_n|+|\ln|y-y_n||}{1+|y|^{4-\delta}}\rd y\\
&=O(\f{1}{p_n^{1-\delta}})+C\int_{B_{\f{r}{\mu_n}}(0)\setminus B_{p_n}(0)}\f{|\ln|y-y_n||}{1+|y|^{4-\delta}}\rd y,
\end{aligned}$$
where $y_n:=\f{x-x_n}{\mu_n}$. Since it was proved in \cite[(4.40)]{LE-1} that $\int_{B_{\f{r}{\mu_n}}(0)\setminus B_{p_n}(0)}\f{|\ln|y-y_n||}{1+|y|^{4-\delta}}\rd y=O(\f{1}{p_n^{1-\delta}})$, so $I_{n,2}=O(\f{1}{p_n^{1-\delta}})$. Substituting the estimates of $I_{n,1}$ and $I_{n,2}$ into \eqref{tem-102}, we obtain \eqref{tem-54}. By a similar argument we can prove \eqref{tem-55}.
\ep

Next we show that $\xi_n,\eta_n$ are also small at points near the blow up point $x_n$. Before that, we need some preliminaries, which will be used also in the next section. For fixed $r\in(0,d_0)$, we define the quadratic forms
\be\lab{function-P}
	P(u,v):=-2r\int_{\pa B_r(x_n)} \abr{\nabla u,\nu}\abr{\nabla v,\nu}\rd\sigma+r\int_{\pa B_r(x_n)}\abr{\nabla u,\nabla v}\rd\sigma,
\ee
and for $i=1,2$,
\be\lab{function-Q}
	Q_i(u,v):=-\int_{\pa B_r(x_n)} \abr{\nabla u,\nu}\pa_i v+\abr{\nabla v,\nu}\pa_i u\rd \sigma+\int_{\pa B_r(x_n)}\abr{\nabla u,\nabla v}\nu_i\rd\sigma,
\ee
where $\nu=(\nu_1,\nu_2)$ denotes the outer normal vector.
We recall the following Pohozaev identities.
\begin{lemma} (\cite[Lemma 2.1]{Chen-Li-Zou})
	It holds
	\be\lab{formula-P-1}
		P(u_n,v_n)=r\int_{\pa B_r(x_n)} \f{u_n^{q_n+1}}{q_n+1}+\f{v_n^{p_n+1}}{p_n+1}\rd\sigma-2\int_{B_r(x_n)}\f{u_n^{q_n+1}}{q_n+1}+\f{v_n^{p_n+1}}{p_n+1}\rd x,
	\ee
	and
	\be\lab{formula-Q-1}
		Q_i(u_n,v_n)=\int_{\pa B_r(x_n)} \sbr{ \f{u_n^{q_n+1}}{q_n+1}+\f{v_n^{p_n+1}}{p_n+1} }\nu_i\rd\sigma.
	\ee
\end{lemma}

\bl
	It holds
	\be\lab{formula-P-2}
		P(u_n,\eta_n)+P(v_n,\xi_n)=r\int_{\pa B_r(x_n)} u_n^{q_n}\xi_n+v_n^{p_n}\eta_n \rd\sigma-2\int_{B_r(x_n)}u_n^{q_n}\xi_n+v_n^{p_n}\eta_n \rd x,
	\ee
	and
	\be\lab{formula-Q-2}
		Q_i(u_n,\eta_n)+Q_i(v_n,\xi_n)=\int_{\pa B_r(x_n)} \sbr{ u_n^{q_n}\xi_n+v_n^{p_n}\eta_n }\nu_i\rd\sigma.
	\ee
\el
\bp
From the equation satisfied by $(u_n,v_n)$ and $(\xi_n,\eta_n)$, we see that
$$\begin{aligned}
	&\quad\text{div}[(x-x_n)\nabla u_n\cdot \nabla\eta_n-\nabla u_n (x-x_n)\cdot \nabla\eta_n-\nabla \eta_n (x-x_n)\cdot \nabla u_n]\\
&\quad+\text{div}[(x-x_n)\nabla v_n\cdot \nabla\xi_n-\nabla v_n (x-x_n)\cdot \nabla\xi_n-\nabla \xi_n (x-x_n)\cdot \nabla v_n]\\
&=(x-x_n)\cdot\mbr{ \nabla\eta_n(-\Delta u_n)+\nabla u_n(-\Delta \eta_n)+\nabla\xi_n(-\Delta v_n)+\nabla v_n(-\Delta \xi_n)}\\
	&=(x-x_n)\cdot\mbr{ \nabla(\eta_nv_n^{p_n})+\nabla(\xi_nu_n^{q_n}) }\\
&=\text{div}((x-x_n)(\eta_nv_n^{p_n}+\xi_nu_n^{q_n}))-2(\eta_nv_n^{p_n}+\xi_nu_n^{q_n}).
\end{aligned}$$
Then by integrating in $B_r(x_n)$ and using the divergence theorem, we obtain \eqref{formula-P-2}.

Similarly, for $i=1,2$, we have
	\begin{align*}
&\partial_{i}(\nabla u_n\cdot\nabla \eta_n+\nabla v_n\cdot\nabla \xi_n)-\text{div}(\nabla u_n\partial_i\eta_n+\nabla \eta_n\partial_i u_n+\nabla v_n\partial_i\xi_n+\nabla \xi_n\partial_i v_n)\\
=&\pa_i\eta_n(-\Delta u_n)+\pa_i u_n(-\Delta \eta_n)+\pa_i\xi_n(-\Delta v_n)+\pa_i v_n(-\Delta \xi_n)\\
=&\pa_i(\eta_nv_n^{p_n})+\pa_i(\xi_nu_n^{q_n}),\end{align*}
from which we immediately obtain \eqref{formula-Q-2}.
\ep

\bl\label{lemma-5-6}
	\begin{align*}&P(G(x_n,\cdot),G(x_n,\cdot))=-\f{1}{2\pi},\\
	&P(G(x_n,\cdot),\pa_hG(x_n,\cdot))=-\f{1}{2}\pa_h R(x_n),\quad h=1,2,\\
	&Q_i(G(x_n,\cdot),G(x_n,\cdot))=-\pa_i R(x_n),\\
	&Q_i(G(x_n,\cdot),\pa_hG(x_n,\cdot))=-\f{1}{2}\pa^2_{ih} R(x_n),\quad h=1,2.
\end{align*}
\el
\bp
See \cite[Proposition 2.2]{LE-1}. Note that there are some misprints about the coefficients in \cite[Proposition 2.2]{LE-1} and we correct them here.
\ep

Using Theorem \ref{thm1}, we have a sharper estimate of $u_n,v_n$ than Lemma \ref{expansion-1}.
\bl\lab{expansion-2}
	For fixed $r\in(0,d_0)$, it holds
	\be\lab{tem-48}\bcs
		u_n(x)=C_nG(x_n,x)+O(\f{\mu_n}{p_n^{2-\delta}}),\\
		v_n(x)=\wt C_n G(x_n,x)+O(\f{\mu_n}{p_n^{2-\delta}}),
	\ecs\quad\text{in}~\CR^1(\Omega\setminus B_{2r}(x_n)), \ee
	for any given small $\delta>0$, where
	\be\lab{tem-52}\begin{aligned}
		C_n&=\int_{B_r(x_n)}v_n^{p_n}\rd x\\
&=\f{8\pi\sqrt e}{p_n}\mbr{ 1 -\f{\ln p_n}{p_n-1} +\f{1}{p_n}\sbr{ 4\pi R(x_n)+3\ln2-1-\f{3}{4}\q_n }+O(\f{1}{p_n^{2-\delta}}) },
	\end{aligned}\ee
	\be\lab{tem-53}\begin{aligned}
		\widetilde C_n&=\int_{B_r(x_n)}u_n^{q_n}\rd x\\
&=\f{8\pi\sqrt e}{p_n}\mbr{ 1 -\f{\ln p_n}{p_n-1} +\f{1}{p_n}\sbr{ 4\pi R(x_n)+3\ln2-1-\f{1}{4}\q_n }+O(\f{1}{p_n^{2-\delta}}) }.
	\end{aligned}\ee
\el
\bp
Recalling the proof of Lemma \ref{expansion-1}, to obtain \eqref{tem-48} we only need to prove
\be\lab{tem-50}
	\int_{B_{r}(x_n)} (y-x_n)_h v_n^{p_n}(y)\rd y=O(\f{\mu_n}{p_n^{2-\delta}}),
\ee
for $h=1,2$. Similarly to the proof of \eqref{est-1}, we actually have
$$\begin{aligned}
	\int_{B_{r}(x_n)} (y-x_n)_hv_n^{p_n}(y)\rd y
	&=\f{\mu_nv_n(x_n)}{p_n}\int_{B_{\f{r}{\mu_n}}(0)}y_h\sbr{1+\f{z_n(y)}{p_n}}^{p_n}\rd y\\
	&=\f{\mu_nv_n(x_n)}{p_n}\sbr{\int_{B_{p_n}(0)}+\int_{B_{\f{r}{\mu_n}}(0)\setminus B_{p_n}(0)}}y_h\sbr{1+\f{z_n(y)}{p_n}}^{p_n}\rd y.
\end{aligned}$$
By Remark \ref{decay-10}, we have
	$$\begin{aligned}
	\abs{\int_{B_{\f{r}{\mu_n}}(0)\setminus B_{p_n}(0)}y_h\sbr{1+\f{z_n(y)}{p_n}}^{p_n}\rd y}
	&\le C\int_{p_n}^{\f{r}{\mu_n}}\f{r^2}{1+r^{4-\delta}}\rd r=O(\f{1}{p_n^{1-\delta}}).
	\end{aligned}$$
 Meanwhile, by \eqref{tem-49} and the symmetry of $y_he^U\sbr{ 1+\f{1}{p_n}(t_{n}^*-\f{1}{2}U^2) }$ we get
$$\begin{aligned}
	\int_{B_{p_n}(0)}y_h\sbr{1+\f{z_n(y)}{p_n}}^{p_n}\rd y
	&=\int_{B_{p_n}(0)}y_he^U\sbr{ 1+\f{1}{p_n}(t_{n}^*-\f{1}{2}U^2) }\rd y+O(\f{1}{p_n^{2-\delta}})\\
	&=\int_{\R^2}y_he^U\sbr{ 1+\f{1}{p_n}(t_{n}^*-\f{1}{2}U^2) }\rd y+O(\f{1}{p_n^{1-\delta}})\\
	&=0+O(\f{1}{p_n^{1-\delta}}).
\end{aligned}$$
This proves \eqref{tem-50} and so \eqref{tem-48}. Finally, since
\begin{align*}C_n&=\int_{B_r(x_n)}v_n^{p_n}\rd x=\frac{v_n(x_n)}{p_n}\int_{B_{\frac{r}{\mu_n}}(0)}\left(1+\frac{z_n}{p_n}\right)^{p_n}dx,\\
\wt C_{n}&=\int_{B_r(x_n)}u_n^{q_n}\rd x=\f{v_n(x_{n})^{\q_n+1}}{p_n}\int_{B_{\f{r}{\mu_{n}}}(0)}\sbr{1+\f{w_{n}}{p_n}}^{q_n}\rd x,
\end{align*}
and note from $k=i=1$ that $\Phi_{k,i}(\boldsymbol x)=R(x)$, the assertions \eqref{tem-52}-\eqref{tem-53} follow readily from \eqref{est-1}-\eqref{est-2}, \eqref{c4-80}-\eqref{c4-81} and \eqref{est-v}.
\ep

By this sharp estimate, we prove that
\bl\label{lemma-5-8}
	Recall that $R(x)=H(x,x)$ is the Robin function. Then
		$$\nabla R(x_n)=O(\f{\mu_n}{p_n^{1-\delta}}), \quad\text{for small}~\delta>0.$$
\el
\bp
Using \eqref{formula-Q-1} and \eqref{c2-31}, we see that
	$$Q_i(u_n,v_n)=\int_{\pa B_r(x_n)} \sbr{ \f{u_n^{q_n+1}}{q_n+1}+\f{v_n^{p_n+1}}{p_n+1} }\nu_i\rd\sigma=O(\f{C^{p_n+1}}{p_n^{p_n+2}}).$$
Using the expansions in Lemma \ref{expansion-2}, \eqref{function-Q} and Lemma \ref{lemma-5-6}, we obtain
$$\begin{aligned}
	Q_i(u_n,v_n)
	&=C_n\wt C_n Q_i(G(x_n,\cdot),G(x_n,\cdot))+O(\f{\mu_n}{p_n^{3-\delta}})\\
	&=-C_n\wt C_n \pa_iR(x_n)+O(\f{\mu_n}{p_n^{3-\delta}}),\quad\text{for any fixed small}~\delta>0.
\end{aligned}$$
Since $C_n,\wt C_n=\frac{8\pi\sqrt{e}}{p_n}(1+o(1))$, we deduce that $\nabla R(x_n)=O(\f{\mu_n}{p_n^{1-\delta}})$.
\ep

\bl\lab{lemma-in}
	It holds $a_0=0$, and further if $x_\iy$ is a nondegenerate critical point of $R(x)$, then $a_1=a_2=0$, where $a_j$ are constants in \eqref{tem-51}. Thus $\wt\xi_\iy=\wt\eta_\iy=0$.
\el
\bp
By \eqref{formula-P-2}, \eqref{tem-46} and \eqref{c2-31}, we have
\be \label{c5-50} P(u_n,\eta_n)+P(v_n,\xi_n)=O(\f{C^{p_n}}{p_n^{p_n}})+O(\f{\mu_n}{p_n^2})=O(\f{\mu_n}{p_n^2}).\ee
While using \eqref{function-P}, Lemma \ref{lemma-out} and Lemmas \ref{lemma-5-6}-\ref{lemma-5-8}, we obtain
$$\begin{aligned}
	P(u_n,\eta_n)
	&=C_n\wt A_{n,0}P(G(x_n,\cdot),G(x_n,\cdot))+\sum_{j=1}^2C_n\wt A_{n,j}P(G(x_n,\cdot),\pa_jG(x_n,\cdot))+o(\f{\mu_n}{p_n})\\
	&=-\f{1}{2\pi}C_n\wt A_{n,0}-\f{1}{2}\sum_{j=1}^2C_n\wt A_{n,j}\pa_jR(x_n)+o(\f{\mu_n}{p_n})\\
	&=-\f{32\pi\sqrt e}{p_n^2}a_0+o(\f{1}{p_n^2}),
\end{aligned}$$
and
$$\begin{aligned}
	P(v_n,\xi_n)
	&=\wt C_nA_{n,0}P(G(x_n,\cdot),G(x_n,\cdot))+\sum_{j=1}^2\wt C_n A_{n,j}P(G(x_n,\cdot),\pa_jG(x_n,\cdot))+o(\f{\mu_n}{p_n})\\
	&=-\f{1}{2\pi}\wt C_nA_{n,0}-\f{1}{2}\sum_{j=1}^2\wt C_n A_{n,j}\pa_jR(x_n)+o(\f{\mu_n}{p_n})\\
	&=-\f{32\pi\sqrt e}{p_n^2}a_0+o(\f{1}{p_n^2}).
\end{aligned}$$
Inserting these estimates into \eqref{c5-50} leads to
	$$-\f{64\pi\sqrt e}{p_n^2}a_0+o(\f{1}{p_n^2})=O(\f{\mu_n}{p_n^2}),$$
which implies $a_0=0$.

Now we prove $a_1=a_2=0$. By \eqref{formula-Q-2} and \eqref{c2-31} we have
	\be\label{c5-51}Q_i(u_n,\eta_n)+Q_i(v_n,\xi_n)=O(\f{C^{p_n}}{p_n^{p_n}}),\quad\text{for }i=1,2.\ee
Since $a_0=0$, it follows from Lemma \ref{lemma-out} that $A_{n,0},\wt A_{n,0}=o(\f{1}{p_n})$. Then similarly,
$$\begin{aligned}
	Q_i(u_n,\eta_n)
	&=C_n\wt A_{n,0}Q_i(G(x_n,\cdot),G(x_n,\cdot))+\sum_{j=1}^2C_n\wt A_{n,j}Q_i(G(x_n,\cdot),\pa_jG(x_n,\cdot))+o(\f{\mu_n}{p_n})\\
	&=-C_n\wt A_{n,0}\pa_iR(x_n)-\f{1}{2}\sum_{j=1}^2C_n\wt A_{n,j}\pa_{ij}^2R(x_n)+o(\f{\mu_n}{p_n})\\
	&=-\f{1}{2}\sum_{j=1}^2C_n\wt A_{n,j}\pa_{ij}^2R(x_n)+o(\f{\mu_n}{p_n}),
\end{aligned}$$
and
$$\begin{aligned}
	Q_i(v_n,\xi_n)
	&=\wt C_n A_{n,0}Q_i(G(x_n,\cdot),G(x_n,\cdot))+\sum_{j=1}^2\wt C_n A_{n,j}Q_i(G(x_n,\cdot),\pa_jG(x_n,\cdot))+o(\f{\mu_n}{p_n})\\
	&=-\f{1}{2}\sum_{j=1}^2\wt C_n A_{n,j}\pa_{ij}^2R(x_n)+o(\f{\mu_n}{p_n}).
\end{aligned}$$
Inserting these into \eqref{c5-51} we get
	$$\nabla^2R(x_n)\cdot\sbr{ C_n\wt A_{n,1}+\wt C_n A_{n,1}, C_n\wt A_{n,2}+\wt C_n A_{n,2} }^T=o(\f{\mu_n}{p_n}).$$
Since $x_\infty$ is a nondegenerate critical point of $R(x)$, i.e. $\nabla^2R(x_\iy)$ is invertible, and
since $\nabla^2R(x_n)\to \nabla^2R(x_\iy)$, we conclude that
	$$C_n\wt A_{n,j}+\wt C_n A_{n,j}=o(\f{\mu_n}{p_n}),\quad j=1,2,$$
which together with \eqref{c5-5} and \eqref{tem-52}-\eqref{tem-53} implies $a_1=a_2=0$.
\ep

\bc\lab{in}
	For any $d>0$, it holds
	$$\nm{\xi_n}_{L^\iy(B_{\mu_nd}(x_n))}=o(1),\quad \nm{\eta_n}_{L^\iy(B_{\mu_nd}(x_n))}=o(1).$$
\ec

\begin{proof}
Since Lemma \ref{lemma-in} implies
	$$(\wt\xi_n,\wt\eta_n)\to(0,0)\quad \text{in}~ \CR_{loc}^2(\R^2),$$
so $\nm{\xi_n}_{L^\iy(B_{\mu_nd}(x_n))}=\|\wt \xi_n\|_{L^\iy(B_{d}(0))}=o(1)$ and similarly $\nm{\eta_n}_{L^\iy(B_{\mu_nd}(x_n))}=o(1)$.
\end{proof}

Now we are in position to prove Theorem \ref{thm2}.
\bp[Proof of Theorem \ref{thm2}]
Recall that by contradiction, we have supposed that $(\xi_n,\eta_n)$ is a solution sequence of the linearized system \eqref{linear-0} satisfying
	$$\nm{\xi_n}_{L^\iy(\Omega)}=\max\lbr{\nm{\xi_n}_{L^\iy(\Omega)},~\nm{\eta_n}_{L^\iy(\Omega)}}=1,\quad\text{for any}~n\ge1.$$
Without loss of generality, we assume
	$$\xi_n(x_n^*)=1=\max_\Omega\lbr{|\xi_n|,|\eta_n|}.$$
For any $d>0$, since $A_{n,0},\wt A_{n,0}=o(\f{1}{p_n})$, Lemma \ref{out-2} and Corollary \ref{in} imply $x_n^*\in B_{2\mu_np_n}(x_n)\setminus B_{\mu_nd}(x_n)$ for $n$ large, i.e.
\be\lab{constant-rn}
	r_n:=|x_n^*-x_n|\in [\mu_nd, 2\mu_np_n].
\ee
Recalling $\wt \xi_n,\wt \eta_n$ defined in \eqref{function-wtxi}, we
define
	$$\xi_n^{**}(x):=\xi_n(x_n+r_nx)=\wt \xi_n(\frac{r_n}{\mu_n}x),\quad \eta_n^{**}(x):=\eta_n(x_n+r_nx)=\wt \eta_n(\frac{r_n}{\mu_n}x).$$
Then we see from \eqref{c5-2} that
\be\bcs
	-\Delta \xi_n^{**}=\f{r_n^2}{\mu_n^2}\sbr{1+\f{z_n(\f{r_n}{\mu_n}x)}{p_n}}^{p_n-1}\eta_n^{**},\\
	-\Delta \eta_n^{**}=\f{r_n^2}{\mu_n^2}\f{q_n}{p_n}v_n^{\q_n}(x_n)\sbr{1+\f{w_n(\f{r_n}{\mu_n}x)}{p_n}}^{q_n-1}\xi_n^{**}.
\ecs\ee
Since for $x\neq0$ and $n$ large, we see from \eqref{constant-rn} and Remark \ref{decay-10} that
	$$\abs{ \f{r_n^2}{\mu_n^2}\sbr{1+\f{z_n(\f{r_n}{\mu_n}x)}{p_n}}^{p_n-1}\eta_n^{**}(x) }\le \f{C}{d^{2-\delta}|x|^{4-\delta}},$$
	$$\abs{ \f{r_n^2}{\mu_n^2}\f{q_n}{p_n}v_n^{\q_n}(x_n)\sbr{1+\f{w_n(\f{r_n}{\mu_n}x)}{p_n}}^{q_n-1}\xi_n^{**}(x) }\le \f{C}{d^{2-\delta}|x|^{4-\delta}},$$
by letting $n\to\iy$ first and then $d\to\iy$, we obtain that up to a subsequence, $(\xi_n^{**},\eta_n^{**})\to(\xi_\iy^{**},\eta_\iy^{**})$ in $\CR_{loc}^2(\R^2\setminus\{0\})$, where
	$$-\Delta \xi_\iy^{**}=0,\quad -\Delta \eta_\iy^{**}=0,\quad\text{in}~\R^2\setminus\{0\}.$$
Since $|\f{x_n^*-x_n}{r_n}|=1$ and $\xi_n^{**}(\f{x_n^*-x_n}{r_n})=1=\max |\xi_n^{**}|$, we obtain $\xi_\iy\equiv1$. Thus $\inf_{|x|=1}\xi_n^{**}(x)\ge C>0$ for $n$ large. Define the average functions
\be\lab{tem-60}
	\xi_n^*(r):=\f{1}{2\pi}\int_0^{2\pi}\wt\xi_n(r,\rho)\rd\rho,\quad \eta_n^*(r):=\f{1}{2\pi}\int_0^{2\pi}\wt\eta_n(r,\rho)\rd\rho.
\ee
We claim
\be\lab{tem-6000}
	\xi_n^*(\f{r_n}{\mu_n})=o(1),\quad \eta_n^*(\f{r_n}{\mu_n})=o(1),
\ee
where $r_n$ is defined in \eqref{constant-rn}. Using this claim, we see that
	$$o(1)=\xi_n^*(\f{r_n}{\mu_n})=\f{1}{2\pi}\int_0^{2\pi}\wt\xi_n(\f{r_n}{\mu_n},\rho)\rd\rho
=\f{1}{2\pi}\int_0^{2\pi}\xi_n^{**}(1,\rho)\rd\rho\ge C>0,$$
a contradiction, which proves Theorem \ref{thm2}.

\vskip0.1in
So it suffices to prove the claim \eqref{tem-60}. From \eqref{c5-2}, we have
\be\label{c5-60} \bcs
	-\Delta\wt\xi_n-e^U\wt\eta_n=f_n:=\sbr{\sbr{1+\f{z_n}{p_n}}^{p_n}-e^U}\wt\eta_n,\\
	-\Delta\wt\eta_n-e^U\wt\xi_n=\wt f_n:=\sbr{ \f{q_n}{p_n}v_n^{\q_n}(x_n)\sbr{1+\f{w_n}{p_n}}^{q_n}-e^U } \wt\xi_n,
\ecs \text{in}~\Omega_n=\frac{\Omega-x_n}{\mu_n}, \ee
First we note that
\be\lab{tem-112}
	f_n(x),\wt f_n(x)=O\sbr{\f{1}{p_n(1+|x|^{4-\delta})}},\quad\text{for}~x\in \Omega_n.
\ee
Indeed, for $x\in B_{\f{d_0}{2\mu_n}}(0)$, by \eqref{c4-2}-\eqref{c4-8}, Lemmas \ref{decay-1}-\ref{decay-12} and Lemma \ref{bound-1}, we get	
\begin{align*}
f_n(x)&=\frac{\wt \eta_n}{p_n}[t_{n}(e^U+g_n)+h_n]\\
&=\frac{\wt \eta_n}{p_n}\left[t_{n}\left(O(\frac{1}{1+|x|^{4-\frac{\delta}{2}}})+O(\frac{\ln(2+|x|)}{1+|x|^{4-\frac{\delta}{2}}})\right)+O(\frac{1}{1+|x|^{4-\frac{\delta}{2}}})\right]\\
&=O\sbr{\f{1}{p_n(1+|x|^{4-\delta})}},
\end{align*}
and similarly
\begin{align*}
	\wt f_n(x)&= \f{\theta_n \wt\xi_n}{p_n}v_n^{\q_n}(x_n)\sbr{1+\f{w_n}{p_n}}^{q_n}+
\frac{\wt \xi_n}{p_n}[s_{n}(e^U+\wt g_n)+\wt h_n]  \\
&=O\sbr{\f{1}{p_n(1+|x|^{4-\delta})}}.
\end{align*}
While for $x\in \Omega_n\setminus B_{\f{d_0}{2\mu_n}}(0)$, from Lemma \ref{lemma-out} we know $\wt\xi_n,\wt\eta_n=O(\f{1}{p_n})$ and hence Remark \ref{decay-10} implies
	$$\max\lbr{|f_n(x)|,~|\wt f_n(x)|}\le \f{C\max\lbr{|\wt\xi_n(x)|,|\wt\eta_n(x)|}}{1+|x|^{4-\delta}}=O\sbr{\f{1}{p_n(1+|x|^{4-\delta})}}.$$
Thus \eqref{tem-112} holds.

By \eqref{tem-60} and \eqref{c5-60}, we see that $\xi_n^*,\eta_n^*$ solves the ODE system
\be\lab{tem-61} \bcs
	-(\xi_n^*)''-\f{1}{r}(\xi_n^*)'-e^U\eta_n^*=f_n^*:=\f{1}{2\pi}\int_0^{2\pi}f_n(r,\rho)\rd\rho,\\
	-(\eta_n^*)''-\f{1}{r}(\eta_n^*)'-e^U\xi_n^*=\wt f_n^*:=\f{1}{2\pi}\int_0^{2\pi}\wt f_n(r,\rho)\rd\rho.
\ecs \ee
Since $\vp_1(r)=\f{8-r^2}{8+r^2}$ and $\vp_2(r)=\f{(8-r^2)\ln r+16}{8+r^2}$ are two linearly independent solutions of $-\vp''-\f{1}{r}\vp'-e^U\vp=0$ (see e.g. \cite[Page 185]{LE-1}), by ODE theory we know that
	$$\vp_n(r):=\vp_1(r)\int_0^rs\vp_2(s)(f_n^*(s)+\wt f_n^*(s))\rd s -\vp_2(r)\int_0^rs\vp_1(s)(f_n^*(s)+\wt f_n^*(s))\rd s$$
is a special solution of $-\vp''-\f{1}{r}\vp'-e^U\vp=f_n^*+\wt f_n^*$. Then it follows from \eqref{tem-61} that $\xi_n^*+\eta_n^*-\vp_n$ is a solution of $-\vp''-\f{1}{r}\vp'-e^U\vp=0$, so
	$$\xi_n^*+\eta_n^*=\vp_n+c_{n,1}\vp_1+c_{n,2}\vp_2,$$
where $c_{n,1},c_{n,2}$ are constants. Since $\vp_n(0)=0$, the continuity of $\xi_n^*+\eta_n^*$ at $r=0$ implies $c_{n,2}=0$ for any $n\ge1$. Then computing at $r=0$ and using Corollary \ref{in} we get
	$$c_{n,1}=\xi_n^*(0)+\eta_n^*(0)=\wt \xi_n(0)+\wt \eta_n(0)=\xi_n(x_n)+\eta_n(x_n)=o(1).$$
Since for $r>1$, it follows from \eqref{tem-112} that
	$$\abs{\vp_n(r)}\le \f{C}{p_n}\int_0^r\f{s\ln s}{1+s^{4-\delta}}\rd s+\f{C\ln r}{p_n}\int_0^r\f{s}{1+s^{4-\delta}}\rd s\le C\f{\ln r}{p_n},$$
we obtain (Recall from \eqref{constant-rn} that $\f{r_n}{\mu_n}\in [d, 2p_n]$ and $d$ can be arbitrary large)
\be\lab{tem-62} \xi_n^*(\f{r_n}{\mu_n})+\eta_n^*(\f{r_n}{\mu_n})=\vp_n(\f{r_n}{\mu_n})+c_{n,1}\vp_1(\f{r_n}{\mu_n})=O(\f{\ln p_n}{p_n})+o(1)=o(1).
\ee

Now we assume by contradiction that the claim \eqref{tem-60} does not hold, then \eqref{tem-62} implies the existence of $C_0>0$ such that
	$$|\xi_n^*(\f{r_n}{\mu_n})|\ge C_0>0 \quad\text{and}\quad |\eta_n^*(\f{r_n}{\mu_n})|\ge C_0>0.$$
Without loss of generality, suppose $\xi_n^*(\f{r_n}{\mu_n})\ge C_0$, then \eqref{tem-62} implies $\eta_n^*(\f{r_n}{\mu_n})\le -C_0$. Since Lemma \ref{out-2} and Corollary \ref{in} imply that for large $n$,
	$$\abs{\xi_n^*(r)-\eta_n^*(r)}\le C_0,\quad\text{uniformly for}~r\in(0,d)\cup(2p_n,+\iy),$$ so there is $s_n\in [d, 2p_n]$ such that
	$$\xi_n^*(s_n)-\eta_n^*(s_n)=\max_{r>0}~(\xi_n^*(r)-\eta_n^*(r))\ge 2C_0,$$
and then
	$$(\xi_n^*)'(s_n)-(\eta_n^*)'(s_n)=0,\quad (\xi_n^*)''(s_n)-(\eta_n^*)''(s_n)\le0.$$
It follows from \eqref{tem-61} and \eqref{tem-112} that
	$$\frac{C}{1+s_n^4}\leq e^{U(s_n)}(\xi_n^*(s_n)-\eta_n^*(s_n))\le f_n^*(s_n)-\wt f_n^*(s_n)\le \f{C}{p_n(1+s_n^{4-\delta})},$$
which implies $p_n\le Cs_n^\delta\le C2^\delta p_n^\delta$, a contradiction with $0<\delta<1$ and $p_n\to\infty$. This proves the claim \eqref{tem-60} and hence completes the proof of  Theorem \ref{thm2}.
\ep

\section{Local uniqueness of \texorpdfstring{$1$}{}-bubble solutions}
In this section, we prove Theorem \ref{thm3} by contradiction. Suppose $(u_n^{(1)},v_n^{(1)})$ and $(u_n^{(2)},v_n^{(2)})$ be two sequences of different $1$-bubble solutions concentrating at the same point $x_\iy$, where $x_\iy$ is a nondegenerate critical point of the Robin function $R(x)$. Let $d_0$ be defined by \eqref{constant-d0}. For $l=1,2$, let $x_n^{(l)}$ be the local maximum points
	$$v_n^{(l)}(x_n^{(l)})=\max_{B_{d_0}(x_\iy)} v_n^{(l)},$$
the scaling parameters be
	$$\mu_n^{(l)}=p_n^{-\f{1}{2}} \sbr{v_n^{(l)}(x_n^{(l)})}^{-\f{p_n-1}{2}} ,$$
and the corresponding scaling functions $(w_n^{(l)},z_n^{(l)})$ be defined as \eqref{function-wz}.
Using Theorem \ref{thm1}, we get
\bl\lab{tem-70}
	We have
		$$\abs{x_n^{(1)}-x_n^{(2)}}=O(\f{\mu_n}{p_n^{1-\delta}}),\qquad
		\f{\mu_n^{(1)}}{\mu_n^{(2)}}-1=O(\f{1}{p_n^{1-\delta}}),$$
		$$\f{v_n^{(1)}(x_n^{(1)})}{v_n^{(2)}(x_n^{(2)})}-1=O(\f{1}{p_n^{2-\delta}}),$$
	for any given small $\delta>0$, where $\mu_n=\max\lbr{\mu_n^{(1)},\mu_n^{(2)}}$.
\el
\bp
Applying Lemma \ref{lemma-5-8} to $(u_n^{(l)}, v_n^{(l)})$, we have $\nabla R(x_n^{(l)})=O(\f{\mu_n^{(l)}}{p_n^{1-\delta}})$ for $l=1,2$. Since $\nabla^2 R(x_\iy)$ is invertible and $x_n^{(l)}\to x_\infty$, we deduce from
	$$\nabla R(x_n^{(l)})=\nabla^2 R(x_\iy)\cdot(x_n^{(l)}-x_\iy)+O(|x_n^{(l)}-x_\iy|^2),\quad l=1,2,$$
that $|x_n^{(l)}-x_\iy|=O(\f{\mu_n}{p_n^{1-\delta}})$, so $|x_n^{(1)}-x_n^{(2)}|=O(\f{\mu_n}{p_n^{1-\delta}})$.

By the Taylor's expansion and $\nabla R(x_\iy)=0$, we have that
	$$R(x_n^{(l)})=R(x_\iy)+O(|x_n^{(l)}-x_\iy|^2),$$
and hence
	$$R(x_n^{(1)})-R(x_n^{(2)})=O\sbr{|x_n^{(1)}-x_\iy|^2+|x_n^{(2)}-x_\iy|^2}=O(\f{\mu_n^2}{p_n^{2-2\delta}}).$$
Then \eqref{rate-1}-\eqref{rate-2} imply that (note from $k=i=1$ that $\Phi_{k,i}(\boldsymbol x)=R(x)$) $$\f{\mu_n^{(1)}}{\mu_n^{(2)}}=e^{-2\pi(R(x_n^{(1)})-R(x_n^{(2)}))}+O(\f{1}{p_n^{1-\delta}})=1+O(\f{1}{p_n^{1-\delta}}),$$
and $$\f{v_n^{(1)}(x_n^{(1)})}{v_n^{(2)}(x_n^{(2)})}=1+\f{4\pi}{p_n}(R(x_n^{(1)})-R(x_n^{(2)}))+O(\f{1}{p_n^{2-\delta}})=1+O(\f{1}{p_n^{2-\delta}}).$$
This completes the proof.\ep

\bl\lab{expansion-3}
	For fixed $r\in(0,d_0)$ and $l=1,2$, it holds
	\be\label{c6-31} \bcs
		u_n^{(l)}(x)=C_n^{(l)} G(x_n^{(1)},x)+O(\f{\mu_n}{p_n^{2-\delta}}),\\
		v_n^{(l)}(x)=\wt C_n^{(l)} G(x_n^{(1)},x)+O(\f{\mu_n}{p_n^{2-\delta}}),
	\ecs\quad\text{in}~\CR^1(\Omega\setminus B_{2r}(x_n^{(1)})), \ee
	for any given small $\delta>0$, where
	\be\label{c6-32}
		C_n^{(l)}=\f{8\pi\sqrt e}{p_n}\mbr{ 1 -\f{\ln p_n}{p_n-1} +\f{1}{p_n}\sbr{ 4\pi R(x_n^{(l)})+3\ln2-1-\f{3}{4}\q_n }+O(\f{1}{p_n^{2-\delta}}) },
	\ee
	and
	\be
		\wt C_n^{(l)}=\f{8\pi\sqrt e}{p_n}\mbr{ 1 -\f{\ln p_n}{p_n-1} +\f{1}{p_n}\sbr{ 4\pi R(x_n^{(l)})+3\ln2-1-\f{1}{4}\q_n }+O(\f{1}{p_n^{2-\delta}}) }.
	\ee
\el
\bp
This lemma follows directly from Lemma \ref{expansion-2} and
$$\begin{aligned}
	u_n^{(2)}(x)
	&=C_n^{(2)} G(x_n^{(2)},x)+O(\f{\mu_n^{(2)}}{p_n^{2-\delta}})\\
	&=C_n^{(2)} \sbr{G(x_n^{(1)},x)+O(|x_n^{(1)}-x_n^{(2)}|)}+O(\f{\mu_n^{(2)}}{p_n^{2-\delta}})\\
	&=C_n^{(2)} G(x_n^{(1)},x)+O(\f{\mu_n}{p_n^{2-\delta}}),
\end{aligned}$$
in $\CR^1(\Omega\setminus B_{2r}(x_n^{(1)}))$.
\ep

\vskip0.1in

Since $(u_n^{(1)},v_n^{(1)})\neq (u_n^{(2)},v_n^{(2)})$ for any $n\ge1$, we suppose without loss of generality that
\be\lab{constant-Ln}	L_n:=\nm{u_n^{(1)}-u_n^{(2)}}_{L^\iy(\Omega)}=\max\lbr{\nm{u_n^{(1)}-u_n^{(2)}}_{L^\iy(\Omega)},\nm{v_n^{(1)}-v_n^{(2)}}_{L^\iy(\Omega)}}>0.
\ee
Set
\be\lab{function-zeta}
	\zeta_n:=\f{u_n^{(1)}-u_n^{(2)}}{L_n},\quad \kappa_n:=\f{v_n^{(1)}-v_n^{(2)}}{L_n},
\ee
then
\be \bcs -\Delta\zeta_n=\f{(v_n^{(1)})^{p_n}-(v_n^{(2)})^{p_n}}{L_n}=D_n\kappa_n,\quad\text{in}~\Omega,\\
	-\Delta\kappa_n=\f{(u_n^{(1)})^{q_n}-(u_n^{(2)})^{q_n}}{L_n}=\wt D_n\zeta_n,\quad\text{in}~\Omega,\\
|\zeta_n|,|\kappa_n|\le 1,\quad\text{in}~\Omega,\\
	\zeta_n=\kappa_n=0,\quad\text{on}~\pa\Omega,
\ecs \ee
where
\be\label{c6-7}
	D_n(x):=p_n\int_0^1\sbr{tv_n^{(1)}(x)+(1-t)v_n^{(2)}(x)}^{p_n-1}\rd t,
\ee
and
\be\label{c6-8}
	\wt D_n(x):=q_n\int_0^1\sbr{tu_n^{(1)}(x)+(1-t)u_n^{(2)}(x)}^{q_n-1}\rd t.
\ee
Now we define
\be\lab{function-wtzeta}
	\wt\zeta_n(x):=\zeta_n(x_n^{(1)}+\mu_n^{(1)}x),\quad \wt\kappa_n(x):=\kappa_n(x_n^{(1)}+\mu_n^{(1)}x).
\ee
Then
\be\lab{tem-71} \bcs
	-\Delta\wt\zeta_n=(\mu_n^{(1)})^2D_n(x_n^{(1)}+\mu_n^{(1)}x) \wt\kappa_n,\quad\text{in}~\Omega_n^{(1)}:=\f{\Omega-x_n^{(1)}}{\mu_n^{(1)}},\\
	-\Delta\wt\kappa_n=(\mu_n^{(1)})^2\wt D_n(x_n^{(1)}+\mu_n^{(1)}x) \wt\zeta_n,\quad\text{in}~\Omega_n^{(1)},\\
|\wt \zeta_n|,|\wt \kappa_n|\le 1,\quad\text{in}~\Omega_n^{(1)}.
\ecs \ee

\bl\lab{tem-76}
	Up to a subsequence, it holds $(\wt\zeta_n,\wt\kappa_n)\to(\wt\zeta_\iy,\wt\kappa_\iy)$ in $\CR_{\loc}^2(\R^2)$, where
		$$(\wt\zeta_\iy,\wt\kappa_\iy)=\sum_{j=0}^2b_j(\phi_j,\phi_j),$$
	for some constants $b_j\in\R$, $j=0,1,2$.
\el
\bp
To use the elliptic estimates, we have to study $D_n$ and $\wt D_n$. Set
\be
	E(t,x):=tv_n^{(1)}(x)+(1-t)v_n^{(2)}(x),
\ee
then $D_n(x)=p_n\int_0^1E(t,x)^{p_n-1}\rd t$. By direct computations, we have
\be\label{c6-2}\begin{aligned}
	&\quad E(t,x_n^{(1)}+\mu_n^{(1)}x)\\	&=tv_n^{(1)}(x_n^{(1)})\sbr{1+\f{z_n^{(1)}(x)}{p_n}}+(1-t)v_n^{(2)}(x_n^{(2)})\sbr{1+\f{z_n^{(2)}\sbr{\f{\mu_n^{(1)}}{\mu_n^{(2)}}x+\f{x_n^{(1)}-x_n^{(2)}}{\mu_n^{(2)}}}}{p_n}}\\
	&=v_n^{(1)}(x_n^{(1)})\sbr{ t+(1-t)\f{v_n^{(2)}(x_n^{(2)})}{v_n^{(1)}(x_n^{(1)})} }\\
		&\quad +\f{v_n^{(1)}(x_n^{(1)})}{p_n}\sbr{ tz_n^{(1)}(x)+(1-t)\f{v_n^{(2)}(x_n^{(2)})}{v_n^{(1)}(x_n^{(1)})}z_n^{(2)}\sbr{ \f{\mu_n^{(1)}}{\mu_n^{(2)}}x+\f{x_n^{(1)}-x_n^{(2)}}{\mu_n^{(2)}} } }\\
	&=:v_n^{(1)}(x_n^{(1)})K_{n,1}+\f{v_n^{(1)}(x_n^{(1)})}{p_n}K_{n,2}.
\end{aligned}\ee
From Lemma \ref{tem-70}, we see that $$K_{n,1}=t+(1-t)\sbr{1+O(\f{1}{p_n^{2-\delta}})}=1+O(\f{1}{p_n^{2-\delta}}),\quad\text{uniformly in}~t\in[0,1].$$

Fix any $r\in (0, d_0)$.
Recalling the corresponding $t_n^{(l)}:=p_n(z_n^{(l)}-U)$ in \eqref{function-stn} for $l=1,2$, by Lemma \ref{bound-1} and Lemma \ref{tem-70}, we get that for $x\in B_{\f{r}{\mu_n^{(1)}}}(0)$,
	$$z_n^{(1)}(x)=U(x)+\f{t_n^{(1)}(x)}{p_n}=U(x)+O\sbr{\f{1+|x|^\tau}{p_n}},$$
and
$$\begin{aligned}
	z_n^{(2)}\sbr{ \f{\mu_n^{(1)}}{\mu_n^{(2)}}x+\f{x_n^{(1)}-x_n^{(2)}}{\mu_n^{(2)}} }
	&=U\sbr{ \f{\mu_n^{(1)}}{\mu_n^{(2)}}x+\f{x_n^{(1)}-x_n^{(2)}}{\mu_n^{(2)}} }+\f{1}{p_n}t_n^{(2)}\sbr{ \f{\mu_n^{(1)}}{\mu_n^{(2)}}x+\f{x_n^{(1)}-x_n^{(2)}}{\mu_n^{(2)}} }\\
	&=U(x)+\nabla U(\xi_{n,x})\cdot\sbr{ \sbr{\f{\mu_n^{(1)}}{\mu_n^{(2)}}-1} x+\f{x_n^{(1)}-x_n^{(2)}}{\mu_n^{(2)}} }\\
		&\quad +O\sbr{\f{1}{p_n}\sbr{1+\abs{\f{\mu_n^{(1)}}{\mu_n^{(2)}}x+\f{x_n^{(1)}-x_n^{(2)}}{\mu_n^{(2)}}}^\tau} }\\
	&=U(x)+O(\f{1}{p_n^{1-\delta}})+O\sbr{\f{1+|x|^\tau}{p_n}}.
\end{aligned}$$
for any given small $\tau>0$ (here $\xi_{n,x}$ lies on the segment determined by $x$ and $\f{\mu_n^{(1)}}{\mu_n^{(2)}}x+\f{x_n^{(1)}-x_n^{(2)}}{\mu_n^{(2)}}$). Thus
$$\begin{aligned}
	K_{n,2}
	&=t\sbr{U(x)+O\sbr{\f{1+|x|^\tau}{p_n}}}\\
		&\quad +(1-t)\sbr{1+O(\f{1}{p_n^{2-\delta}})}\sbr{U(x)+O(\f{1}{p_n^{1-\delta}})+O\sbr{\f{1+|x|^\tau}{p_n}}}\\
	&=U(x)+O\sbr{\f{1+|x|^\tau}{p_n^{1-\delta}}},\quad\text{uniformly in}~t\in[0,1],~x\in B_{\f{r}{\mu_n^{(1)}}}(0).
\end{aligned}$$
Inserting these estimates into \eqref{c6-2}, we obtain
	$$E(t,x_n^{(1)}+\mu_n^{(1)}x)=v_n^{(1)}(x_n^{(1)})\sbr{ 1+\f{1}{p_n}U(x)+O(\f{1+|x|^\tau}{p_n^{2-\delta}}) },$$
uniformly in $t\in[0,1],~x\in B_{\f{r}{\mu_n^{(1)}}}(0)$. It follows that
\be\lab{tem-90} \begin{aligned}
	(\mu_n^{(1)})^2D_n(x_n^{(1)}+\mu_n^{(1)}x)
	&=(\mu_n^{(1)})^2p_n\int_0^1 E(t,x_n^{(1)}+\mu_n^{(1)}x)^{p_n-1}\rd t\\
	&=\sbr{ 1+\f{1}{p_n}U(x)+O(\f{1+|x|^\tau}{p_n^{2-\delta}}) }^{p_n-1},
\end{aligned} \ee
uniformly in $x\in B_{\f{r}{\mu_n^{(1)}}}(0)$. In particular,
\be\label{c6-4}(\mu_n^{(1)})^2D_n(x_n^{(1)}+\mu_n^{(1)}x)\to e^{U(x)}\quad\text{uniformly in }\mathcal{C}_{loc}(\R^2).\ee

Similarly, we can give an estimate for $\wt D_n$. Set
\be
	\wt E(t,x):=tu_n^{(1)}(x)+(1-t)u_n^{(2)}(x),
\ee
then $\wt D_n(x)=q_n\int_0^1\wt E(t,x)^{q_n-1}\rd t$. By direct computations, we have
\be\label{c6-3}\begin{aligned}
	&\quad \wt E(t,x_n^{(1)}+\mu_n^{(1)}x)\\ &=tv_n^{(1)}(x_n^{(1)})\sbr{1+\f{w_n^{(1)}(x)}{p_n}}+(1-t)v_n^{(2)}(x_n^{(2)})\sbr{1+\f{w_n^{(2)}\sbr{\f{\mu_n^{(1)}}{\mu_n^{(2)}}x+\f{x_n^{(1)}-x_n^{(2)}}{\mu_n^{(2)}}}}{p_n}}\\
	&=v_n^{(1)}(x_n^{(1)})\sbr{ t +(1-t)\f{v_n^{(2)}(x_n^{(2)})}{v_n^{(1)}(x_n^{(1)})} }\\
		&\quad +\f{v_n^{(1)}(x_n^{(1)})}{p_n}\sbr{ tw_n^{(1)}(x)+(1-t)\f{v_n^{(2)}(x_n^{(2)})}{v_n^{(1)}(x_n^{(1)})}w_n^{(2)}\sbr{ \f{\mu_n^{(1)}}{\mu_n^{(2)}}x+\f{x_n^{(1)}-x_n^{(2)}}{\mu_n^{(2)}} } }\\
	&=:v_n^{(1)}(x_n^{(1)})K_{n,1}+\f{v_n^{(1)}(x_n^{(1)})}{p_n}K_{n,3}.
\end{aligned}\ee
Recalling the corresponding $s_n^{(l)}:=p_n(w_n^{(l)}+\theta_n\ln v_n^{(l)}(x_n^{(l)})-U)$ in \eqref{function-stn} for $l=1,2$, by Lemma \ref{bound-1} and Lemma \ref{tem-70}, we get that for $x\in B_{\f{r}{\mu_n^{(1)}}}(0)$,
$$\begin{aligned}
	w_n^{(1)}(x)
	&=U(x)-\q_n\ln v_n^{(1)}(x_n^{(1)})+\f{s_n^{(1)}(x)}{p_n}\\
	&=U(x)-\q_n\ln v_n^{(1)}(x_n^{(1)})+O\sbr{\f{1+|x|^\tau}{p_n}},
\end{aligned}$$
and
$$\begin{aligned}
	w_n^{(2)}\sbr{ \f{\mu_n^{(1)}}{\mu_n^{(2)}}x+\f{x_n^{(1)}-x_n^{(2)}}{\mu_n^{(2)}} }
	&=U\sbr{ \f{\mu_n^{(1)}}{\mu_n^{(2)}}x+\f{x_n^{(1)}-x_n^{(2)}}{\mu_n^{(2)}} }-\q_n\ln v_n^{(2)}(x_n^{(2)}) \\
		&\quad +\f{1}{p_n}s_n^{(2)}\sbr{ \f{\mu_n^{(1)}}{\mu_n^{(2)}}x+\f{x_n^{(1)}-x_n^{(2)}}{\mu_n^{(2)}} }\\
	&=U(x)-\q_n\ln v_n^{(1)}(x_n^{(1)})+O(\f{1}{p_n^{1-\delta}})+O\sbr{\f{1+|x|^\tau}{p_n}}.
\end{aligned}$$
Thus
$$\begin{aligned}
	K_{n,3}
	=U(x)-\q_n\ln v_n^{(1)}(x_n^{(1)})+O(\f{1+|x|^\tau}{p_n^{1-\delta}}),\;\text{uniformly in}~t\in[0,1],~x\in B_{\f{r}{\mu_n^{(1)}}}(0).
\end{aligned}$$
Inserting these estimate into \eqref{c6-3}, we obtain
	$$\wt E(t,x_n^{(1)}+\mu_n^{(1)}x)=v_n^{(1)}(x_n^{(1)})\sbr{ 1+\f{1}{p_n}(U(x)-\q_n\ln v_n^{(1)}(x_n^{(1)}))+O(\f{1+|x|^\tau}{p_n^{2-\delta}}) },$$
uniformly in $t\in[0,1],~x\in B_{\f{r}{\mu_n^{(1)}}}(0)~$. It follows that
\be\lab{tem-91} \begin{aligned}
	&\quad (\mu_n^{(1)})^2\wt D_n(x_n^{(1)}+\mu_n^{(1)}x)\\
	&=(\mu_n^{(1)})^2q_n\int_0^1 \wt E(t,x_n^{(1)}+\mu_n^{(1)}x)^{q_n-1}\rd t\\
	&=\f{q_n}{p_n}(v_n^{(1)}(x_n^{(1)}))^{\q_n} \sbr{ 1+\f{1}{p_n}(U(x)-\q_n\ln v_n^{(1)}(x_n^{(1)}))+O(\f{1+|x|^\tau}{p_n^{2-\delta}}) }^{q_n-1},
\end{aligned} \ee
uniformly in $x\in B_{\f{r}{\mu_n^{(1)}}}(0)~$. In particular,
\be\label{c6-5}(\mu_n^{(1)})^2\wt D_n(x_n^{(1)}+\mu_n^{(1)}x)\to e^{U(x)}\quad\text{in }\mathcal{C}_{loc}(\R^2).\ee

Thanks to \eqref{c6-4} and \eqref{c6-5}, we can apply the standard elliptic estimates to the system \eqref{tem-71}, and obtain that up to a subsequence, $(\wt\zeta_n,\wt\kappa_n)\to(\wt\zeta_\iy,\wt\kappa_\iy)$ in $\CR_{\loc}^2(\R^2)$, where $(\wt\zeta_\iy,\wt\kappa_\iy)$ solves the linearized system \eqref{linearsys-0}. Then by Lemma \ref{linear} and $\|\wt\zeta_\iy\|_{L^\infty}, \|\wt \kappa_\iy\|_{L^\infty}\leq 1$,  it holds
	$$(\wt\zeta_\iy,\wt\kappa_\iy)=\sum_{j=0}^2b_j(\phi_j,\phi_j),$$
for some constants $b_j\in\R$.
\ep

Like Section 5, we want to prove that $b_j=0$ for all $j$ and obtain a contradiction. To this goal, we need the following estimates.

\bl\lab{lemma-out-1}
	For fixed $r\in(0,d_0)$, it holds
	\be\label{c6-11} \bcs \zeta_n(x)=B_{n,0}G(x_n^{(1)},x)+\sum_{j=1}^2B_{n,j}\pa_jG(x_n^{(1)},x)+o(\mu_n^{(1)}),\\
		\kappa_n(x)=\wt B_{n,0}G(x_n^{(1)},x)+\sum_{j=1}^2\wt B_{n,j}\pa_jG(x_n^{(1)},x)+o(\mu_n^{(1)}),\\
	\ecs \ee
	in $\CR^1(\Omega\setminus B_{2r}(x_n^{(1)}))$. Moreover, we have
		$$B_{n,0},~\wt B_{n,0}=\f{8\pi}{p_n}b_0+o(\f{1}{p_n}),$$
		\be\label{c6-15}B_{n,j},~\wt B_{n,j}=2\pi b_j\mu_n^{(1)}+o(\mu_n^{(1)}),\quad\text{for }j=1,2.\ee
\el
\bp
The proof is similar to that of Lemma \ref{lemma-out}.
First, by $u_n^{(l)}, v_n^{(l)}=O(\frac{1}{p_n})$ in $\Omega\setminus B_{r}(x_n^{(1)})$, we see from \eqref{c6-7}-\eqref{c6-8} that
\be\label{c6-12}D_n=O(\frac{C^{p_n-1}}{p_n^{p_n-2}}),\quad \wt D_n=O(\frac{C^{q_n-1}}{p_n^{q_n-2}})\quad\text{in $\Omega\setminus B_{r}(x_n^{(1)})$}.\ee
Second,
by Remark \ref{decay-10} and Lemma \ref{tem-70}, we see that for any $x\in\Omega_n^{(1)}$, large $n>1$ and $l=1,2$,
	$$\begin{aligned}
	u_n^{(l)}(x_n^{(1)}+\mu_n^{(1)}x)
	&=v_n^{(l)}(x_n^{(l)})\sbr{ 1+\f{1}{p_n}w_n^{(l)}\sbr{ \f{\mu_n^{(1)}}{\mu_n^{(l)}}x+\f{x_n^{(1)}-x_n^{(l)}}{\mu_n^{(l)}}} }\\
	&\le v_n^{(l)}(x_n^{(l)})\sbr{\f{C}{1+\abs{ \f{\mu_n^{(1)}}{\mu_n^{(l)}}x+\f{x_n^{(1)}-x_n^{(l)}}{\mu_n^{(l)}}}^{4-\delta}}}^{\f{1}{q_n}}\\
	&\le v_n^{(1)}(x_n^{(1)})\sbr{1+O(\f{1}{p_n^{2-\delta}})}\sbr{\f{C}{1+|x|^{4-\delta}}}^{\f{1}{q_n}},
	\end{aligned}$$
and
	$$\begin{aligned}
	v_n^{(l)}(x_n^{(1)}+\mu_n^{(1)}x)
	&=v_n^{(l)}(x_n^{(l)})\sbr{ 1+\f{1}{p_n}z_n^{(l)}\sbr{ \f{\mu_n^{(1)}}{\mu_n^{(l)}}x+\f{x_n^{(1)}-x_n^{(l)}}{\mu_n^{(l)}}} }\\
	&\le v_n^{(1)}(x_n^{(1)})\sbr{1+O(\f{1}{p_n^{2-\delta}})}\sbr{\f{C}{1+|x|^{4-\delta}}}^{\f{1}{p_n}}.
	\end{aligned}$$
Inserting these estimates into the expressions \eqref{c6-7}-\eqref{c6-8} of $D_n,\wt D_n$ and using $(\mu_n^{(1)})^2p_n(v_n^{(1)}(x_n^{(1)}))^{p_n-1}=1$, we obtain that for $x\in\Omega_n^{(1)}$,
	\be\label{c6-9}\abs{(\mu_n^{(1)})^2D_n(x_n^{(1)}+\mu_n^{(1)}x)}\le \f{C}{1+|x|^{4-\delta}},\ee
	\be\label{c6-10}\abs{(\mu_n^{(1)})^2\wt D_n(x_n^{(1)}+\mu_n^{(1)}x)}\le \f{C}{1+|x|^{4-\delta}}.\ee
Thanks to \eqref{c6-12} and \eqref{c6-9}-\eqref{c6-10}, we can follow the proof of the first part in Lemma \ref{lemma-out} and obtain \eqref{c6-11}, with the corresponding
\be\label{c6-13}
	B_{n,0}:=\int_{B_r(x_n^{(1)})}D_n\kappa_n\rd y, \quad  B_{n,j}:=\int_{B_r(x_n^{(1)})}(y-x_n)_jD_n\kappa_n\rd y,
\ee
\be\label{c6-14}
	\wt B_{n,0}:=\int_{B_r(x_n^{(1)})}\wt D_n\zeta_n\rd y, \quad  \wt B_{n,j}:=\int_{B_r(x_n^{(1)})}(y-x_n)_j\wt D_n\zeta_n\rd y.
\ee

Now we compute the estimates of $B_{n,j}$ and $\wt B_{n,j}$ for $j=1,2$. Similarly as the proof of \eqref{c5-91}, by \eqref{c6-4}, \eqref{c6-5} and \eqref{c6-9}-\eqref{c6-10}, we can apply the Dominated Convergence Theorem to obtain
$$\begin{aligned}
	B_{n,j}
	&=\mu_n^{(1)}\int_{B_{\f{r}{\mu_n^{(1)}}}(0)} y_j (\mu_n^{(1)})^2D_n(x_n^{(1)}+\mu_n^{(1)}y)\wt\kappa_n(y)\rd y\\
	&=\mu_n^{(1)}\int_{\R^2}y_je^{U(y)}\wt\kappa_\iy(y)\rd y +o(\mu_n)\\
	&=\mu_n^{(1)}\sum_{i=0}^2b_i\int_{\R^2}y_je^{U(y)}\phi_i(y)\rd y +o(\mu_n)\\
	&=2\pi b_j\mu_n^{(1)}+o(\mu_n),\quad j=1,2,
\end{aligned}$$
and
$$\begin{aligned}
	\wt B_{n,j}
	&=\f{q_n}{p_n}\mu_n^{(1)}\int_{B_{\f{r}{\mu_n^{(1)}}}(0)} y_j(\mu_n^{(1)})^2\wt D_n(x_n^{(1)}+\mu_n^{(1)}y)\wt\zeta_n(y)\rd y\\
	&=\mu_n^{(1)}\int_{\R^2} y_je^{U(y)}\wt\zeta_\iy(y)\rd y +o(\mu_n)\\
	&=2\pi b_j\mu_n^{(1)}+o(\mu_n),\quad j=1,2.
\end{aligned}$$

Finally, we compute the estimates of $B_{n,0}$ and $\wt B_{n,0}$. Observe from \eqref{c6-13}-\eqref{c6-14} that
\be\lab{tem-40-1}\begin{aligned}
	B_{n,0}
	&=\int_{B_r(x_n^{(1)})}\sbr{1-\f{v_n^{(1)}}{v_n^{(1)}(x_n^{(1)})}}D_n\kappa_n\rd x +\f{1}{v_n^{(1)}(x_n^{(1)})}\int_{B_r(x_n^{(1)})} v_n^{(1)}D_n\kappa_n\rd x,
\end{aligned}\ee
and
\be\lab{tem-41-1}\begin{aligned}
	\wt B_{n,0}
	&=\int_{B_r(x_n^{(1)})}\sbr{1-\f{u_n^{(1)}}{v_n^{(1)}(x_n^{(1)})}}\wt D_n\zeta_n\rd x +\f{1}{v_n^{(1)}(x_n^{(1)})}\int_{B_r(x_n^{(1)})} u_n^{(1)}\wt D_n\zeta_n\rd x.
\end{aligned}\ee
Similarly as the proof of \eqref{tem-42}-\eqref{tem-43}, again by the Dominated Convergence Theorem we see that
\be\lab{tem-42-1}\begin{aligned}
	&\quad \int_{B_r(x_n^{(1)})}\sbr{1-\f{v_n^{(1)}}{v_n^{(1)}(x_n^{(1)})}}D_n\kappa_n\rd x \\
	&=-\f{1}{p_n}\int_{B_{\f{r}{\mu_n^{(1)}}}(0)}z_n^{(1)}(y)(\mu_n^{(1)})^2D_n(x_n^{(1)}+\mu_n^{(1)}y)\wt\kappa_n(y)\rd y\\
	&=-\f{1}{p_n}\int_{\R^2}U(y)e^{U(y)}\wt\kappa_\iy(y)\rd y+o(\f{1}{p_n})\\
	&=-\f{1}{p_n}\sum_{i=0}^2b_i\int_{\R^2}U(y)e^{U(y)}\phi_i(y)\rd y +o(\f{1}{p_n})\\
	&=\f{8\pi}{p_n}b_0+o(\f{1}{p_n}),
\end{aligned}\ee
and
\be\lab{tem-43-1}\begin{aligned}
	&\quad \int_{B_r(x_n^{(1)})}\sbr{1-\f{u_n^{(1)}}{v_n^{(1)}(x_n^{(1)})}}\wt D_n\zeta_n\rd x \\
	&=-\f{1}{p_n}\int_{B_{\f{r}{\mu_n^{(1)}}}(0)}w_n^{(1)}(y)(\mu_n^{(1)})^2\wt D_n(x_n^{(1)}+\mu_n^{(1)}y)\wt\zeta_n(y)\rd y\\
	&=-\f{1}{p_n}\int_{\R^2}(U(y)-\f{\q}{2})e^{U(y)}\wt\zeta_\iy(y)\rd y+o(\f{1}{p_n})\\
	&=-\f{1}{p_n}\sum_{i=0}^2b_i\int_{\R^2}(U(y)-\f{\q}{2})e^{U(y)}\phi_i(y)\rd y +o(\f{1}{p_n})\\
	&=\f{8\pi}{p_n}b_0+o(\f{1}{p_n}).
\end{aligned}\ee
Substituting \eqref{tem-42-1}-\eqref{tem-43-1} into \eqref{tem-40-1}-\eqref{tem-41-1}, we get
\be\lab{tem-44-1} B_{n,0}=\f{8\pi}{p_n}b_0+o(\f{1}{p_n})+\f{1}{v_n^{(1)}(x_n^{(1)})}\int_{B_r(x_n^{(1)})} v_n^{(1)}D_n\kappa_n\rd x,
\ee
\be\lab{tem-45-1}
	\wt B_{n,0}=\f{8\pi}{p_n}b_0+o(\f{1}{p_n})+\f{1}{v_n^{(1)}(x_n^{(1)})}\int_{B_r(x_n^{(1)})} u_n^{(1)}\wt D_n\zeta_n\rd x.
\ee
So it suffices to estimate $\int_{B_r(x_n^{(1)})} v_n^{(1)}D_n\kappa_n\rd x$ and $\int_{B_r(x_n^{(1)})} u_n^{(1)}\wt D_n\zeta_n\rd x$.

Since Lemma \ref{decay-12} and \eqref{c6-9}-\eqref{c6-10} imply
	{\allowdisplaybreaks
\begin{align*}
	&\abs{\f{1}{v_n^{(1)}(x_n^{(1)})}\int_{B_r(x_n^{(1)})} v_n^{(1)}D_n\kappa_n\rd x}\\
=&\int_{B_{\f{r}{\mu_n^{(1)}}}(0)}\abs{\sbr{1+\f{z_n^{(1)}(y)}{p_n}}(\mu_n^{(1)})^2 D_n(x_n^{(1)}+\mu_n^{(1)}y) \wt\kappa_n(y) }\rd y\\
\le& C\int_{B_{\f{r}{\mu_n^{(1)}}}(0)}\f{\ln(2+|y|)}{1+|y|^{4-\delta}}\rd y=O(1),
	\end{align*}
}%
and
	$$\begin{aligned}
	&\abs{\f{1}{v_n^{(1)}(x_n^{(1)})}\int_{B_r(x_n^{(1)})} u_n^{(1)}\wt D_n\zeta_n\rd x}\\
=&\int_{B_{\f{r}{\mu_n^{(1)}}}(0)}\abs{\sbr{1+\f{w_n^{(1)}(y)}{p_n}}(\mu_n^{(1)})^2 \wt D_n(x_n^{(1)}+\mu_n^{(1)}y) \wt\zeta_n(y) }\rd y\\
\le& C\int_{B_{\f{r}{\mu_n^{(1)}}}(0)}\f{\ln(2+|y|)}{1+|y|^{4-\delta}}\rd y=O(1),
	\end{aligned}$$
we have that $B_{n,0},\wt B_{n,0}=O(1)$. On the other hand, similarly as \eqref{tem-74} we obtain
\be\label{c6-18}\begin{aligned}
	&\int_{B_{2r}(x_n^{(1)})} v_n^{(1)}D_n\kappa_n\rd x=\int_{B_{2r}(x_n^{(1)})} -\Delta\zeta_nv_n^{(1)}\rd x\\
	&=\int_{B_{2r}(x_n^{(1)})} \Delta v_n^{(1)}\zeta_n-\Delta\zeta_nv_n^{(1)}\rd x+\int_{B_{2r}(x_n^{(1)})} (u_n^{(1)})^{q_n}\zeta_n\rd x\\
	&=\int_{\pa B_{2r}(x_n^{(1)})}\abr{\nabla v_n^{(1)},\nu}\zeta_n-\abr{\nabla\zeta_n,\nu}v_n^{(1)}\rd\sigma_x+\int_{B_{2r}(x_n^{(1)})} (u_n^{(1)})^{q_n}\zeta_n\rd x,\\
\end{aligned}\ee
Using Lemma \ref{expansion-3}, \eqref{c6-11} and \eqref{c6-15}, we see that
	$$\int_{\pa B_{2r}(x_n^{(1)})}\abr{\nabla v_n^{(1)},\nu}\zeta_n-\abr{\nabla\zeta_n,\nu}v_n^{(1)}\rd\sigma_x=O\sbr{\f{\wt B_{n,1}+\wt B_{n,2}}{p_n}+\f{\mu_n\wt B_{n,0}}{p_n^{2-\delta}}}=O(\f{\mu_n}{p_n}).$$
We also observe from \eqref{c5-93} that
	$$\begin{aligned}
	\int_{B_{2r}(x_n^{(1)})} (u_n^{(1)})^{q_n}\zeta_n\rd x	&=\f{1}{p_n}v_n^{(1)}(x_n^{(1)})^{\q_n+1}\int_{B_{\f{2r}{\mu_n^{(1)}}}(0)}\sbr{1+\f{w_n^{(1)}(y)}{p_n} }^{q_n}\wt\zeta_n(y)\rd y\\	&=\f{1}{p_n}v_n^{(1)}(x_n^{(1)})^{\q_n+1}\int_{\R^2}e^{U-\frac{\theta}{2}}\sum_{j=0}^2b_j\phi_j\rd y+o(\f{1}{p_n})\\
	&=0+o(\f{1}{p_n})=o(\f{1}{p_n}).
	\end{aligned}$$
Inserting these estimates into \eqref{c6-18} we get
	$$\int_{B_{2r}(x_n^{(1)})} v_n^{(1)}D_n\kappa_n\rd x=o(\f{1}{p_n}).$$
From here and $v_n^{(l)}=O(\f{1}{p_n})$, $D_n=O(\f{C^{p_n-1}}{p_n^{p_n-2}})$ in $B_{2r}(x_n^{(1)})\setminus B_{r}(x_n^{(1)})$ (see \eqref{c6-12}), we obtain
\be\lab{tem-77}
	\int_{B_{r}(x_n^{(1)})} v_n^{(1)}D_n\kappa_n\rd x=\int_{B_{2r}(x_n^{(1)})} v_n^{(1)}D_n\kappa_n\rd x+O(\f{C^{p_n-1}}{p_n^{p_n-1}})=o(\f{1}{p_n}).
\ee
By a similar argument one can prove
\be\lab{tem-78}\begin{aligned}
	\int_{B_r(x_n^{(1)})} u_n^{(1)}\wt D_n\zeta_n\rd x=o(\f{1}{p_n}).
\end{aligned}\ee
Inserting \eqref{tem-77}-\eqref{tem-78} into \eqref{tem-44-1}-\eqref{tem-45-1}, we finally obtain
	$$B_{n,0},~\wt B_{n,0}=\f{8\pi}{p_n}b_0+o(\f{1}{p_n}).$$
This completes the proof.
\ep

\bl\lab{out-2-1}
	It holds
	\be
		\nm{\zeta_n}_{L^\iy(\Omega\setminus B_{2\mu_n^{(1)}p_n}(x_n))}=O\sbr{p_nB_{n,0}}+O(\f{1}{p_n^{1-\delta}}),
	\ee
	\be
		\nm{\kappa_n}_{L^\iy(\Omega\setminus B_{2\mu_n^{(1)}p_n}(x_n))}=O\sbr{p_n\wt B_{n,0}}+O(\f{1}{p_n^{1-\delta}}).
	\ee
\el

\begin{proof}
This lemma can be proved exactly as Lemma \ref{out-2}.
\end{proof}

\bl\lab{lemma-in-1}
	It holds $b_0=b_1=b_2=0$, where $b_j$ are constants in Lemma \ref{tem-76}, i.e. $\wt \xi_\infty=\wt \kappa_\iy=0$.
\el
\bp
For fixed $r\in(0,d_0)$, we recall the quadratic forms (i.e. replace $x_n$ with $x_n^{(1)}$ in \eqref{function-P}-\eqref{function-Q})
\be\label{c6-37}
	P(u,v)=-2r\int_{\pa B_r(x_n^{(1)})} \abr{\nabla u,\nu}\abr{\nabla v,\nu}\rd\sigma+r\int_{\pa B_r(x_n^{(1)})}\abr{\nabla u,\nabla v}\rd\sigma,
\ee
\be
	Q_i(u,v)=-\int_{\pa B_r(x_n^{(1)})} \abr{\nabla u,\nu}\pa_i v+\abr{\nabla v,\nu}\pa_i u\rd \sigma+\int_{\pa B_r(x_n^{(1)})}\abr{\nabla u,\nabla v}\nu_i\rd\sigma.
\ee
Then \eqref{formula-P-1} implies
\be\lab{tem-81}\begin{aligned}
		P(u_n^{(l)},v_n^{(l)})
		&=r\int_{\pa B_r(x_n^{(1)})} \f{(u_n^{(l)})^{q_n+1}}{q_n+1}+\f{(v_n^{(l)})^{p_n+1}}{p_n+1}\rd\sigma-2\int_{B_r(x_n^{(1)})}\f{(u_n^{(l)})^{q_n+1}}{q_n+1}+\f{(v_n^{(l)})^{p_n+1}}{p_n+1}\rd x.
\end{aligned}\ee
By \eqref{constant-Ln}-\eqref{function-zeta} and \eqref{c6-37}, we have
\be\label{c6-22}\f{P(u_n^{(1)},v_n^{(1)})}{L_n}-\f{P(u_n^{(2)},v_n^{(2)})}{L_n}=P(v_n^{(1)},\zeta_n)+P(u_n^{(2)},\kappa_n).\ee
By $u_n^{(l)},~v_n^{(l)}=O(\f{1}{p_n})$ and $D_n=O(\f{C^{p_n-1}}{p_n^{p_n-2}})$, $\wt D_n=O(\f{C^{q_n-1}}{p_n^{q_n-2}})$ in $\pa B_{r}(x_n^{(1)})$ (see \eqref{c6-12}), we have
	\be\label{c6-24}\int_{\pa B_r(x_n^{(1)})}\f{(u_n^{(1)})^{q_n+1}-(u_n^{(2)})^{q_n+1}}{L_n}\rd x=\int_{\pa B_r(x_n^{(1)})}(u_n^{(1)})^{q_n}\zeta_n +u_n^{(2)}\wt D_n\zeta_n\rd x=O(\f{C^{q_n-1}}{p_n^{q_n-1}}),\ee
	\be\label{c6-25}\int_{\pa B_r(x_n^{(1)})}\f{(v_n^{(1)})^{p_n+1}-(v_n^{(2)})^{p_n+1}}{L_n}\rd x=\int_{\pa B_r(x_n^{(1)})}(v_n^{(1)})^{p_n}\kappa_n +v_n^{(2)} D_n\kappa_n\rd x=O(\f{C^{p_n-1}}{p_n^{p_n-1}}).\ee
Besides, by similar arguments as \eqref{tem-77}-\eqref{tem-78} we obtain
\[
	\int_{B_r(x_n^{(1)})}\f{(u_n^{(1)})^{q_n+1}-(u_n^{(2)})^{q_n+1}}{L_n}\rd x=\int_{B_r(x_n^{(1)})}(u_n^{(1)})^{q_n}\zeta_n +u_n^{(2)}\wt D_n\zeta_n\rd x=o(\f{1}{p_n}),
\]
\[
	\int_{B_r(x_n^{(1)})}\f{(v_n^{(1)})^{p_n+1}-(v_n^{(2)})^{p_n+1}}{L_n}\rd x=\int_{B_r(x_n^{(1)})}(v_n^{(1)})^{p_n}\kappa_n +v_n^{(2)}D_n\kappa_n\rd x=o(\f{1}{p_n}).
\]
Thus \eqref{tem-81}-\eqref{c6-22} imply
	$$P(v_n^{(1)},\zeta_n)+P(u_n^{(2)},\kappa_n)=o(\f{1}{p_n^2}).$$
On the other hand, using the expansions in Lemma \ref{expansion-3} and Lemma \ref{lemma-out-1}, we get
$$\begin{aligned}
	P(v_n^{(1)},\zeta_n)
	&=\wt C_n^{(1)}B_{n,0}P(G(x_n^{(1)},\cdot),G(x_n^{(1)},\cdot))+o(\f{\mu_n}{p_n})\\
	&=-\f{1}{2\pi}\wt C_n^{(1)}B_{n,0}+o(\f{\mu_n}{p_n})=-\f{32\pi\sqrt e}{p_n^2}b_0+o(\f{1}{p_n^2}),
\end{aligned}$$
and
$$\begin{aligned}
	P(u_n^{(2)},\kappa_n)
	&=-\f{1}{2\pi}C_n^{(2)}\wt B_{n,0}+o(\f{\mu_n}{p_n})=-\f{32\pi\sqrt e}{p_n^2}b_0+o(\f{1}{p_n^2}).
\end{aligned}$$
It follows that
	$$\f{64\pi\sqrt e}{p_n^2}b_0+o(\f{1}{p_n^2})=o(\f{1}{p_n^2}),$$
which implies $b_0=0$.

Recall \eqref{formula-Q-1} that
\be
	Q_i(u_n^{(l)},v_n^{(l)})=\int_{\pa B_r(x_n^{(1)})} \sbr{ \f{(u_n^{(l)})^{q_n+1}}{q_n+1}+\f{(v_n^{(l)})^{p_n+1}}{p_n+1} }\nu_i\rd\sigma,
\ee
which, together with \eqref{c6-24}-\eqref{c6-25}, gives
\be\label{c6-26}\begin{aligned}
	&\quad Q_i(v_n^{(1)},\zeta_n)+Q_i(u_n^{(2)},\kappa_n)=\f{Q_i(u_n^{(1)},v_n^{(1)})}{L_n}-\f{Q_i(u_n^{(2)},v_n^{(2)})}{L_n}\\
&=\int_{\pa B_r(x_n^{(1)})}\f{(u_n^{(1)})^{q_n+1}-(u_n^{(2)})^{q_n+1}}{L_n (q_n+1)}\rd x	+\int_{\pa B_r(x_n^{(1)})}\f{(v_n^{(1)})^{p_n+1}-(v_n^{(2)})^{p_n+1}}{L_n (p_n+1)}\rd x\\
	&=O(\f{C^{p_n-1}}{p_n^{p_n}}).
\end{aligned}\ee
Again by using Lemma \ref{expansion-3} and Lemma \ref{lemma-out-1}, we get
$$\begin{aligned}
	Q_i(v_n^{(1)},\zeta_n)
	&=\wt C_n^{(1)}B_{n,0}Q_i(G(x_n^{(1)},\cdot),G(x_n^{(1)},\cdot))+\sum_{j=1}^2\wt C_n^{(1)}B_{n,j}Q_i(G(x_n^{(1)},\cdot),\pa_jG(x_n^{(1)},\cdot))+o(\f{\mu_n}{p_n})\\
	&=-\wt C_n^{(1)}B_{n,0}\pa_iR(x_n^{(1)})-\f{1}{2}\sum_{j=1}^2\wt C_n^{(1)}B_{n,j}\pa_{ij}^2R(x_n^{(1)})+o(\f{\mu_n}{p_n})\\
	&=-\f{1}{2}\sum_{j=1}^2\wt C_n^{(1)}B_{n,j}\pa_{ij}^2R(x_n^{(1)})+o(\f{\mu_n}{p_n}),
\end{aligned}$$
and similarly
$$\begin{aligned}
	Q_i(u_n^{(2)},\kappa_n)
	&=-\f{1}{2}\sum_{j=1}^2 C_n^{(2)}\wt B_{n,j}\pa_{ij}^2R(x_n^{(1)})+o(\f{\mu_n}{p_n}),
\end{aligned}$$
Inserting these into \eqref{c6-26} leads to
	$$\nabla^2R(x_n^{(1)})\cdot\sbr{ C_n^{(2)}\wt B_{n,1}+\wt C_n^{(1)} B_{n,1}, C_n^{(2)}\wt B_{n,2}+\wt C_n^{(1)} B_{n,2} }^T=o(\f{\mu_n}{p_n}).$$
Since $\nabla^2R(x_n)\to\nabla^2R(x_\iy)$ and $\nabla^2R(x_\iy)$ is invertible, we know that for $n$ large,
	$$C_n^{(2)}\wt B_{n,j}+\wt C_n^{(1)} B_{n,j}=o(\f{\mu_n}{p_n}),\quad j=1,2,$$
which implies $b_1=b_2=0$ by using \eqref{c6-31}-\eqref{c6-32} and \eqref{c6-15}.
\ep

\bc\lab{in-1}
	For any $d>0$, it holds
	$$\nm{\zeta_n}_{L^\iy(B_{\mu_n^{(1)}d}(x_n^{(1)}))}=o(1),\quad \nm{\kappa_n}_{L^\iy(B_{\mu_n^{(1)}d}(x_n^{(1)}))}=o(1).$$
\ec

\bp
Thanks to Lemma \ref{lemma-in-1}, the proof is the same as Corollary \ref{in}.
\ep

Now we are in the position to prove Theorem \ref{thm3}.
\bp[Proof of Theorem \ref{thm3}] Recall that
by contradiction, we have supposed that $(u_n^{(1)},v_n^{(1)})$ and $(u_n^{(2)},v_n^{(2)})$ are two sequences of different $1$-bubble solutions concentrating at the same $x_\iy$. Let $\zeta_n,\kappa_n$ be defined by \eqref{function-zeta}. Without loss of generality, we assume
	$$\zeta_n(x_n^*)=1=\max_\Omega\lbr{|\zeta_n|,|\kappa_n|}.$$
For any $d>0$, since $B_{n,0},\wt B_{n,0}=o(\f{1}{p_n})$, Lemma \ref{out-2-1} and Corollary \ref{in-1} imply $x_n^*\in B_{2\mu_n^{(1)}p_n}(x_n^{(1)})\setminus B_{\mu_n^{(1)}d}(x_n^{(1)})$ for $n$ large. Let
\be\lab{constant-rn-1}
	r_n:=|x_n^*-x_n^{(1)}|\in [\mu_n^{(1)}d, 2\mu_n^{(1)}p_n].
\ee
Recalling $\wt \zeta_n,\wt \kappa_n$ defined in \eqref{function-wtzeta}, 
	$$\zeta_n^{**}(x):=\zeta_n(x_n^{(1)}+r_nx)=\wt \zeta_n(\frac{r_n}{\mu_n^{(1)}}x),\quad \kappa_n^{**}(x):=\kappa_n(x_n^{(1)}+r_nx)=\wt \kappa_n(\frac{r_n}{\mu_n^{(1)}}x).$$
Then a similar argument as in the proof of Theorem \ref{thm2} implies $\zeta_n^{**}\to1$ in $\CR_{loc}^2(\R^2\setminus\{0\})$. Thus $\inf_{|x|=1}\zeta_n^{**}(x)\ge C_0>0$ for $n$ large. Define the average functions
\be
	\zeta_n^*(r):=\f{1}{2\pi}\int_0^{2\pi}\wt\zeta_n(r,\rho)\rd\rho,\quad \kappa_n^*(r):=\f{1}{2\pi}\int_0^{2\pi}\wt\kappa_n(r,\rho)\rd\rho.
\ee
Then using \eqref{tem-90}, \eqref{tem-91} and \eqref{c6-9}-\eqref{c6-10}, we can prove exactly as \eqref{tem-6000} that
	$$\zeta_n^*(\f{r_n}{\mu_n^{(1)}})=o(1),\quad \kappa_n^*(\f{r_n}{\mu_n^{(1)}})=o(1).$$
Consequently,	$$o(1)=\zeta_n^*(\f{r_n}{\mu_n^{(1)}})=\f{1}{2\pi}\int_0^{2\pi}\wt\zeta_n(\f{r_n}{\mu_n^{(1)}},\rho)\rd\rho=\f{1}{2\pi}\int_0^{2\pi}\zeta_n^{**}(1,\rho)\rd\rho\ge C_0>0,$$
a contradiction. This proves Theorem \ref{thm3}.
\ep

\section{On convex domains}
In this section, we assume $\Omega$ is convex and prove Theorem \ref{thm4}. We first colloect some known results.
\bl\lab{tem-100}(\cite{priori-1})
	Let $\Omega\subset\R^2$ be any star-shaped bounded domain with $\CR^2$-boundary and assume the exponents $p,q$ satisfies
		$$p\ge1,\quad q\ge1,\quad pq-1\ge\kappa,\quad \f{1}{K}q\le p\le Kq,$$
	for some $\kappa>0$ and $K\ge1$. There exists a positive constant $C=C(\kappa,K,\Omega)$ such that any classical solution $(u,v)=(u_{p,q},v_{p,q})$ of \eqref{equ1} satisfies
		$$p\int_\Omega\nabla u_{p,q}\cdot\nabla v_{p,q}\rd x\le C.$$
\el

\bl\lab{tem-97}(\cite{cp-3})
	Let $\Omega\subset\R^2$ be any convex bounded domain, then the Robin function $R(x)$ is strictly convex and has a unique critical point. Moreover the corresponding Hessian matrix $\nabla^2R$ in this point is positive definite.
\el

\bl\lab{tem-98}(\cite{cp-1})
	Let $\Omega\subset\R^2$ be any convex bounded domain and $k\ge2$. Set
		$$\Omega_0^k=\lbr{(x_1,\cdots,x_k)\in\Omega^k:~x_i=x_j~\text{for some}~i\neq j}.$$
	Then there does not exist $(x_1,\cdots,x_k)\in\Omega^k\setminus\Omega_0^k$ such that
		$$\nabla R(x_i)-2\sum_{j=1,j\neq i}^k\nabla G(x_i,x_j)=0,\quad \forall i=1,\cdots,k.$$
\el

\vskip0.2in

\bp[Proof of Theorem \ref{thm4}]
The existence of positive solutions has been proved in \cite{exist1} by the degree method, and later proved in \cite{unique-2} by the variational method. 

Now we prove the uniqueness. Let $\Omega$ be convex and fix any $\Lambda>0$. By contradiction, we assume there exist two number sequences $p_n,q_n\to\iy$ satisfying $\sup_{n\ge1}|p_n-q_n|\le \Lambda$ such that the corresponding Lane-Emden system \eqref{tem-95} with $(p,q)=(p_n, q_n)$ has two different solutions $(u_n^{(1)},v_n^{(1)})$ and $(u_n^{(2)},v_n^{(2)})$, i.e.
\be\lab{tem-99}
	(u_n^{(1)},~v_n^{(1)})\neq (u_n^{(2)},~v_n^{(2)}),\quad\text{for any}~n\ge1.
\ee
From Lemma \ref{tem-100}, we see that
	$$\sup_{n\ge1}p_n\int_\Omega\nabla u_n^{(l)}\cdot\nabla v_n^{(l)}\rd x<+\iy,$$
for $l=1,2$. Applying Theorem \ref{thm0}, we deduce that there exist two integers $k_l\ge1$, two sets $\SR_l\subset\Omega$ of $k_l$ concerntration points such that conclusions (1)-(5) in Theorem \ref{thm0} hold for $(u_n^{(l)},v_n^{(l)})$, $l=1,2$.

Since $\Omega$ is convex, we can apply Lemma \ref{tem-98} to deduce that $k_1=k_2=1$. So the concentration sets $\SR_l=\{x_{\iy}^{(l)}\}$ and $\nabla R(x_\iy^{(l)})=0$ for $l=1,2$. Then Lemma \ref{tem-97} implies that $x_\iy^{(1)}=x_\iy^{(2)}$ is the unique critical point of the Robin function $R(x)$, which is also nondegenerate.
So by Theorem \ref{thm3}, there exists $n_0>1$ such that
	$$u_n^{(1)}=u_n^{(2)}\quad\text{and}\quad v_n^{(1)}=v_n^{(2)},\quad\text{for any}~n\ge n_0.$$
This clearly contradicts with \eqref{tem-99}. This proves the uniqueness, and the nondegeneracy follows from Theorem \ref{thm2}.
\ep

\vskip0.26in
\begin{appendices}

\end{appendices}

\vskip 0.1in\noindent
{\bf Acknowledgements} \quad This work is partially supported by NSFC (No.12171265,12071240).


\vskip0.2in
{\sc Address of the authors:}

\vskip0.15in
\indent Zhijie Chen\\
\indent Yau Mathematical Sciences Center\\
\indent Tsinghua University\\
\indent Beijing 100084\\
\indent China\\
\indent {\it E-mail address: zjchen2016@tsinghua.edu.cn}

\vskip0.15in
\indent Houwang Li\\
\indent Yanqi Lake Beijing Institute of Mathematical Sciences and Applications\\
\indent Beijing 101408\\
\indent China\\
\indent {\it E-mail address: lhwmath@bimsa.cn}

\vskip0.15in
\indent Wenming Zou\\
\indent Department of Mathematical Sciences\\
\indent Tsinghua University\\
\indent Beijing 100084\\
\indent China\\
\indent {\it E-mail address: zou-wm@mail.tsinghua.edu.cn}


\end{document}